\documentclass[11pt,french]{amsart}
\begin{document}

%ceci est la version finale de l'article intitule : 
%Correspondance de Jacquet-Langlands pour les corps locaux de 
%caracteristique non nulle.
%16.12.2001

\newtheorem{theo}{Th\'eor\`eme}[section]
\newtheorem{prop}[theo]{Proposition}
\newtheorem{lemme}[theo]{Lemme}
\newtheorem{cor}[theo]{Corollaire}
\newtheorem{rem}[theo]{Remarque}
\numberwithin{equation}{section}

\def\cf{{\bf C}_L}
\def\r{{\mathbb R}}
\def\ppl{{\bf P}}

\def\ssr{\'el\'ement semisimple r\'egulier}
\def\ssrs{\'el\'ements semisimples r\'eguliers}

\def\sgls{sous-groupe de Levi standard}
\def\sglss{sous-groupes de Levi standard}

\def\eci{repr\'esentation essentiellement de carr\'e int\'egrable}
\def\ecis{repr\'esentations essentiellement de carr\'e int\'egrable}

\def\ett{repr\'esentation essentiellement temp\'er\'ee}
\def\etts{repr\'esentations essentiellement temp\'er\'ees}

\def\rir{repr\'esentation irr\'eductible}
\def\rirs{repr\'esentations irr\'eductibles}

\def\rli{repr\'esentation lisse irr\'eductible}
\def\rlis{repr\'esentations lisses irr\'eductibles}

\def\cusp{repr\'esentation cuspidale}
\def\cusps{repr\'esentations cuspidales}

\def\care{repr\'esentation de carr\'e int\'egrable}
\def\cares{repr\'esentations de carr\'e int\'egrable}

\def\temp{repr\'esentation temp\'er\'ee}
\def\temps{repr\'esentations temp\'er\'ees}

\def\aa{\mathcal A}
\def\a{\mathbb A}
\def\o{\mathcal O}
\def\F{{\mathbb F}}
\def\G{\tilde{G}}
\def\g{{\mathbb G}}
\def\D{\mathbb D}
\def\lra{\leftrightarrow}

\def\zzzf{\zeta_{FL}^m}
\def\zzz{\bar{\zeta}^m_{FL}}
\def\zzd{\bar{\zeta}_{D_FD_L}^m}
\def\zzzd{\zeta_{D_FD_L}^m}
\def\zzzz{\bar{\zeta}^{m''}_{FL}}
\def\zzzdl{\zeta_{D_FD_L}^l}
\def\zzdl{\bar{\zeta}_{D_FD_L}^l}

\def\lam{\lambda_{FL}^m}
\def\lll{\bar{\lambda}^m_{FL}}
\def\llam{\lambda_{EK}^m}
\def\llll{\bar{\lambda}_{EK}^m}
\def\llld{\bar{\lambda}_{D_FD_L}^m}
\def\lamd{\lambda_{D_FD_L}^m}
\def\lamdl{\lambda_{D_FD_L}^l}
\def\llldl{\bar{\lambda}_{D_FD_L}^l}

\def\vvf{v_{M_n(F)}}
\def\vvl{v_{M_n(L)}}

\def\T{{\bf{T}}_{F,l}}
\def\TT{\tilde{\bf{T}}_{F,l,A}}
\def\Tl{{\bf{T}}_{L,l}}
\def\TTl{\tilde{\bf{T}}_{L,l,\zzzdl(A)}}

\def\rep{re\-pr\'e\-sen\-ta\-ti\-on}

\def\hfl#1#2{\smash{\mathop{\hbox to 12mm{\rightarrowfill}}
\limits^{\scriptstyle#1}_{\scriptstyle#2}}}
\def\diagram#1{\def\normalbaselines{\baselineskip=0pt
\lineskip=10pt\lineskiplimit=1pt}   \begin{matrix}#1 \end{matrix}}
\def\vfl#1#2{\llap{$\scriptstyle #1$}\left\uparrow
\vbox to 6mm{}\right.\rlap{$\scriptstyle #2$}}

\title[Jacquet-Langlands en caract\'eristique non nulle]{Correspondance
de  Jacquet-Langlands pour les corps locaux de caract\'eristique non nulle}
\maketitle
\centerline{\bf Jacquet-Langlands correspondence in non zero
characteristic}
\centerline{\footnotesize{par Alexandru Ioan
BADULESCU}{\footnote{Alexandru Ioan BADULESCU, Universit\'e de
Poitiers, UFR Sciences SP2MI, D\'epartement de Math\'ematiques,
T\'el\'eport 2, Boulevard Marie et Pierre Curie, BP 30179, 86962
FUTUROSCOPE CHASSENEUIL CEDEX\\
E-mail : badulesc\makeatletter @\makeatother
wallis.sp2mi.univ-poitiers.fr}}}
\ \\
\ \\

{\footnotesize{
\ \\
{\bf Abstract :} In this article we prove the Jacquet-Langlands
local correspondence in non-zero characteristic. Let $F$ be a local
field of non-zero charactersitic and $G'$ an inner form of
$GL_n(F)$; then, following [Ka], we prove relations between the
representation theory of $G'$ and the representation theory of an inner
form of $GL_n(L)$, where $L$ is a local field of zero characteristic
close to  $F$. The proof of the Jacquet-Langlands correspondence between
$G'$ and $GL_n(F)$ is done using the above results and ideas from the
proof by Deligne, Kazhdan and Vign\'eras ([DKV]) of the zero characteristic
case.
We also get the following, already known in zero
characteristic~: orthogonality relations for $G'$, inequality
involving conductor and level for representations of $G'$ and finiteness
for automorphic cuspidal representations with fixed
component at almost every place for an inner form of $GL_n$ over a
global field of non-zero characteristic.   
\ \\
\ \\
{\bf R\'esum\'e :} Le but de cet article est la preuve de la
correspondance de Jacquet-Langlands locale en caract\'eristique non
nulle. Si $F$ est un corps local de caract\'eristique
non nulle et 
$G'$ est une forme int\'erieure de $GL_n(F)$, on construit un
parall\`ele entre la th\'eorie des repr\'esentations de $G'$ et celle 
d'une forme int\'erieure de $GL_n(L)$ respectivement, 
o\`u $L$ est un corps local de caract\'eristique nulle, proche de $F$
au sens de [Ka]. La correspondance de Jacquet-Langlands entre $G'$
et $GL_n(F)$ est prouv\'ee en utilisant cette construction et les
id\'ees d\'evelopp\'ees  d\'ej\`a par Deligne, Kazhdan et Vign\'eras
([DKV]) pour la
preuve en caract\'eristique nulle. Nous obtenons au passage les
r\'esultats suivants, connus jusqu'ici uniquement en
caract\'eristique nulle~: le
th\'eor\`eme d'orthogonalit\'e des caract\`eres pour $G'$, des
relations entre le conducteur et le niveau d'une 
repr\'esentation de $G'$, ainsi qu'un th\'eor\`eme de finitude pour les
repr\'esentations automorphes cuspidales d'une forme int\'erieure de
$GL_n$ sur un corps global de caract\'eristique non nulle.

\tableofcontents}}
\section{Introduction}

Soient $F$ un corps local non archim\'edien et $A$ une alg\`ebre
 centrale simple de dimension finie sur $F$. Soit $A^*$ le groupe des
 \'el\'ements 
 inversibles de $A$. On sait que $A^*$ est isomorphe \`a $GL_r(D)$,
 o\`u $D$ est une alg\`ebre \`a division sur $F$, et on l'identifiera
 avec ce groupe. La dimension de $D$ en tant qu'espace vectoriel sur
 $F$ est un carr\'e $d^2$. On pose $n=rd$, $G=GL_n(F)$ et
 $G'=GL_r(D)$. Si $O_F$ (resp.$O_D$) est l'anneau des entiers de $F$
 (resp.$D$), on fixe une mesure de Haar sur $G$ (resp.$G'$) telle que le
 volume de $GL_n(O_F)$ (resp.$GL_r(O_D)$) soit \'egal \`a 1.
Un \'el\'ement de
 $G$ ou $G'$ dont le polyn\^ome caract\'eristique est
 s\'eparable (i.e. sans racine multiple sur une cl\^oture alg\'ebrique
 de $F$) est dit {\it semisimple r\'egulier}. Si $g$ est
 un \'el\'ement de $G$ et $g'$ un \'el\'ement de $G'$ on \'ecrit
 $g\lra g'$ si $g$ et $g'$ sont semisimples r\'eguliers et ont le
 m\^eme polyn\^ome caract\'eristique. On dit alors que $g$ 
 et $g'$ {\it se correspondent}. Si $\pi$
 est une repr\'esentation lisse de longueur finie de $G$ ou $G'$, 
 on note $\chi_{\pi}$ le caract\`ere de $\pi$.

On note $E^2(G)$ l'ensemble des classes d'\'equivalence de
\ecis\  de $G$ et
$E^2(G')$ l'ensemble des classes d'\'equivalence de \ecis\ de $G'$. La
correspondance de Jacquet-Langlands s'\'enonce de la fa\c{c}on
suivante~:\\ 

\begin{theo}\label{JL} Il existe une unique bijection~:
$${\bf{C}}:E^2(G)\to E^2(G')$$
telle que pour tout $\pi\in E^2(G)$ on ait 
\begin{equation}\label{JLe}
\chi_{\pi}(g)=(-1)^{n-r} \chi_{{\bf C}(\pi)}(g'),
\end{equation}
\`a chaque fois que $g\in G$ et $g'\in G'$ se correspondent.
\end{theo}

Ce r\'esultat a \'et\'e prouv\'e pour $F$ de
caract\'eristique nulle dans [DKV] et le but de cet article est de
donner une d\'emonstration en caract\'eristique non nulle. Comme nous
le verrons plus bas, un des principaux obstacles {\it en caract\'eristique
non nulle} est le fait que les relations d'orthogonalit\'e des
caract\`eres des \cares\ ne sont pas prouv\'ees \`a ce jour pour le
groupe $G'$. Bien entendu, si nous d\'emontrons la correspondance de
Jacquet-Langlands d'une autre fa\c{c}on, ces relations
d'orthogonalit\'e pour $G'$ en d\'ecouleront puisqu'on les a d\'ej\`a
montr\'ees pour $G$ ([Ba2]). C'est en effet ce qui va se
passer. Attaquer la correspondance en essayant de prouver
l'orthogonalit\'e pour $G'$ en caract\'eristique non nulle semble
depuis longtemps moins raisonnable que de l'attaquer par la m\'ethode
des ``corps
proches'' en suivant une id\'ee de Kazhdan ([Ka]). C'est ce chemin
que nous allons suivre. Expliquons la d\'emarche pr\'econis\'ee et les
difficult\'es rencontr\'ees. 

Soient $F$ un corps local de caract\'eristique non nulle et $D_F$ une
alg\`ebre \`a division centrale sur $F$ de dimension $d^2$. Soit $r$
un entier strictement positif. On pose $n=rd$. On note $G'_F$ le
groupe $GL_r(D_F)$. Soient $O_F$ l'anneau des entiers de $F$, $P_F$
l'id\'eal maximal de $O_F$ et $\pi_F$ une uniformisante de $F$. Soient
$O_{D_F}$ l'anneau des entiers de $D_F$ et $P_{D_F}$ l'id\'eal
maximal de $O_{D_F}$. On fixe une uniformisante $\pi_{D_F}$ de $D_F$. On pose $K_F=GL_r(O_{D_F})$ et, pour tout entier
positif non nul $l$, $K_F^l=Id+M_r(P^{ld}_{D_F})$. Si $K$ est un
sous-groupe ouvert 
compact de $G'_F$, on note $H(G'_F;K)$ l'alg\`ebre de Hecke des
fonctions localement constantes \`a support compact sur $G'_F$ qui sont
bi-invariantes par $K$. Si $\pi$ est une repr\'esentation lisse
irr\'eductible de $G'_F$, alors le {\it niveau} de $\pi$
(notation $niv(\pi)$) est par d\'efinition 0 si $\pi$ a un vecteur
fixe non nul sous $K_F$ et le plus petit entier positif $l$ tel que $\pi$ ait
un vecteur fixe non nul sous $K_F^l$ sinon. Il est alors connu qu'il y a une
correspondance bijective entre les repr\'esentations lisses
irr\'eductibles de $G'_F$ de niveau inf\'erieur ou \'egal \`a $l$ et les $H(G'_F,K_F^l)$-modules
irr\'eductibles. Pour \'etudier $H(G'_F,K_F^l)$, on part de la
d\'ecomposition de 
Bruhat 
$$G'_F=\coprod_{A\in \aa_F}K_FAK_F,$$ 
o\`u $\aa_F$ est l'ensemble des matrices $A=(a_{ij})_{1\leq i,j\leq r}$
telles que pour tout $i,j$ on ait $a_{ij}=\delta_{i,j}\pi_{D_F}^{a_i}$, avec
$a_{1}\leq a_2\leq...\leq a_r$. Si $l$ est un entier strictement
positif, alors $K_F^l$ est un sous-groupe distingu\'e de $K_F$ et on peut
donc obtenir une d\'ecomposition du type
$$G'_F=\coprod_{A\in \aa_F}\coprod_{(x,y)}K^l_FxAyK^l_F.$$
Alors, l'ensemble des fonctions caract\'eristiques des ensembles
$K^l_FxAyK^l_F$ forme une base de l'espace vectoriel
$H(G'_F,K_F^l)$. L'ensemble des couples $(x,y)$ s'identifie de
fa\c{c}on naturelle \`a un sous-ensemble de $(K_F/K_F^l)\times
(K_F/K_F^l)$. Or, on a 
$K_F/K_F^l\cong GL_r(O_{D_F}/P^{ld}_{D_F})$. L'id\'ee de Kazhdan se
traduit, dans notre situation, de la fa\c{c}on 
suivante~: si $L$ est un corps local non archim\'edien de
caract\'eristique non nulle, et $\pi_L$ une uniformisante fix\'ee de
$L$, si $m$ est un entier strictement positif, on dit que $L$ est
$m$-proche de $F$ (ou que $L$ et $F$ sont $m$-proches) s'il existe un
isomorphisme d'anneaux
$$\lll: O_{F}/P^{m}_{F} \to O_{L}/P^{m}_{L}$$
qui envoie la classe de $\pi_F$ sur la classe de $\pi_L$. Dans cette
situation, on peut consid\'erer sur $L$ une alg\`ebre \`a division
centrale de dimension $d^2$ qui ait le m\^eme invariant de Hasse que
$D_F$, et, si on adopte toutes les notations plus haut aussi pour le corps
$L$, obtenir  que $\lll$ induit naturellement un isomorphisme de
$O_{D_F}/P^{md}_{D_F}$ sur 
$O_{D_L}/P^{md}_{D_L}$. Ensuite on utilise les d\'ecompositions  
$$G'_F=\coprod_{A_F\in \aa_F}\coprod_{(x_F,y_F)}K^m_Fx_FA_Fy_FK^m_F$$
et
$$G'_L=\coprod_{A_L\in \aa_L}\coprod_{(x_L,y_L)}K^m_Lx_LA_Ly_LK^m_L.$$
L'isomorphisme de  $O_{D_F}/P^{md}_{D_F}$ sur 
$O_{D_L}/P^{md}_{D_L}$ induit un isomorphisme de
$GL_r(O_{D_F}/P^{md}_{D_F})\times GL_r(O_{D_F}/P^{md}_{D_F})$ sur
$GL_r(O_{D_L}/P^{md}_{D_L})\times GL_r(O_{D_L}/P^{md}_{D_L})$, et on
peut esp\'erer que cet isomorphisme r\'ealise  une bijection entre les
sous-ensembles 
$\{(x_F,y_F)\}$ et  $\{(x_L,y_L)\}$ de
ces deux groupes sur lesquels sont ind\'ex\'ees les partitions plus
haut. Ceci est une premi\`ere difficult\'e. On montre (lemme
\ref{3.2.7}) que, en
choisissant comme il faut les uniformisantes de $D_F$ et $D_L$,  on
trouve bien une telle bijection. On obtient alors, par les
consid\'erations plus haut, une application
bijective d'une base de  $H(G'_F,K_F^m)$ sur une base de
$H(G'_L,K_L^m)$, et donc un isomorphisme naturel d'{\it espaces vectoriels}  
entre les deux alg\`ebres de Hecke. Dans une situation id\'eale cette
application serait en fait un isomorphisme d'{\it alg\`ebres}
(conjecture de Kazhdan), et induirait
ainsi une bijection naturelle entre l'ensemble des repr\'esentations
lisses irr\'eductibles de niveau inf\'erieur ou \'egal \`a $m$ de
$G'_F$ et l'ensemble des 
repr\'esentations  
lisses irr\'eductibles de niveau inf\'erieur ou \'egal \`a $m$ de
$G'_L$. Remarquons que, si $L$ 
et $F$ sont $m$ proches, et si $m\geq l>0$, alors $L$ et $F$ sont
trivialement $l$-proches. Nous montrons dans la
section 2 (th.\ref{isomalg}) le r\'esultat suivant, plus faible que la
conjecture de Kazhdan~: pour tout entier strictement positif
$l$, il existe un 
entier $m\geq l$ tel que, si $L$ est un corps $m$ proche de $F$,
alors l'isomorphisme d'espaces vectoriels entre  $H(G'_F,K_F^l)$ et
$H(G'_F,K_F^l)$ induit \`a partir de cette proximit\'e est un
isomorphisme d'alg\`ebres.  Nous prouvons aussi que la bijection
induite entre l'ensemble des repr\'esentations lisses irr\'eductibles
de niveau inf\'erieur ou \'egal \`a $l$ de $G'_F$ et l'ensemble des
repr\'esentations lisses irr\'eductibles 
de niveau inf\'erieur ou \'egal \`a $l$ de $G'_L$
envoie les \cares\ sur des \cares\ et les \cusps\ sur des
\cusps\ (th.\ref{careseridica}), et qu'elle pr\'eserve les facteurs
$\epsilon'$\ (th.\ref{epsilonseridica}). Comme elle
pr\'eserve \'egalement le niveau, il s'ensuit que les r\'esultats
que nous avons prouv\'es dans [Ba3] uniquement en caract\'eristique
nulle (in\'egalit\'es concernant le niveau et le conducteur d'une
repr\'esentation), ainsi que leurs cons\'equences (th\'eor\`eme de
finitude pour les repr\'esentations automorphes cuspidales
d'une forme int\'erieure de $GL_n$), sont valables aussi en
caract\'eristique non nulle (section 3, th.\ref{niveau} et th.\ref{finitude}).

Supposons maintenant que nous voulions prouver
le th\'eor\`eme \ref{JL} quand le corps $F$ est de caract\'eristique
non nulle. Nous voulons d'abord fixer une \eci\ $\pi$ de $G$ et trouver
une \eci\ $\pi'$ de $G'$ telle qu'en posant ${\bf C}(\pi)=\pi'$ on ait
la relation \ref{JLe}. Consid\'erons un entier strictement positif
$l$. On peut trouver $m$ suffisamment grand pour que, si $L$ est un
corps local non archim\'edien de caract\'eristique nulle  $m$-proche
de $F$, on ait un
diagramme~:    
$$\diagram
{G_L&\hfl{1}{}&G'_L\cr
 && \cr
\vfl{2}{}&&\vfl{2'}{}\cr
 && \cr
G&........&G'\cr}
$$
o\`u

- la fl\`eche 1 est la correspondance d\'ej\`a \'etablie en
caract\'eristique nulle par [DKV], 

- la fl\`eche 2, cons\'equence de la
proximit\'e des corps $F$ et $L$, est une application
qui r\'ealise une bijection  
entre l'ensemble de classes de \rlis\  
de niveau inf\'erieur ou \'egal \`a $l$ de $G$ et l'ensemble de
classes de 
\rlis\ 
de niveau inf\'erieur ou \'egal \`a $l$ de $G_L$,

- la
fl\`eche 2' joue un r\^ole similaire 
pour $G'$ et $G'_L$,

- les pointill\'es en bas repr\'esentent pour l'instant notre d\'esir
de construire 
une 
correspondance (Jacquet-Langlands en
caract\'eristique non nulle). 

Fixons $l\geq niv(\pi)$ (sans quoi la fl\`eche 2
n'est pas d\'efinie pour $\pi$). On pense, naturellement, monter $\pi$ 
le long de $2$, transf\'erer le long de $1$, redescendre le long de
$2'$ et prouver ensuite que le r\'esultat ainsi obtenu est le bon candidat
au poste de $\pi'$. Voici maintenant les 
principales difficult\'es~:

- nous pouvons d\'eplacer $\pi$ le long de $2$ et ensuite le long
  de $1$, mais on n'est pas s\^ur de pouvoir descendre le r\'esultat le
  long de $2'$ ; en effet, il se peut qu'il soit de niveau sup\'erieur
  \`a $l$.  C'est vrai que la fl\`eche $2$ pr\'eserve le niveau, mais
  $1$ ne le pr\'eserve 
  certainement pas dans le sens o\`u il est d\'efini ici (m\^eme si on
  esp\`ere qu'il conserve le niveau normalis\'e tel qu'il est d\'efini
  dans la th\'eorie des types). Prendre $L$ plus proche revient \`a
  changer compl\`etement de diagramme, et les m\^emes probl\`emes
  recommencent. La solution est de borner uniform\'ement, \`a partir
  uniquement de $\pi$, le niveau de
  toutes les repr\'esentations susceptibles d'intervenir dans la
  d\'emonstration. 

- m\^eme si on pouvait descendre le long de $2'$ le r\'esultat ainsi
  obtenu, il est difficile de 
  prouver que ce qu'on trouve convient pour d\'efinir une correspondance
  de Jacquet-Langlands entre $G$ et $G'$. Les fl\`eches $2$ et $2'$
  sont de nature essentiellement lin\'eaire alors que la
  correspondance est de nature essentiellement harmonique.\\

Nous r\'eglerons le premier probl\`eme en montrant que le facteur
$\epsilon'$ d'une repr\'esentation (tel que d\'efinit dans
[GJ])
est 
conserv\'e \`a la fois 
par 1, 2, et 2'. Donc, le facteur $\epsilon'$ de la repr\'esentation
obtenue en appliquant $2$ et ensuite 
$1$ \`a $\pi$  est \'egal au facteur $\epsilon'$ de $\pi$.
Mais on sait borner le niveau d'une repr\'esentation \`a partir du
conducteur qui se lit sur son
facteur $\epsilon$ (voir la section 3).  Or, \`a facteur $\epsilon'$ fix\'e, ce
conducteur est born\'e. Tout cela montre  qu'en se donnant $\pi$ (et
m\^eme en se donnant uniquement le facteur $\epsilon'$ de $\pi$), on
se donne automatiquement un entier $l'$ tel que le niveau de l'image
par la fl\`eche 1 de l'image par la fl\`eche 2 de $\pi$ soit
inf\'erieur \`a $l'$. En fixant donc $\pi$, on choisit d\`es le d\'epart
$m$ (et donc $L$) de fa\c{c}on a ce qu'il fonctionne pour $l'$ (\`a la
place de $l$), et on est assur\'e du fait qu'on puisse promener $\pi$
sur le diagramme jusqu'\`a $G'$.

Pour r\'esoudre le deuxi\`eme probl\`eme, c'est plus difficile
  techniquement. Pour une relation du type \ref{JLe} il faut regarder
  des 
  couples $g\lra g'$ o\`u $g\in G$ et $g'\in G'$, mais les constructions
  faites avec des corps proches ne voient que des ouverts assez gros,
  pas des points. Nous serons donc
forc\'es de raisonner sur des voisinages de $g$ et $g'$ sur lesquels
les caract\`eres fonction des repr\'esentations qui nous int\'eressent
sont constants. La difficult\'e, une fois de plus, sera un probl\`eme
  de borne 
uniforme (de la ``taille'' des ouverts). Ceci explique la machinerie 
  calculatoire d\'evelopp\'ee 
\`a la section 4.
    
\ \\ 

Dans la deuxi\`eme section on \'etudie les groupes $GL_r(D)$ d\'efinis
sur des corps proches, en suivant les indications de Kazhdan ([Ka]). 
La construction est p\'enible ; les r\'esultats sont les th\'eor\`emes 
\ref{careseridica} et \ref{epsilonseridica}. 
Nous tirons dans la troisi\`eme section des cons\'equences
imm\'ediates de ces th\'eor\`emes~: la validit\'e des r\'esultats de
[Ba3] en caract\'eristique non nulle, ici les th.\ref{niveau} et
\ref{finitude}. Ces r\'esultats seront par
ailleurs utilis\'es 
dans la preuve de la correspondance (section 5). Dans la quatri\`eme
section, nous donnons quelques r\'esultats concernant 
l'analyse harmonique de deux groupes du type $GL_r(D)$ d\'efinis 
sur des corps proches et qui se correspondent comme dans la
section 2. Dans la cinqui\`eme section nous d\'emontrons enfin le
th\'eor\`eme 
\ref{JL} 
en caract\'eristique non nulle. Nous nous inspirons fortement de
[DKV]. Cet article a \'et\'e r\'edig\'e en caract\'eristique
nulle, et les auteurs utilisent plusieurs techniques que nous ne
pouvons pas nous permettre dans le cas de caract\'eristique non nulle.
Notamment, la
correspondance est prouv\'ee par une r\'ecurrence qui fait intervenir
en m\^eme temps un
transfert d'int\'egrales orbitales qu'on ne peut pas obtenir ici. Nous
avons s\'epar\'e la d\'emonstration de la correspondance, en la
rendant d\'ependante d'un seul r\'esultat \`a prouver~:
l'orthogonalit\'e des caract\`eres (connu en caract\'eristique
nulle mais non en caract\'eristique non nulle). Ici nous utilisons 
la th\'eorie des corps proches d\'ecrite aux sections 2 et 4 et nous
pratiquons un va-et-vient entre le th\'eor\`eme  
\ref{JL}   
et le th\'eor\`eme d'orthogonalit\'e des caract\`eres sur les groupes
$GL_r(D)$ en caract\'eristique non nulle (th.\ref{ortogonalitate}), en
prouvant, au passage, 
la validit\'e de ce dernier r\'esultat.
\ \\

Ce travail a \'et\'e men\'e dans le cadre d'une th\`ese de doctorat
sous la direction de Guy Henniart ; je lui suis
tr\`es  reconnaissant de m'avoir propos\'e un si beau
sujet, et de m'avoir encourag\'e et aid\'e tout au long de mes
recherches, avec la gentillesse et la disponibilit\'e qui le
caract\'erisent.  Je consid\'ererai toujours un honneur que d'avoir
\'et\'e accueilli ce temps durant par le Laboratoire
d'Arithm\'etique et G\'eom\'etrie Alg\'ebrique de l'Universit\'e Paris
Sud, Orsay, o\`u j'ai pu c\^otoyer des chercheurs exceptionnels,
merveilleux exemples pour tout jeune math\'ematicien. Je remercie
Bertrand Lemaire qui m'a beaucoup aid\'e par ses conseils et tous ceux
avec lequels j'ai eu des discussions math\'ematiques, et qui m'ont
souvent donn\'e des r\'eponses \`a des questions vitales pour mon
travail~: Marie-France Vign\'eras, Fran\c{c}ois Court\`es, Anne-Marie
Aubert, et les autres membres du s\'eminaire Groupes
R\'eductifs et Formes Automorphes. Je remercie ceux qui ont lu le
manuscrit avec attention, et qui m'ont fait des commentaires tr\`es
utiles~: Guy Henniart, Herv\'e Jacquet, Colette M{\oe}glin et Bertrand
Lemaire.

\newpage
\section{Corps locaux proches et formes int\'erieures du groupe lin\'eaire}

Soient $F$ un corps local de caract\'eristique non nulle, $D_F$ une
alg\`ebre \`a division centrale de dimension $d^2$ sur $F$ et $G'_F$
le groupe $GL_r(D_F)$. Soient $m\in \mathbb{N}$ et $L$ un corps local
$m$-proche de $F$ au sens de [Ka]. On voudrait d\'efinir un groupe
$G'_L$ sur $L$ tel qu'il y ait une ressemblance entre la th\'eorie des
repr\'esentations de $G'_F$ et celle de $G'_L$ \`a la mani\`ere dont
Kazhdan l'a fait pour les groupes de Chevalley dans [Ka] et Lemaire
pour $GL_n$ dans [Le]. L'id\'ee est simple~: choisir une alg\`ebre
\`a division $D_L$ centrale sur $L$ de dimension $n^2$ qui ait le
m\^eme invariant de Hasse que $D_F$, et poser $G'_L=GL_r(D_L)$. Dans
ce qui suit nous montrons que ce groupe est effectivement solution de
notre probl\`eme ; on construit $G'_L$ de fa\c{c}on \`a pouvoir
v\'erifier les th\'eor\`emes que Kazhdan a montr\'es pour les groupes
de Chevalley, et montrer aussi quelques autres r\'esultats utiles pour
la suite. Les r\'esultats de \ref{extnonram2} et \ref{algadiv2} sont
d\'emontr\'es dans un cadre g\'en\'eral dans [De]. On les a
red\'emontr\'es ici d'une fa\c{c}on concr\`ete dans ce cas particulier
pour fixer les notations  qu'on utilise par la suite.

\subsection{Pr\'eliminaires sur les alg\`ebres \`a
division}{\label{algadiv1}} 
 
Soient $F$ un corps local non archim\'edien, $v_F$ la valuation
normalis\'ee de $F$, $O_F$ l'anneau des entiers de $F$, $P_F$
l'id\'eal maximal de $O_F$ form\'e des \'el\'ements de $F$ de
valuation strictement positive, $k_F=O_F/P_F$ le corps r\'esiduel de $F$ (de
cardinal fini $q$) et $\pi_F$ une uniformisante de $F$. Soit
$D$ une alg\`ebre \`a division centrale sur $F$. On identifie $F$ au
sous-corps $1_DF$ de $D$. On sait que $dim_F(D)$ est un carr\'e
$d^2$. On sait \'egalement qu'il existe une valuation $v_D$ sur $D$
\`a valeurs dans $\mathbb{Z}$ surjective qui prolonge $dv_F$ ; on note
$O_D$ l'ensemble des \'el\'ements de $D$ de valuation positive ou
nulle et $P_D$ l'ensemble des \'el\'ements de $D$ de valuation
strictement positive ; $O_D$ est un anneau local non commutatif et $P_D$ est
l'id\'eal (bilat\`ere) maximal de $O_D$.  

On dispose d'une classification des alg\`ebres \`a division centrales
sur $F$ de dimension $d^2$~: si $D$ est une telle alg\`ebre, il existe
un  sous-corps commmutatif maximal $E$ de $D$ qui soit une extension
non ramifi\'ee de degr\'e $d$ de $F$ ; il existe \'egalement un
\'el\'ement $\pi_D$ de $D$ et un g\'en\'erateur $\sigma$ du groupe de
Galois de l'extension $E/F$ tels que~: 

- $\pi_D^d=\pi_F$

- $D=\bigoplus_{i=0}^{d-1} \pi_D^i E$

- pour tout $e\in E$ on a $\pi_D^{-1}e\pi_D=\sigma(e)$.
L'\'el\'ement $\sigma$ est uniquement d\'etermin\'e et permet de
calculer l'invariant de Hasse de $D$. 

R\'eciproquement, en se donnant un g\'en\'erateur $\sigma$ du groupe
de Galois de $E/F$ o\`u $E$ est une
extension non ramifi\'ee de degr\'e $d$ de $F$ (unique \`a
isomorphisme pr\`es), on peut construire une
alg\`ebre \`a division centrale sur $F$ de dimension $d^2$ en lui
imposant les trois conditions plus haut. Si on fixe $E$, deux telles
alg\`ebres sont 
isomorphes si et seulement si $\sigma\in Gal(E/F)$ est le m\^eme. Ce
sont les r\'esultats de la Proposition a, page 277, et Corollaire b,
page 335 de [Pi]. 

Nous commen\c{c}ons par la relation entre les
extensions non ramifi\'ees de m\^eme degr\'e sur deux corps locaux
proches. 

\subsection{Pr\'eliminaires sur les extensions non ramifi\'ees d'un
corps local}{\label{extnonram1}} 

Soient $F,v_F,O_F,P_F,\pi_F, k_F, q=card(k_F)$ comme plus haut. Soit
$d\in {\mathbb{N}}^*$. Si on fixe une cl\^oture alg\'ebrique $\bar{F}$
de $F$, alors il existe une unique extension $E$ non ramifi\'ee de
degr\'e $d$ de $F$ incluse dans $\bar{F}$. On a les propri\'et\'es
suivantes~: 

- la valuation normalis\'ee de $E$ prolonge celle de $F$ et $\pi_F$
  est une uniformisante de $E$ ; 
 
- si $l=q^d$ et $P$ est le polyn\^ome $P(X)=X^{l-1}-1\in F[X]$, alors
  $E$ est un corps de d\'ecomposition de $P$ ;  en particulier, pour
  toute racine primitive d'ordre $l-1$ de l'unit\'e $y_E$ dans
  $\bar{F}$, on a $y_E\in E$ et $E=F[y_E]$ ; 

- le cardinal du corps r\'esiduel $k_E=O_E/P_E$ de $E$ est \'egal \`a
$q^d$ et on a un isomorphisme de groupes $g_{E/F}$ de $Gal(E/F)$ sur
$Gal(k_E/k_F)$ ; cet isomorphisme est donn\'e par l'application
suivante~: si $\sigma\in Gal(E/F)$ alors $\sigma$ induit un
automorphisme $\sigma'$ de $O_E$ qui envoie $P_E$ sur $P_E$ et qui
agit comme l'identit\'e sur $O_F$ ; $\sigma'$ induit par cons\'equent
un isomorphisme $\sigma''$ de $k_E$ sur $k_E$ dont la restriction \`a
$k_F$ est l'identit\'e et on pose $g_{E/F}(\sigma)=\sigma''$.\\ 
C'est la Proposition 17.8, page 334, [Pi].

Soit $m\in {\mathbb{N}}^*$. On pose $O_{Fm}=O_F/P_F^m$ et
$O_{Em}=O_E/P_E^m$. Dans ce qui suit, une barre au-dessus d'un symbole
rappelle que ce symbole est rattach\'e (d'une fa\c{c}on ou d'une
autre) \`a $k_F$ ou $k_E$, et un chapeau au-dessus d'un symbole
rappelle qu'il est rattach\'e  \`a $O_{Fm}$ ou $O_{Em}$. Soient $P$ le
polyn\^ome $X^{l-1}-1\in F[X]$, $\hat{P}$ le polyn\^ome
$X^{l-1}-\hat{1}\in O_{Fm}[X]$ et $\bar{P}$ le polyn\^ome
$X^{l-1}-\bar{1}\in k_F[X]$. Soit
$\bar{P}=\bar{P_1}\bar{P_2}...\bar{P_k}$ une d\'ecomposition en
produit de polyn\^omes unitaires irr\'eductibles de $\bar{P}$ dans
$k_F[X]$. Supposons (sans restreindre la g\'en\'eralit\'e) que
$\bar{P_1}$ ait  une racine primitive d'ordre $l-1$ de l'unit\'e
$\bar{y}$ dans l'extension $k_E$ de $k_F$. On sait qu'alors le degr\'e
de $\bar{P_1}$ est \'egal \`a $d$. 

\begin{prop}\label{polynomes} a) Il existe un unique polyn\^ome
unitaire $\hat{P_1}\in O_{Fm}[X]$ tel que  

A) $\hat{P_1}$ divise $\hat{P}$ et 

B) l'image de $\hat{P_1}$ dans $k_F[X]$ soit $\bar{P_1}$ (en
particulier $\hat{P_1}$ est irr\'eductible). 

\ \ \ \ b) Soit $y$ l'unique racine de $P$ dans $E$ telle que l'image
de $y$ dans $k_E$ soit $\bar{y}$. Notons $\hat{y}$ l'image de $y$ dans
$O_{Em}$. Il existe un isomorphisme de $O_{Fm}$-alg\`ebres $\hat{f}_m
: O_{Fm}[X]/(\hat{P_1})\rightarrow O_{Em}$ induit par le morphisme
$f_m:O_{Fm}[X]\rightarrow O_{Em}$ donn\'e par $Q\mapsto Q(\hat{y})$. 
\end{prop}

\ \\
{\bf{D\'emonstration.}} a) On montre l'existence et l'unicit\'e. 

EXISTENCE~: Dans la d\'emonstration de la Proposition a, page 277 de
Pierce on montre que la d\'ecomposition de $P$ en produit de
polyn\^omes unitaires irr\'eductibles est du type $P=P_1P_2...P_k$
o\`u pour tout $i\in\{1,2...k\}$, l'image facteur de $P_i$ dans
$k_F[X]$ est $\bar{P_i}$ (en outre, les $\bar{P_i}$ sont tous distincts
parce qe $\bar{P}$ est sans racine multiple dans une cl\^oture
alg\'ebrique de $k_F$). On v\'erifie aussit\^ot que l'image facteur
de $P_1$ dans $O_{Fm}[X]$ satisfait \`a A) et B). 

UNICIT\'E~: Soit $\hat{M}$ un  polyn\^ome unitaire de
$O_{Fm}[X]$ qui satisfait \`a A) et B). Il suffit de montrer que
$\hat{M}$ a un repr\'esentant $M$ dans $O_F[X]$ qui divise
$P=P_1P_2...P_k$ (voir l'EXISTENCE) ; on aurait alors que $M$ est un
produit des $P_i$, et, par B), que $M=P_1$. 

\'Ecrivons  $\hat{P}=\hat{M}\hat{V}$ o\`u, $\hat{M}$ et $\hat{P}$
\'etant unitaires, $\hat{V}$ est unitaire. Soient $N$ un
repr\'esentant unitaire de $\hat{M}$ et $U$ un
repr\'esentant unitaire de $\hat{V}$ dans $O_F[X]$.  
On applique le
th\'eor\`eme 1 de [BC] IV, §3, plus pr\'ecis\'ement la remarque qui suit
la d\'emonstration. Cette variante plus fine du lemme de Hensel assure
que si $N$ et $U$ sont deux polyn\^omes unitaires tels que 
$$P\equiv
NU\mod(\pi_F^m),$$
et par ailleurs le r\'esultant des polyn\^omes $N$
et 
$U$ est {\it de valuation nulle}, alors la
d\'ecomposition ``se rel\`eve'' au sens o\`u il existe une
d\'ecomposition $P=MV$ dans $O_F[X]$ o\`u $M$ et $V$ sont tels que   
$$M\equiv
N\mod(\pi_F^m)$$
et
$$V\equiv
U\mod(\pi_F^m).$$
Or, le r\'esultant $\hat{R}$ de $\hat{M}$ et $\hat{V}$ est de
valuation nulle parce que, les polyn\^omes \'etant unitaires, la
r\'eduction de $\hat{R}$ modulo $\pi_F$ est \'egale au r\'esultant des
r\'eductions $\bar{M}$ et $\bar{V}$ modulo $\pi_F$, et ce dernier est
non nul dans $k_F$, parce que $\bar{P}=\bar{M}\bar{V}$ est sans
facteurs multiples, soit $\bar{M}$ et $\bar{V}$ sont forc\'ement
premiers entre eux. (Le lecteur remarquera que la difficult\'e de cet
exercice r\'eside dans le  fait que $O_{mF}$ n'est pas principal. En
utilisant que
$O_F$ et $k_F$ le sont on coince le probl\`eme entre les deux).

\ \\

b) L'application $f_m$ est visiblement un morphisme de
$O_{Fm}$-alg\`ebres. Montrons qu'elle est surjective. L'alg\`ebre
$O_{Em}$ est un $O_{Fm}$-module libre de rang $d$ dont une base est
${\mathcal{B}}=\{\hat{1},\hat{y},\hat{y}^2...\hat{y}^{d-1}\}$ ; en
effet, $\mathcal{B}$ est une famille g\'en\'eratrice car
$\{1;y;y^2...y^{d-1}\}$ 
\'etait une famille g\'en\'eratrice de $O_E$ sur $O_F$ ([We], prop.5,
page 20) et $\mathcal{B}$ est une famille libre car, $\hat{P_1}$
\'etant irr\'eductible, c'est le polyn\^ome minimal de $\hat{y}$. Or,
la base $\mathcal{B}$ se trouve dans l'image de $f_m$ donc $f_m$ est
bien surjective. Comme le noyau de $f_m$ est clairement l'id\'eal
principal engendr\'e par $\hat{P_1}$, le r\'esultat en d\'ecoule.
\qed

\begin{prop}\label{isom} On note $G_m$ l'ensemble des
$O_{Fm}$-automorphismes de l'alg\`ebre $O_{Em}$. Alors le morphisme
naturel $g_{m,E/F}:Gal(E/F)\rightarrow G_m$ est un isomorphisme de groupes.
\end{prop}

\ \\
{\bf{D\'emonstration.}} On rappelle qu'on a not\'e $g_{E/F}$
l'isomorphisme canonique  
$$Gal(E/F)\simeq Gal(k_E/k_F).$$ 
 S'il existait $\sigma$ et $\sigma'$ distincts dans $Gal(E/F)$ tels
que $g_{m,E/F}(\sigma)=g_{m,E/F}(\sigma')$ alors on aurait \'egalement
$g_{E/F}(\sigma)=g_{E/F}(\sigma')$ ce qui est impossible car $g_{E/F}$
est un isomorphisme. Donc $g_{m,E/F}$ est injective. D'autre part,
d'apr\`es  la proposition \ref{polynomes}b), un $O_{Fm}$-automorphisme
de l'alg\`ebre $O_{Em}$ est uniquement d\'etermin\'e par l'image de
$\hat{y}$, et par ailleurs cette image doit \^etre une racine de
$\hat{P_1}$. Comme le degr\'e de $\hat{P_1}$ est $d$, le cardinal de
$G_m$ est inf\'erieur ou \'egal \`a $d$ qui est le cardinal de
$Gal(E/F)$. Mais alors, l'application $g_{m,E/F}$ \'etant injective elle est
forc\'ement surjective et finalement bijective.
\qed

\subsection{Corps locaux proches et extensions non
ramifi\'ees}{\label{extnonram2}} 

Soit $F$ est un corps local de caract\'eristique non nulle. Si $L$ est un
corps local de carac\'eristique nulle et $m$ est un entier
sup\'erieur ou \'egal \`a 1, on dit que $F$ et $L$ sont
{\it $m$-proches} si $O_{Fm}$ et $O_{Lm}$ sont isomorphes en tant
qu'anneaux. On appelle alors un {\it triplet de $m$-proximit\'e} un
triplet $(\pi_F;\pi_L;\lll)$ o\`u $\pi_F$ est une uniformisante de
$F$, $\pi_L$ est une uniformisante de
$L$, et $\lll$ est un isomorphisme de $O_{Fm}$ sur $O_{Lm}$ qui envoie
la classe de $\pi_F$ sur la classe de $\pi_L$.
Soient $F$ et $L$ deux corps locaux $m$-proches ($m\geq 1$) et
$(\pi_F;\pi_L;\lll)$ un triplet de $m$-proximit\'e
correspondant. Soient $d\geq 1$, $E$ une extension non ramifi\'ee de
dimension $d$ de $F$ et $K$ une extension non ramifi\'ee de dimension
$d$ de $L$. 

\begin{theo}\label{extnonram} Les corps $E$ et $K$ sont $m$-proches.
\end{theo}

\ \\
{\bf{D\'emonstration.}} L'isomorphisme $\lll:O_{Fm}\rightarrow O_{Lm}$
s'\'etend de fa\c{c}on naturelle en un isomorphisme $\lll~:
O_{Fm}[X]\rightarrow O_{Lm}[X]$. Soient $\bar{P}_F \in k_F[X]$ et
$\hat{P}_F\in O_{Fm}[X]$ choisis comme $\bar{P}_1$ et $\hat{P}_1$ dans
la proposition \ref{polynomes}. Posons $\bar{P}_L=\lll(\bar{P}_F) \in
k_L[X]$ et $\hat{P}_L=\lll(\hat{P}_F)\in O_{Lm}[X]$. Alors on a un
isomorphisme  

\begin{equation}\label{iso1}
\lll~: O_{Fm}[X]/(\hat{P}_F)\simeq O_{Lm}/(\hat{P}_L).
\end{equation}
Maintenant on sait par le point b) de la proposition \ref{polynomes} que  

\begin{equation}\label{iso2}
O_{Fm}[X]/(\hat{P}_F)\simeq O_{Em},
\end{equation}
isomorphisme qui d\'epend du choix d'une racine de $\hat{P}_F$.  
D'autre part $\hat{P}_L$ est un polyn\^ome qui v\'erifie les
conditions A) et B) de la proposition (avec corps de base cette fois
le corps $L$). Par unicit\'e et par le point b) de la proposition
\ref{polynomes} (appliqu\'ee cette fois sur $L$) on a un isomorphisme  

\begin{equation}\label{iso3}
O_{Lm}[X]/(\hat{P}_L)\simeq O_{Km}
\end{equation}
qui d\'epend du choix d'une racine primitive de $\hat{P}_F$.  
Des isomorphismes \ref{iso1}, \ref{iso2} et \ref{iso3} on d\'eduit un
isomorphisme 
$$\llll:O_{Em}\simeq O_{Km}.$$
Le triplet $(\pi_F;\pi_L;\llll)$ est un triplet de $m$-proximit\'e
pour  les corps $E$ et $K$.
\qed 
\ \\
\ \\

\subsection{Corps locaux proches et alg\`ebres \`a division}{\label{algadiv2}}

Soient $F$ et $L$ deux corps locaux $m$-proches ($m\geq 1$) et soit
$D_F$ une alg\`ebre \`a division centrale de dimension $d^2$ sur
$F$. On choisit un sous-corps non ramifi\'e maximal $E$ de $D_F$ et
notons $\sigma_E$ le g\'en\'erateur du groupe de Galois de
l'extension $E/F$ qui correspond \`a $D_F$ comme dans
l'introduction. Soit $K$ une extension non ramifi\'ee de degr\'e $d$
de $L$. Il existe un isomorphisme canonique
$g_{E/F,K/L}:Gal(E/F)\simeq Gal(K/L)$, image de l'isomorphisme
$g:Gal(k_E/k_F)\simeq Gal(k_K/k_L)$ qui envoie le Frobenius sur le
Frobenius. A son tour, l'isomorphisme $g_{E/F,K/L}$ induit (par la
proposition \ref{isom} un isomorphisme $g_{m,E/F,K/L}:G_{m,E/F}\simeq
G_{m,K/L}$. Par transport de structure on a pour $\hat{\sigma}\in
G_{m,E/F}$ et $\hat{x}\in O_{Em}$~: 

\begin{equation}\label{galcom}
g_{m,E/F,K/L}(\hat{\sigma})(\llll(\hat{x}))=\llll(\hat{\sigma}(\hat{x})).
\end{equation}

On note $D_L$ l'alg\`ebre \`a division centrale sur $L$ de dimension
$d^2$ qui correspond \`a l'extension $K/L$ et \`a l'\'el\'ement
$\sigma_K=g_{E/F,K/L}(\sigma_E)$ du groupe de Galois $Gal(K/L)$. On
fixe une uniformisante $\pi_{D_L}$ de $D_L$ avec les propri\'et\'es de
\ref{algadiv1} (relatives \`a $\pi_L$). 

En reprenant les notations de \ref{algadiv1}  on a

\begin{equation}
D_F=\bigoplus_{i=1}^{d-1}\pi_{D_F}^iE
\end{equation} 

\begin{equation}
O_{D_F}=\bigoplus_{i=0}^{d-1}\pi_{D_F}^iO_E
\end{equation}

\begin{equation}
P_{D_F}^{md}=\bigoplus_{i=0}^{d-1}\pi_{D_F}^iP_E^{m}
\end{equation} 
d'o\`u

\begin{equation}
O_{D_F}/P_{D_F}^{md}=\bigoplus_{i=0}^{d-1}\pi_{D_F}^iO_E/P_E^{m}. 
\end{equation} 
On a donc un isomorphisme de groupes additifs 
$\llld:O_{D_F}/P_{D_F}^{md}\simeq O_{D_L}/P_{D_L}^{md}$ 
qui envoie  
$\sum_{i=0}^{d-1}\pi_{D_F}\hat{\lambda}_i$ sur 
$\sum_{i=0}^{d-1}\pi_{D_L}\llll(\hat{\lambda}_i)$ pour tout $d$-uplet
$\{\hat{\lambda} _0;\hat{\lambda}_1...\hat{\lambda}_{d-1}\}\in 
O_{Em}^d$. Montrons que cet isomorphisme est compatible avec la 
multiplication (raison pour laquelle on a choisi $\sigma_K$
correspondant \`a $\sigma_E$). En effet, il suffit de v\'erifier la
compatibilit\'e avec la multiplication sur deux \'el\'ements du type
$\pi_{D_F}^i\hat{x}$ et $\pi_{D_F}^j\hat{x}'$ o\`u $0\leq i,j\leq n-1$
et $\hat{x}, \hat{x}'\in O_{Em}$. Soit $\hat{\sigma}_E$ l'image de
$\sigma_E$ dans $G_{m,E/F}$.  
On a 

\begin{equation}
\pi_{D_F}^i\hat{x}\pi_{D_F}^j\hat{x}'=\pi_{D_F}^{i+j}
g_{m,E/F,K/L}(\hat{\sigma}_E^j)(\hat{x}) \hat{x}'\ \text{si}\ i+j<n
\end{equation}
et                   

\begin{equation}
\pi_{D_F}^i\hat{x}\pi_{D_F}^j\hat{x}'=\pi_{D_F}^{i+j-n}\hat{\pi}_F
g_{m,E/F,K/L}(\hat{\sigma}_E^j)(\hat{x})\hat{x}'\  \text{si}\ i+j\geq
n   
\end{equation} 
Finalement, en utilisant la relation \ref{galcom}  plus haut on a
obtenu le~: 
 
\begin{theo}\label{algadiv} La fl\`eche 
$\llld:O_{D_F}/P_{D_F}^{md}\rightarrow O_{D_L}/P_{D_L}^{md}$ d\'efinie
plus haut est un isomorphisme   
d'anneaux. 
\end{theo}
\ \\

\subsection{Bijections formelles}{\label{bijform}}

Prenons le cas le plus simple du groupe lin\'eaire sur deux corps
locaux proches $F$ et $L$. Ce cas a d\'ej\`a \'et\'e trait\'e par
Lemaire qui a montr\'e comment on peut ``transf\'erer'' certaines
parties ouvertes et compactes de $GL_n(F)$ \`a $GL_n(L)$ et les
implications qu'a ce transfert pour les th\'eories 
des repr\'esentations des deux groupes. Supposons maintenant qu'on
se soit donn\'e un \'el\'ement d'un de ces sous-ensembles de $GL_n(F)$
et un \'el\'ement du sous-ensemble de $GL_n(L)$ correspondant et qu'on
veuille comparer leurs polyn\^omes caract\'eristiques qui sont,
certes, \`a coefficients dans des corps diff\'erents, mais
proches. C'est un probl\`eme tr\`es concret si on se repr\'esente les
deux \'el\'ements comme des matrices dont on sait comparer les
\'el\'ements. Seulement, la fa\c{c}on abstraite dont
on d\'efinit les bijections entre des objets attach\'es \`a $GL_n(F)$
et $GL_n(L)$ respectivement ne le permet pas. On va alors d\'efinir
ici des {\it{bijections formelles}} (formelles parce qu'elles ne
respectent pas les op\'erations) qui nous permettront de traiter ce
genre de situation concr\`ete. Ce probl\`eme peut \^etre vu aussi
comme celui d'un choix {\it pr\'ecis} de repr\'esentant dans une
classe d'\'equivalence. 
\ \\

\subsubsection{Les corps locaux proches}{\label{bijformcorps}}

Soit $F$ un corps local de caract\'eristique non nulle. On choisit un
syst\`eme de repr\'esentants $S_F$ de $O_F/P_F$ dans $O_F$. 

Soit $m\in \mathbb{N}^*$, $L$ un corps $m$-proche de $F$ et
$(\pi_F;\pi_L;\lll)$  le triplet de $m$-proximit\'e associ\'e. Pour
tout  $x\in S_F$ on 
choisit un repr\'esentant $y(x)$ dans $O_L$ de l'image par $\lll$
(dans $O_L/P_L^m$) de la classe de $x$ dans $O_F/P_F^m$. L'ensemble
$S_L=\{y(x) : x\in S_F\}$ est un syst\`eme de repr\'esentants de
$O_L/P_L^m$ dans $O_L$. Pour simplifier les calculs, on impose la
condition suivante~: $0_F$ et $1_F$ font partie de $S_F$, et on a
$y(0_F)=0_L$ et $y(1_F)=1_L$. 

On note $\lam$ la bijection de $S_F$ sur $S_L$ qui envoie $x$ sur
$y(x)$. Si on se repr\'esente les \'el\'ements de $F$ et de $L$ par
des s\'eries \`a l'aide des uniformisantes $\pi_F$ et $\pi_L$ et des
syst\`emes de repr\'esentants $S_F$ et $S_L$ respectivement, on
obtient une bijection de $F$ dans $L$ qui prolonge $\lam$ -- et pour
laquelle on utilisera donc la m\^eme notation -- donn\'ee par~: 

$$\lam(\displaystyle\sum_{j=j_0}^{\infty}s_j\pi_F^j)= 
\sum_{j=j_0}^{\infty}\lam(s_j)\pi_L^j.$$ 
La bijection $\lam$ induit une bijection (bien d\'efinie) de
$O_F/P_F^m$ sur $O_L/P_L^m$ et cette bijection n'est autre que
l'isomorphisme $\lll$. On appelle $\lam$ une {\it{bijection
formelle}}. Notons les propri\'et\'es suivantes qui sont imm\'ediates~: 

\begin{equation}
\forall i\in \mathbb{N},\  \forall x\in F,\ \lam(\pi_F^ix)=\pi_L^i\lam(x),
\end{equation}

\begin{equation}
\forall x\in F,\ v_L(\lam(x))=v_F(x).
\end{equation}
\ \\

\subsubsection{Les extensions non ramifi\'ees}{\label{bijformextnonram}}

Soient $F$ et $L$ comme plus haut, $E$ une extension non ramifi\'ee de
$F$ de degr\'e $d$ et $K$ une extension non ramifi\'ee de $L$ de
degr\'e $d$. Nous reprenons les notations de \ref{extnonram1}.
Les corps $E$ et $K$ \'etant $m$-proches par le
th\'eor\`eme \ref{extnonram}, on peut bien entendu d\'efinir une
bijection formelle entre $E$ et $K$ comme plus haut en oubliant
compl\`etement les corps $F$ et $L$. Mais, pour d\'efinir une
bijection formelle entre deux alg\`ebres \`a division sur $F$ et $L$
respectivement, bijection qui ait une certaine propri\'et\'e utile par
la suite (voir \ref{bijformalgadiv}), on d\'efinit une bijection formelle entre $E$ et $K$ de la
fa\c{c}on suivante, plus pr\'ecise~: au lieu de partir comme dans la
sous-section 
\ref{bijformcorps} avec un syst\`eme de repr\'esentants $S_E$
{\it{quelconque}} de $O_E/P_E$ dans $O_E$, on fixe un isomorphisme
$\llam:k_E\simeq k_K$ compatible avec $\lll:k_F\simeq k_L$ et on
impose que 

- $S_E$ soit form\'e de 0 et de toutes les racines du polyn\^ome
  $P_E=X^{l-1}-1\in F[X]$, 

- $S_K$ soit form\'e  de 0 et de toutes les racines du polyn\^ome
  $P_K=X^{l-1}-1\in L[X]$, 

- la bijection $\llam$ entre $S_E$ et $S_K$ soit donn\'ee par
  l'application suivante~: si $y$ est une racine de $P_E$, et
  $\bar{y}$ est l'image de $y$ dans $k_E$, alors $\llam(y)=z$ o\`u $z$
  est l'unique racine de $P_K$ qui se trouve au-dessus de
  $\llll(\bar{y})\in k_K$.  

C'est \`a partir de ce choix (qui est, remarquons-le), compatible
avec la condition $y(0_F)=0_L$ et $y(1_F)=1_L$) qu'on \'etend
l'application $\llam$ en une bijection $\llam:E\simeq K$ comme dans la
sous-section \ref{bijformcorps}. 

On obtient ainsi une bijection qui a la propri\'et\'e d'\^etre
compatible avec l'action des \'el\'ements du groupe de Galois
$Gal(E/F)$. En effet, quel que soit $\sigma\in Gal(E/F)$ on v\'erifie
facilement que, pour tout $x\in E$,  

\begin{equation}\label{**}
\llam(\sigma(x))=g_{E/F,K/L}(\sigma)(\llam(x)).
\end{equation}
L'application $\llam$ est une bijection formelle entre les corps
locaux $E$ et $K$. Elle induit un isomorphisme $\llll:O_{Em}\simeq
O_{Km}$.  
\ \\

\subsubsection{Les alg\`ebres \`a division}{\label{bijformalgadiv}} 

Soient $F,E,D_F$ et $L,K,D_L$ comme dans la sous-section \ref{algadiv2}. On
suppose qu'on 
a construit des bijections formelles $\lam:F\simeq L$ et $\llam:E\simeq
K$ comme au \ref{bijformcorps} et \ref{bijformextnonram}. On construit
une bijection formelle entre $D_F$ et 
$D_L$ de la  fa\c{c}on (naturelle) suivante~: on pose, pour tout
$d$-uplet $(\alpha _0;\alpha_1...\alpha_{d-1})\in E^d$, 

\begin{equation}
\lamd(\sum_{i=0}^{d-1}\pi^i_{D_F}\alpha_i)=
\sum_{i=0}^{d-1}\pi^i_{D_L}\llam(\alpha_i). 
\end{equation} 
La bijection $\lamd$ induit une bijection (bien d\'efinie) de
$O_{D_F}/P_{D_F}^{md}$ sur $O_{D_L}/P_{D_L}^{md}$ et cette bijection
est l'isomorphisme $\llld$. La bijection $\lamd$ a aussi les
propri\'et\'es suivantes~:

\begin{equation}\label{pro1}
\forall x\in D_F,\ v_{D_L}(\lam(x))=v_{D_F}(x),
\end{equation}

\begin{equation}\label{pro2}
\forall i\in \mathbb{N},\ \forall x\in D_F,\ 
\lamd(\pi_{D_F}^ix)=\pi_{D_L}^i\lamd(x),
\end{equation} 

\begin{equation}\label{pro3}
\forall i\in \mathbb{N},\ \forall x\in D_F,\ 
\lamd(x\pi_{D_F}^i)=\lamd(x)\pi_{D_L}^i.
\end{equation} 
\ \\

Les propri\'et\'es \ref{pro1} et \ref{pro2} sont \'evidentes, seule la
propri\'et\'e \ref{pro3} pose un petit probl\`eme. Il suffit bien
s\^ur de la montrer pour $x\in E$. On sait que si $x\in E$ alors
$x\pi_{D_F}^i=\pi_{D_F}^i\sigma_E^i(x)$. Il suffit donc de d\'emontrer
que pour tout $x\in E$ on a
$\llam(\sigma_E^i(x))=\sigma_K^i(\llam(x))$. Mais c'est la relation
\ref{**} (et c'\'etait justement pour avoir cette propri\'et\'e
\ref{pro3} qu'on avait choisi $\llam$ comme dans \ref{bijformextnonram}). 
\ \\

\subsubsection{Les matrices, les polyn\^omes}{\label{bijformmatrices}}
Si $F,D_F$ et $L,D_L$ sont comme au \ref{bijformalgadiv}, la bijection
formelle entre $F$ et $L$ s'\'etend de fa\c{c}on naturelle en une
bijection entre $F^n$ et $L^n$. Si $U$ est un $F$-espace vectoriel de
dimension finie $n$ muni d'une base et $V$ un $L$-espace vectoriel
de dimension $n$ muni d'une base, une bijection formelle entre $F$ et
$L$ s'\'etend ``composante par composante'' en une bijection entre
$U$ et $V$. C'est pareil pour des espaces vectoriels \`a droite ou \`a
gauche sur $D_F$ et $D_L$. En particulier cela vaut pour les espaces
de matrices $M_n(F)$, $M_n(L)$ et $M_r(D_F)$,
$M_r(D_L)$, et aussi pour des espaces de polyn\^omes. En g\'en\'eral,
on notera $\lam$ la bijection formelle 
donn\'ee entre $F$ et $L$, et $\zzzf$ la bijection induite entre
$M_n(F)$ et $M_n(L)$. Pareillement, la bijection formelle entre $M_r(D_F)$ et
$M_r(D_L)$ induite par $\lamd$ sera not\'ee $\zzzd$. \\  
\ \\

\subsection{La construction}{\label{construction}}

Soient $F,E,D_F,L,K,D_L$ comme plus haut et soit $r\in
{\mathbb{N}}^*$. On se propose d'\'etudier les ressemblances entre les
groupes $G'_F=GL_r(D_F)$ et $G'_L=GL_r(D_L)$. On va suivre le chemin
indiqu\'e par Kazhdan dans [Ka]. On pose $K_F=GL_r(O_{D_F})$,
$K_L=GL_r(O_{D_L})$ et, pour tout $l\in \mathbb{N}^*$, on note $K_F^l$
le noyau de la projection $K_F\rightarrow GL_r(O_{D_F}/P_{D_F}^{ld})$
et $K_L^l$ le noyau de la projection $K_L\rightarrow
GL_r(O_{D_L}/P_{D_L}^{ld})$. Pour tout $l$, $K_F^l$ est un sous-groupe
distingu\'e de $K_F$. C'est aussi le groupe
$Id+M_r(P_{D_F}^{ld})$. Si $K$ est un sous-groupe ouvert compact de
$G'_F$ on note $H(G'_F;K)$ l'alg\`ebre des fonctions complexes localement
constantes \`a support compact sur $G'_F$ et bi-invariantes par 
$K$. On fait de m\^eme pour $G'_L$.\\  
\ \\
{\bf{Notation~:}} Si $g=(g_{ij})\in M_r(D_F)$ alors on note $v_F(g)$
le minimum des valuations en tant qu'\'el\'ements de $D_F$ des
coefficients $g_{ij}$ de $g$. De m\^eme sur~$L$.\\ 

Maintenant, pour tout $l\in \mathbb{N}^*$, si $F$ et $L$ sont
$l$-proches, on construit  un isomorphisme d'espaces vectoriels
$$\zzdl : H(G'_F;K_F^l)\simeq H(G'_L;K_L^l).$$ 
Cet isomorphisme sera d\'ependant du triplet de $l$-proximit\'e
$(\lamd;\pi_{D_F};\pi_{D_L})$ d\'eduit du triplet de $l$-proximit\'e
fix\'e pour $F$ et $L$ comme nous l'avons expliqu\'e au
\ref{algadiv2}.  

Posons 
$${\mathcal{A}}_F=\{A=(a_{ij})_{i,j}\in GL_r(D_F) :
a_{ij}=\delta_{ij}\pi_{D_F}^{a_i},\ a_1\leq a_2\leq ...\leq a_r\}$$ 
o\`u $\delta_{ij}$ est le symbole de Kroneker.\\ 

\begin{lemme} On a $$G'_F={\bf{\coprod}}_{A\in {\mathcal{A}}}K_FAK_F.$$
\end{lemme}
\ \\
{\bf{D\'emonstration.}} Voir [Sa], page 43.
\qed
\ \\

Soit $l\in {\mathbb{N}^*}$ et $A\in {\mathcal{A}}_F$. Comme $K_F^l$ est un sous-groupe de $K_F$, il existe un sous-ensemble $X$ de $K_F\times K_F$ tel que 
$$K_FAK_F={\bf{\coprod}}_{(B;C)\in X}K_F^lBAC^{-1}K_F^l.$$ 
Posons $$\T =GL_r(O_{D_F}/P_{D_F}^{ld})\times GL_r(O_{D_F}/P_{D_F}^{ld}).$$ 
Si $(B;C)\in K_F\times K_F$, alors $K_F^lBAC^{-1}K_F^l$ ne d\'epend
que de la classe $(\hat{B};\hat{C})$ de $(B;C)$ dans $\T$. On peut
donc d\'efinir sans ambiguit\'e $K_F^l\hat{B}A\hat{C}^{-1}K_F^l$ pour
tout $(\hat{B};\hat{C})\in \T$. 
Notons maintenant ${\bf{H}}_{F,l,A}$ l'ensemble des couples
$(\hat{B};\hat{C})\in \T$ tels qu'il existe un repr\'esentant $B$ de
$\hat B$ dans $GL_r(O_{D_F})$ et un repr\'esentant $C$ de $\hat C$
dans $GL_r(O_{D_F})$ tel qu'on ait $BA=AC$. En \'ecrivant cette
relation $BAC^{-1}=A$ on v\'erifie que ${\bf{H}}_{F,l,A}$ est un
sous-groupe de $\T$. Posons enfin  
$$\TT=\T/{\bf{H}}_{F,l,A}.$$\\

\begin{lemme}\label{3.2.6} a) Pour tout $(\hat{B};\hat{C})\in \T$,
l'ensemble $K_F^l\hat{B}A\hat{C}^{-1}K_F^l$ ne d\'epend que de la
classe $(\tilde{B};\tilde{C})$ de $(\hat{B};\hat{C})$ dans $\TT$ ; on
le note $K_F^l\tilde{B}A\tilde{C}^{-1}K_F^l$. 

b) On a $$K_FAK_F={\bf{\coprod}}_{(\tilde{B};\tilde{C})\in
\TT}K_F^l\tilde{B}A\tilde{C}^{-1}K_F^l.$$ 
\end{lemme}
\ \\
{\bf{D\'emonstration.}} a) Soient $(\tilde{B};\tilde{C})\in \TT$ et
$(\hat{B};\hat{C})\in \T$, $(\hat{B}';\hat{C}')\in \T$ deux
repr\'esentants de  $(\tilde{B};\tilde{C})$. Il existe donc
$(\hat{U};\hat{V})\in {\bf{H}}_{F,l,A}$ tel que
$\hat{B}'=\hat{B}\hat{U}$ et $\hat{C}'=\hat{C}\hat{V}$. Soit $(B;C)$
un repr\'esentant de $(\hat{B};\hat{C})$ dans $GL_r(O_{D_F})\times
GL_r(O_{D_F})$ et $U$ et $V$ des repr\'esentants de $\hat{U}$ et de
$\hat{V}$ dans $GL_r(O_{D_F})$ qui v\'erifient $UA=AV$. On a donc
$K_F^l\hat{B}'A\hat{C'}^{-1}K_F^l=K_F^l\hat{B}\hat{U}A\hat{V}^{-1}\hat{C}^{-1}K_F^l=K_F^lBUAV^{-1}C^{-1}K_F^l=K_F^lBAC^{-1}K_F^l=K_F^l\hat{B}A\hat{C}^{-1}K_F^l$. 
\ \\

b) On a 
 
$$K_FAK_F=\bigcup_{(B;C)\in GL_r(O_{D_F})\times
GL_r(O_{D_F})}K_F^lBAC^{-1}K_F^l$$
et pour tout $(B;C)\in GL_r(O_{D_F})\times GL_r(O_{D_F})$ on a 

$$K_F^lBAC^{-1}K_F^l=K_F^l\tilde{B}A\tilde{C}^{-1}K_F^l.$$
Donc, 

$$K_FAK_F=\bigcup_{(\tilde{B};\tilde{C})\in
\TT}K_F^l\tilde{B}A\tilde{C}^{-1}K_F^l.$$
Il suffit de montrer que la r\'eunion est bien disjointe. Soient
$(\tilde{B};\tilde{C})$ et $(\tilde{B'};\tilde{C'})$ deux \'el\'ements
de $\TT$. Alors $K_F^l\tilde{B}A\tilde{C}^{-1}K_F^l$ et
$K_F^l\tilde{B'}A\tilde{C'}^{-1}K_F^l$ sont ou disjointes ou
\'egales. Supposons que  

$$K_F^l\tilde{B}A\tilde{C}^{-1}K_F^l=K_F^l\tilde{B'}A\tilde{C'}^{-1}K_F^l.$$
Comme $K_F^l$ est distingu\'e dans $K_F$ on a 

$$K_F^l\tilde{B}A\tilde{C}^{-1}K_F^l=\tilde{B}K_F^lAK_F^l\tilde{C}^{-1}$$
 et 

$$K_F^l\tilde{B'}A\tilde{C'}^{-1}K_F^l=\tilde{B'}K_F^lAK_F^l\tilde{C'}^{-1}.$$
 En posant  $\tilde{X}=\tilde{B'}^{-1}\tilde{B}$ et
$\tilde{Y}=\tilde{C'}^{-1}\tilde{C}$ et en r\'eutilisant le fait que
$K_F^l$ est distingu\'e dans $K_F$, on obtient 

$$K_F^l\tilde{X}A\tilde{Y}^{-1}K_F^l=K_F^lAK_F^{l}.$$
 Il suffit donc de montrer que dans ce cas on a forc\'ement
$\tilde{X}=\tilde{Y}=\tilde{1}$.  
C'est pareil que de montrer que, si $(\hat{X};\hat{Y})$ est dans $\T$
et v\'erifie $K_F^l\hat{X}A\hat{Y}^{-1}K_F^l=K_F^lAK_F^{l}$, alors
$(\hat{X};\hat{Y})$ a un repr\'esentant $(X;Y)$ dans
$GL_r(O_{D_F})\times GL_r(O_{D_F})$ tel que $XA=AY$. Montrons que
cette assertion est vraie~:
$K_F^l\hat{X}A\hat{Y}^{-1}K_F^l=K_F^lAK_F^{l}$ implique qu'il existe
un repr\'esentant $(X';Y')$ de $(\hat{X};\hat{Y})$ dans
$GL_r(O_{D_F})\times GL_r(O_{D_F})$ et deux \'el\'ements $k_1$ et
$k_2$ de $K_F^l$ tels qu'on ait $X^{'}AY^{'-1}=k_1Ak_2$. Mais alors il
suffit de prendre $X=k_1^{-1}X'$ et $Y=k_2Y'$, car l'appartenance de
$k_1$ et $k_2$ \`a $K_F^l=Id+M_r(P_{D_F}^{ld})$ nous garantit
$\hat{X}=\hat{X}'$ et $\hat{Y}=\hat{Y}'$ et avec ce choix on a bien
$XA=AY$.
\qed 
\ \\

Maintenant on va traiter les corps $F$ et $L$ \`a la fois. Les objets
d\'efinis plus haut sur $F$ se d\'efinissent de m\^eme sur $L$ et
c'est l'indice $F$ ou $L$ qui indique \`a tout moment de quel corps il
s'agit. Les r\'esultats plus haut sont \'evidemment valables aussi sur
$L$.  

On se rappelle l'isomorphisme $\llldl:O_{D_F}/P_{D_F}^{ld}\simeq
O_{D_F}/P_{D_F}^{ld}$. Il induit un isomorphisme $\zzdl~: \T\simeq
\Tl$.  

On se rappelle \'egalement la bijection formelle
$\zzzdl:M_r(D_F)\simeq M_r(D_L)$. Elle induit une bijection $\zzzdl :
{\mathcal{A}}_F\simeq {\mathcal{A}}_L$.\\ 

\begin{lemme}\label{3.2.7} Pour tout $A\in {\mathcal{A}}_F$,
l'isomorphisme $\zzdl : \T\simeq \Tl$ induit un isomorphisme $\zzdl :
{\bf{H}}_{F,l,A}\simeq {\bf{H}}_{L,l,\zzzdl(A)}$ et par cons\'equent
une bijection $\zzdl : \TT\simeq \TTl$. 
\end{lemme} 
\ \\
{\bf{D\'emonstration.}} Supposons que $A$ soit la matrice
$diag(\pi_{D_F}^{a_1},\pi_{D_F}^{a_2}...\pi_{D_F}^{a_r})$. \\ 
Soit $(\hat{B};\hat{C})\in {\bf{H}}_{F,l,A}$. Il existe donc un
repr\'esentant $(B;C)$ de $(\hat{B};\hat{C})$ dans
$GL_r(O_{D_F})\times GL_r(O_{D_F})$ tel que 

\begin{equation}\label{*} 
BA=AC
\end{equation}
\'Ecrivons $B=(b_{ij})$ et $C=(c_{ij})$. La relation \ref{*} se traduit
par~:  

\begin{equation}
\forall i,j,\ b_{ij}\pi_{D_F}^{a_j}=\pi_{D_F}^{a_i}c_{ij}
\end{equation}
Par les propri\'et\'es \ref{pro2} et \ref{pro3} de la bijection
formelle $\zzdl$ on a alors~: 
\begin{equation}
\forall i,j,\
\zzzdl(b_{ij})\pi_{D_L}^{a_j}=\pi_{D_L}^{a_i}\zzzdl(c_{ij})
\end{equation} 
et cette relation implique $\zzzdl(B)\zzzdl(A)=\zzzdl(A)\zzzdl(C)$ et
par cons\'equent $\bigl(\zzdl(\hat{B});\zzdl(\hat{C})\bigl)\in
{\bf{H}}_{L,l,\zzzdl(A)}$. Donc $\zzdl( {\bf{H}}_{F,l,A})\subset
{\bf{H}}_{L,l,\zzzdl(A)}$. Comme les r\^oles de $F$ et $L$ sont
sym\'etriques il est \'evident que ce r\'esultat est suffisant pour en
d\'eduire que $\zzdl : {\bf{H}}_{F,l,A}\rightarrow
{\bf{H}}_{L,l,\zzzdl(A)}$ est un isomorphisme, par exemple parce qu'on
peut lui construire une application r\'eciproque en partant de
$(\zzdl) ^{-1}: \Tl\simeq \T$.
\qed 
\ \\

D'apr\`es les lemmes \ref{3.2.6} et \ref{3.2.7}, si pour tout ensemble
$W$ on note ${\bf{1}}_W$ la fonction caract\'eristique de $W$, alors
l'ensemble~: 

$$\{{\bf{1}}_{K_F^l\tilde{B}A\tilde{C}^{-1}K_F^l} : A\in
{\mathcal{A}}_F, (\tilde{B};\tilde{C})\in \TT\} $$
est une base de l'espace vectoriel $H(G'_F;K_F^l)$. Soit 

$$\zzdl:H(G'_F;K_F^l)\rightarrow H(G'_L;K_L^l)$$ 
l'application lin\'eaire d\'etermin\'ee par 

$$\zzdl({\bf{1}}_{K_F^l\tilde{B}A\tilde{C}^{-1}K_F^l})=
{\bf{1}}_{K_L^l\bigl(\zzdl(\tilde{B})\bigl)\bigl(\zzzdl(A)\bigl)\bigl(\zzdl(\tilde{C})\bigl)^{-1}K_L^l}.$$ 
 
\begin{theo} 
L'application $\zzdl$ est un isomorphisme d'espaces vectoriels.
\end{theo}
\ \\
{\bf D\'emonstration.} Cela d\'ecoule des lemmes \ref{3.2.6} et
\ref{3.2.7}.
\qed
\ \\
\begin{rem}{\label{lcomutativ}}
{\rm Soit $l$ un entier strictement positif.  La surjection canonique 
$$O_{D_F}/P_{D_F}^{(l+1)d}\to O_{D_F}/P_{D_F}^{ld}$$
induit une  surjection canonique~:
$${\bf T}_{F,l+1}\to {\bf T}_{F,l}$$
dont la restriction induit une surjection
$${\bf H}_{F,l+1,A}\to {\bf H}_{F,l,A}.$$
Finalement, on obtient  une surjection canonique
$$\tilde{\bf T}_{F,l+1,A}\to \tilde{\bf T}_{F,l,A}.$$
Supposons que $F$ et $L$ sont $(l+1)$-proches. De la m\^eme fa\c{c}on,
on a une surjection canonique  
$$\tilde{\bf T}_{L,l+1,{\zeta}_{D_FD_L}^{l+1}(A)}\to \tilde{\bf
T}_{F,l,{\zeta}_{D_FD_L}^{l+1}(A)}.$$ 
Il est tr\`es facile \`a v\'erifier que le diagramme
$$\diagram
{\tilde{\bf T}_{L,l+1,{\zeta}_{D_FD_L}^{l+1}(A)}&\hfl{}{}&\tilde{\bf
T}_{F,l,{\zeta}_{D_FD_L}^{l+1}(A)}\cr 
\vfl{\bar{\zeta}_{D_FD_L}^{l+1}}{}&&\vfl{\bar{\zeta}_{D_FD_L}^{l}}{}\cr
\tilde{\bf T}_{F,l+1,A}&\hfl{}{}&\tilde{\bf T}_{F,l,A}\cr} 
$$
o\`u on aura tenu compte du fait que ${\zeta}_{D_FD_L}^{l+1}(A)= {\zeta}_{D_FD_L}^{l}(A)$ (par les d\'efinitions m\^emes) est commutatif. Cela implique que l'isomorphisme d'espaces vectoriels 
$$\zzdl:H(G'_F;K_F^l)\simeq H(G'_L;K_L^l)$$
est induit par la restriction de l'isomorphisme
$$\bar{\zeta}_{D_FD_L}^{l+1}:H(G'_F;K_F^{l+1})\simeq H(G'_L;K_L^{l+1}).$$}
\end{rem}
\ \\

\begin{lemme}\label{prevolumes}
Soit 
$A=diag(\pi_{D_F}^{a_1};\pi_{D_F}^{a_2}...\pi_{D_F}^{a_r})\in
{\mathcal A}_F$. Alors 
$$vol(K_F^lAK_F^l)=q^{d\sum_{i< j} a_j-a_i}vol(K_F^l).$$
\end{lemme}
\ \\
{\bf D\'emonstration.} On a  
$$vol(K_F^lAK_F^l)=card\big( K_F^l/(AK_F^lA^{-1}\cap K_F^l)\big)vol(AK_F^l)$$
$$=card\big(K_F^l/(AK_F^lA^{-1}\cap K_F^l)\big)vol(K_F^l).$$
Il suffit donc de montrer que 
\begin{equation}\label{suffit}
card\big(K_F^l/(AK_F^lA^{-1}\cap
K_F^l)\big)=q^{d\sum_{i< j} a_j-a_i}.
\end{equation}
\def\x{{\mathcal X}}

Posons $\x_A=AM_r(P_{D_F}^{dl})A^{-1}\cap
M_r(P_{D_F}^{dl})$ (en particulier,
$\x_1=M_r(P_{D_F}^{dl})$). L'ensemble $\x_A$ est 
form\'e des matrices $X=(x_{ij})_{1\leq i,j\leq r} \in M_r(P_{D_F}^{dl})$
qui v\'erifient $v_{D_F}(\pi_D^{a_i}x_{ij}\pi_D^{-a_j})\geq dl$ pour
tout $i<j$, soit des matrices $X$ qui v\'erifient $X\in
M_r(P_{D_F}^{dl})$ et pour tout 
$i<j,\ v_{D_F}(x_{ij})\geq dl+a_j-a_i$.\\
Par cons\'equent, le cardinal du groupe additif
$\x_1/\x_A$ est $q^{d\sum_{i< j} a_j-a_i}$. 
Ce que nous nous voulons montrer
(eq. \ref{suffit}) se
traduit par 
$$card(1+\x_1/1+\x_A)=q^{d\sum_{i< j} a_j-a_i}=card(\x_1/\x_A).$$ 
Disons que
$A\in {\mathcal A}_F$ est {\it petite} si on a
$max_{i<j}(a_j-a_i)\leq ld$. Des calculs simples montrent que, si $A$
est petite, alors 

- $\x_1 \ ^.\x_A\subset \x_A$, donc $1+\x_A$ est un sous-groupe
distingu\'e de $1+\x_1$ multiplicatif, 

- $\x_1 \ ^.\x_1\subset \x_A$, donc l'application $1+x\mapsto x$ de $1+\x_1$
dans $\x_1$ induit un isomorphisme de groupes
$$(1+\x_1/1+\x_A;\ ^.)\cong (\x_1/\x_A;+).$$
Nous avons ainsi le r\'esultat voulu (eq. \ref{suffit}) si $A$ est
petite. Comme tout \'el\'ement $A\in {\mathcal A}_F$ s'\'ecrit comme
produit de matrices petites de ${\mathcal A}_F$, pour avoir le
r\'esultat en g\'en\'eral il suffit de remarquer que, pour toute  $A\in
{\mathcal A}_F$ et toute $A'\in
{\mathcal A}_F$ petite on a 

- $\x_A \ ^.\x_{AA'}\subset \x_{AA'}$, donc $1+\x_{AA'}$ est un
sous-groupe distingu\'e de $1+\x_A$ multiplicatif, 

- $\x_A \ ^.\x_A\subset \x_{AA'}$, donc l'application $1+x\mapsto x$ de $1+\x_A$
dans $\x_A$ induit un isomorphisme de groupes
$$(1+\x_A/1+\x_{AA'};\ ^.)\cong (\x_A/\x_{AA'};+).$$
\qed

Pour tout $g\in G'_F$ on pose
$h(g)=(vol(K_F^l))^{-1}{\bf{1}}_{K_F^lgK_F^l}$. On a le~:

\begin{lemme}\label{volumes} 
a) Pour tout $A,A'\in {\mathcal{A}}_F$ on a $h(A)*h(A')=h(AA')$.

b) Pour tout $B,C\in GL_r(O_{D_F})\times GL_r(O_{D_F})$ on a $h(B)*h(A)*h(C)=h(BAC)$.
\end{lemme}
\ \\
{\bf{D\'emonstration.}} a) Par la proposition 2.2, chapitre 3 de [Ho], il suffit de d\'emontrer qu'on a 
$$vol(K_F^lAK_F^l)vol(K_F^lA'K_F^l)=vol(K_F^l)vol(K_F^lAA'K_F^l).$$
Ceci est une cons\'equence directe du lemme \ref{prevolumes}.

b) Par la m\^eme proposition 2.2 de [Ho], il suffit de montrer qu'on a~:
$$vol(K_F^lBK_F^l)vol(K_F^lACK_F^l)=vol(K_F^l)vol(K_F^lBACK_F^l)$$
et
$$vol(K_F^lAK_F^l)vol(K_F^lCK_F^l)=vol(K_F^l)vol(K_F^lACK_F^l).$$
C'est une relation facile \`a obtenir parce que $K_F^l$ est
distingu\'e dans $K_F$ 
(donc on peut ``sortir''$B$ et $C$) et les volumes sont pris par
rapport \`a une mesure de Haar \`a droite et \`a gauche (donc
finalement on peut effacer $B$ et $C$ partout o\`u ils
apparaissent). 
\qed

\begin{rem}{\label{remarquevolume}} 
{\rm Si $W_F$ est un ensemble ouvert et compact de $G'_F$ invariant
\`a gauche et \`a droite par $K_F^l$ et $L$ est un corps $l$-proche de 
$F$, la construction de l'isomorphisme $\zzdl$ implique que l'image de 
la fonction caract\'eristique de $W_F$ est la fonction
caract\'eristique d'un ensemble ouvert compact $W_L$ de $G'_L$ qui est
invariant \`a gauche et \`a droite par $K_L^l$. Le calcul de volumes
qu'on fait dans la d\'emonstration du lemme \ref{volumes} plus haut
montre qu'on a alors~: 
$$vol(W_F)=vol(W_L).$$ On pose $\zzdl (W_F)=W_L$ }
\end{rem}
\ \\ 
\begin{theo}\label{isomalg} 
Soit $l\in \mathbb{N}^*$. Il existe un entier $m$ tel que, si les
corps $F$ et $L$ sont $m$-proches, alors l'isomorphisme d'espaces
vectoriels $\zzdl~: H(G'_F;K_F^l)\simeq H(G'_L;K_L^l)$ est un
isomorphisme d'alg\`ebres. 
\end{theo}
\ \\
{\bf{D\'emonstration.}} On utilise la d\'emarche de  Kazhdan en
pr\'ecisant certains d\'etails~: 
\ \\

\begin{lemme}\label{cpct}  
Soit $\mathcal{C}$ un sous-ensemble fini de ${\mathcal{A}}_F$ et soit  
$$G'_F({\mathcal{C}})=\bigcup_{A\in \mathcal{C}}K_FAK_F.$$   

a) Il existe $m>l$ qui d\'epend de $\mathcal{C}$ tel qu'on ait, pour
tout $g\in G'_F({\mathcal{C}})$, $gK_F^mg^{-1}\subset K_F^l$. 

b) Supposons que $L$ soit $m$-proche de $F$. Alors pour tout $f_1,\
f_2\in H(G'_F;K_F^l)$ \`a support dans $G'_F({\mathcal{C}})$ on a  
$$\zzdl(f_1*f_2)=\zzdl(f_1)*\zzdl(f_2).$$
\end{lemme}
\ \\
{\bf{D\'emonstration.}} a) Il faut trouver un $m$ tel qu'on ait
$K_F^m\subset \bigcap_{g\in G'_F({\mathcal{C}})}g^{-1}K_F^lg$. Il
suffit de prendre 
$$m\geq l+\sup\displaystyle_{g\in
G'_F({\mathcal{C}})}v_F(g^{-1})+\text{sup}_{ g\in
G'_F({\mathcal{C}})}v_F(g)$$
i.e.
$$m\geq l+\text{max}_{A\in
{\mathcal{C}}}v_F(A^{-1})+\text{max}_{A\in {\mathcal{C}}}v_F(A).$$ 
\ \\

b) Il suffit de montrer ce r\'esultat pour
$f_1={\bf{1}}_{K_F^lgK_F^l}$ et $f_2={\bf{1}}_{K_F^lg'K_F^l}$, avec
$g,g'\in G'_F({\mathcal{C}})$. En revenant \`a la d\'efinition du
produit de convolution on trouve~: 
$${\bf{1}}_{K_F^lgK_F^l}*{\bf{1}}_{K_F^lg'K_F^l}(x)=vol\bigl(K^l_FgK_F^l\cap   
K^l_Fg'K_F^l x\bigl).$$  
Par le point a), cette intersection est bi-invariante par $K_F^m$ et on a~:
$${\bf{1}}_{K_F^lgK_F^l}*{\bf{1}}_{K_F^lg'K_F^l}=$$
$$=\sum_{A\in {\mathcal{A}}_F\ }\sum_{(\tilde{B};\tilde{C})\in
\TT}vol\bigl(K^l_FgK_F^l\cap K^l_Fg'K_F^l
\tilde{B}A\tilde{C}\bigl){\bf{1}}_{K_F^m\tilde{B}A\tilde{C}^{-1}K_F^m}.$$ 
Les corps $F$ et $L$ \'etant $m$-proches, $\zzd$ est bien
d\'efinie. D'apr\`es la formule plus haut qui vaut aussi bien sur $L$
que sur $F$, et la remarque \ref{remarquevolume} sur les volumes,
$\zzd(f_1*f_2)=\zzd(f_1)*\zzd(f_2).$ Le r\'esultat de b) est alors une
cons\'equence du fait que, si $m\geq l$, alors $\zzdl$ est induite par
la restriction de $\zzd$ (remarque \ref{lcomutativ} plus haut). 
\ \\

\begin{lemme}\label{generatori} 
Pour tout entier $i$, $0\leq i\leq r$,  on pose
$A_i=diag(\pi_{D_F}^{a_1},\pi_{D_F}^{a_2}...\pi_{D_F}^{a_r})$ o\`u
pour $1\leq j\leq i$, $a_j=0$ et pour $i+1\leq j\leq r$, $a_j=1$. On
pose aussi
$A_{-1}=diag(\pi_{D_F}^{-1};\pi_{D_F}^{-1}...\pi_{D_F}^{-1})$. Soit
$s_F$ un syst\`eme de repr\'esentants de $GL_r(O_{D_F}/P_{D_F}^{ld})$
dans $GL_r(O_{D_F})$. Alors $\{h(x) : x\in s_F\cup
\{A_{-1};A_0;A_1...A_r\}\}$ est une famille g\'en\'eratrice de
$H(G'_F;K_F^l)$ comme $\mathbb{C}$-alg\`ebre.

\end{lemme}
\ \\
{\bf{D\'emonstration.}} On a d\'ej\`a vu que $\{h(BAC) : A\in
{\mathcal{A}}_F, B\in s_F,C\in s_F\}$ \'etait une famille
g\'en\'eratrice de $H(G'_F;K_F^l)$ comme ${\mathbb{C}}$-espace
vectoriel. Si $A\in {\mathcal{A}}_F$,
$A=diag(\pi_{D_F}^{a_1};\pi_{D_F}^{a_2}...\pi_{D_F}^{a_r})$ alors, si
$a_1\geq 0$, $A$ s'\'ecrit $A=A_0^{a_1}\Pi_{1\leq i\leq
r-1}A_i^{a_{i+1}-a_i}A_r$, et si $a_1<0$, $A$ s'\'ecrit
$A=A_{-1}^{-a_1}\Pi_{1\leq i\leq r-1}A_i^{a_{i+1}-a_i}$. Le lemme
\ref{volumes} implique alors le lemme \ref{generatori}.
\qed

\begin{lemme}\label{presfinie}  
L'alg\`ebre $H(G'_F;K_F^l)$ est de pr\'esentation finie.
\end{lemme}
\ \\
{\bf{D\'emonstration.}} D'apr\`es la remarque qui suit le corollaire
3.4. de [Be],  
 $H(G'_F;K_F^l)$ est un module de type fini sur son centre
 ${\mathcal{Z}}(G'_F;K_F^l)$, qui, \`a son tour, est une alg\`ebre de
type fini sur ${\mathbb{C}}$. Il existe donc un
 ${\mathcal{Z}}(G'_F;K_F^l)$-module libre $M$ de rang fini $p$ et un
 sous-module $N$ de $M$ tel que $H(G'_F;K_F^l)\simeq M/N$ en tant que
 ${\mathcal{Z}}(G'_F;K_F^l)$-modules. Soit $\{Y_1,Y_2,...Y_p\}$ une
 base de $M$ sur ${\mathcal{Z}}(G'_F;K_F^l)$. L'alg\`ebre commutaitive
 ${\mathcal{Z}}(G'_F;K_F^l)$, \'etant de type fini, est isomorphe \`a 
 ${\mathbb{C}}[X_1;X_2...X_n]/I$ o\`u $I$ est un id\'eal  donn\'e par
 un nombre fini de relations $R_1,R_2...R_u$ 
 entre les $X_i$. Elle est en particulier noeth\'erienne et donc le
 module $N$ est de type fini. Le ${\mathcal{Z}}(G'_F;K_F^l)$-module
 $M/N$ est donc le module engendr\'e par la famille
 $\{Y_1,Y_2...Y_p\}$ avec un nombre fini de relations
 $R'_1,R'_2...R'_v$ lin\'eaires entre les $Y_i$. En \'ecrivant encore
 pour tout $i,j\in \{1,2...p\}$ le produit $Y_iY_j$ sur la base
 $\{Y_1,Y_2...Y_p\}$ de $M$ on obtient encore une famille finie (de
 cardinal au plus $p^2$) de relations $\{R''_1,R''_2...R''_w\}$. Alors
 $H(G'_F;K_F^l)$ est isomorphe \`a l'alg\`ebre non commutative
 engendr\'ee sur $\mathbb{C}$ par les $n+p$ variables
 $X_1,X_2...X_n,Y_1,Y_2...Y_p$, avec les relations
 $R_1,R_2...R_u,R'_1,R'_2...R'_v,R''_1,R''_2,...R''_w$ et les
 $n(n-1)/2$ relations qui traduisent le fait que les variables
 $X_1,X_2...X_n$ commutent entre elles. Le lemme est d\'emontr\'e. 
\qed
\ \\
\ \\
{\bf Fin de la d\'emonstration du th\'eor\`eme} \ref{isomalg}{\bf.}
Nous explicitons la d\'emonstration dont le principe est d\^u \`a
Kazhdan~:  

Par le lemme \ref{presfinie}, $H(G'_F;K_F^l)$ est de pr\'esentation
finie ; c'est donc l'alg\`ebre non commutative engendr\'ee sur
$\mathbb{C}$ par un nombre fini de g\'en\'erateurs $g_1,g_2...g_n$
avec un nombre fini de relations $R_1,R_2...R_u$ qu'on va regarder
comme des polyn\^omes non commutatifs en $n$ variables qui s'annulent
en $(g_1;g_2...g_n)$. Par ailleurs, le lemme \ref{generatori} nous
fournit une famille finie $\{h_1,h_2...h_p\}$ de g\'en\'erateurs de
$H(G'_F;K_F^l)$. Soient alors $G_i(1\leq i\leq n)$ des polyn\^omes en $p$
variables qui appliqu\'es \`a $(h_1;h_2...h_p)$ donnent les 
$g_1,g_2...g_n$, et $F_i(1\leq i\leq p)$ des polyn\^omes en $n$
variables qui appliqu\'es \`a $(g_1;g_2...g_n)$ donnent les
$h_1,h_2...h_p$. Soit $s$ le plus grand nombre parmi les degr\'es
(totaux) des polyn\^omes   
$$R'_1=R_1\big((G_1;G_2...G_n)\big),R'_2=R_2\big((G_1;G_2...G_n)\big)\
.\ .\ .\ R'_u=R_u\big((G_1;G_2...G_n)\big)$$ 
et
$$F'_1=F_1\big((G_1;G_2...G_n)\big),F'_2=
F_2\big((G_1;G_2...G_n)\big)\ .\ .\ .\
F'_p=F_p\big((G_1;G_2...G_n)\big)$$ 
 et un compact {\bf{K}} suffisamment grand dans $G'_F$ pour qu'il
contienne tous les produits de $s$ \'el\'ements qui se trouvent dans
la r\'eunion des supports de tous les $h_i$. Soit ${\mathcal{C}}\in
{\mathcal{A}}_F$ de cardinal fini tel que ${\bf{K}}\subset
G'_F({\mathcal{C}})$ (voir le lemme \ref{cpct}). Prenons l'entier $m$
associ\'e \`a $\mathcal{C}$ comme dans le lemme \ref{cpct}. Notons
$h'_1,h'_2...h'_p$ les images de $h_1,h_2...h_p$ par
$\zzdl$. D\'efinissons un morphisme d'alg\`ebres
$t:{\mathbb{C}}(g_1,g_2...g_p)\rightarrow H(G'_L;K_F^l)$ d\'efini par
: 
$$t(g_i)=G_i\big((h'_1;h'_2...h'_p)\big)$$ 

Ce morphisme v\'erifie $t(R_i\big((g_1,g_2...g_n)\big))=0$ pour tout
$1\leq i\leq u$ par le choix de $m$ relativement aux polyn\^omes
$R'_1,R'_2...R'_u$ et par le lemme \ref{cpct}. Il induit donc un
morphisme d'alg\`ebres  
$$\bar{t}:H(G'_F;K_F^l)\rightarrow H(G'_L;K_L^l),$$
car $H(G'_F;K_F^l)$ est la $\mathbb{C}$-alg\`ebre non commutative
engendr\'ee par $g_1,g_2...g_n$ avec les relations traduites par
l'annulation des polyn\^omes $R_1,R_2...R_u$ en $(g_1;g_2...g_n)$. Or,
ce morphisme d'alg\`ebres v\'erifie  

\begin{equation}\label{***}
\bar{t}(h_i)=h'_i
\end{equation} 
pour tout $1\leq i\leq p$ par le choix de $m$ relatif aux polyn\^omes
$F'_1,F'_2...F'_p$ et par le lemme \ref{cpct}. Comme c'est un
morphisme d'alg\`ebres, \ref{***} et les lemmes \ref{volumes} et
\ref{generatori} impliquent que pour tout $A\in {\mathcal{A}}_F$, pour
tout $(\tilde{B};\tilde{C})\in \TT$ on a 
$$\bar{t}({\bf{1}}_{K_F^l\tilde{B}A\tilde{C}^{-1}K_F^l})=
{\bf{1}}_{K_L^l\bigl(\zzdl(\tilde{B})\bigl)
\bigl(\zzzdl(A)\bigl)\bigl(\zzdl(\tilde{C})\bigl)^{-1}K_L^l}.$$    
Comme c'est un morphisme d'espaces vectoriels qui co\"\i ncide avec
$\zzd$ sur une base de $H(G'_F;K_F^l)$, on a $\bar{t}=\zzd$, donc
$\zzd$ est un isomorphisme d'alg\`ebres.
\qed 
\\

Nous rappelons que le niveau d'une repr\'esentation lisse
irr\'eductible $\pi$ de $G'_F$ (resp. $G'_L$) est le plus petit entier
$l$ tel que $\pi$ ait un vecteur fixe non nul sous $K_F^l$ (resp. $K_L^l$).
Soit $\pi$ une \rli\  de niveau inf\'erieur
ou \'egal \`a $l$ de 
$G'_F$ et notons $V_{\pi}$ l'espace de la repr\'esentation
$\pi$. Alors $H(G'_F;K_F^l)$ agit sur $V_{\pi}^{K_F^l}$, espace des
vecteurs fixes sous $K_F^l$, par   
$$f(v)=\pi(f)v$$
pour tout $f\in  H(G'_F;K_F^l)$ et tout $v\in
V_{\pi}^{K_F^l}$. L'espace $V_{\pi}^{K_F^l}$ est  ainsi muni d'une
structure de $H(G'_F;K_F^l)$-module. On note $V_{\pi,H}$ le
$H(G'_F;K_F^l)$-module $V_{\pi}^{K_F^l}$ pour le diff\'erencier du
$\mathbb C$-espace $V_{\pi}^{K_F^l}$ avec lequel il co\"\i ncide
ensemblistement. On sait que $V_{\pi,H}$ est un $H(G'_F;K_F^l)$-module
non nul irr\'eductible et que $\pi\mapsto V_{\pi,H}$ induit une
bijection entre l'ensemble des classes  d'\'equivalence  de
repr\'esentations irr\'eductibles de $G'_F$ de niveau inf\'erieur ou
\'egal \`a $l$ et 
l'ensemble des classes d'isomorphie de $H(G'_F;K_F^l)$-modules
irr\'eductibles ([Be] ou [Ca2]). Maintenant, si $m$ est comme dans le th.
\ref{isomalg}, l'isomorphisme $\zzdl :
H(G'_F;K_F^l)\simeq H(G'_L;K_L^l)$ induit une bijection (not\'e
toujours $\zzdl$) entre l'ensemble des classes d'isomorphie de
$H(G'_F;K_F^l)$-modules irr\'eductibles et l'ensemble des classes
d'isomorphie de $H(G'_L;K_L^l)$-modules irr\'eductibles. Donc l'image
par $\zzdl$ de la classe d'isomorphie de $V_{\pi,H}$ est une classe
d'isomorphie de $H(G'_L;K_L^l)$-modules irr\'eductibles et elle
correspond \`a une classe d'\'equivalence $C_L$ de repr\'esentations
irr\'eductibles de niveau inf\'erieur ou \'egal \`a $l$ de $G'_L$. Si
$C_F$ est la classe 
d'\'equivalence de la repr\'esentation $\pi$, on pose
$\zzdl(C_F)=C_L$.

\begin{theo}\label{careseridica} 
a) L'application $\zzdl$ r\'ealise une bijection de l'ensemble des
classes d'\'equivalence  des \rlis\ de 
$G'_F$ de niveau inf\'erieur ou \'egal \`a $l$ (resp. \'egal \`a $l$)
sur l'ensemble des classes 
d'\'equivalence  des \rlis\ de $G'_L$ de  niveau inf\'erieur ou \'egal
\`a $l$ (resp. \'egal \`a $l$).     

b) Soit $\pi$ une \rli\ de $G'_F$ de niveau inf\'erieur ou \'egal \`a
$l$. Alors $\pi$ est 
de carr\'e int\'egrable si et seulement si $\zzdl(\pi)$
est une \care\ de $G'_L$.   

c) Soit $\pi$ une \rli\ de $G'_F$ de niveau inf\'erieur ou \'egal \`a
 $l$. Alors $\pi$ est 
 cuspidale si et seulement si $\zzdl(\pi)$
est une \cusp\ de $G'_L$.  
\end{theo} 

\begin{rem}
{\rm
On parle abusivement de $\zzdl(\pi)$ alors que $\zzdl$
n'est d\'efinie que pour les classes d'\'equivalence. C'est ce qui
arrivera parfois aussi par la suite, puisque toutes les propri\'et\'es
des repr\'esentations sont en fait des propri\'et\'es des classes
d'\'equivalence.} 
\end{rem}
\ \\
{\bf D\'emonstration.} a) Le fait que $\zzdl$ soit une bijection de
l'ensemble des 
classes d'\'equivalence  des \rlis\ de 
$G'_F$ de niveau inf\'erieur ou \'egal \`a $l$ 
sur l'ensemble des classes 
d'\'equivalence  des \rlis\ de $G'_L$ de  niveau inf\'erieur ou \'egal
\`a $l$ r\'esulte de la discussion faite avant l'\'enonc\'e du
th\'eor\`eme.  

Soit maintenant $\Pi(G'_F,l)$ (resp. $\Pi(G'_F,l-1)$) l'ensemble des classes
d'\'equivalence de repr\'esentations lisses irr\'eductibles de niveau 
inf\'erieur ou \'egal \`a $l$ (resp. \`a $l-1$) de $G'_F$. Adoptons
les m\^emes notations pour $G'_L$. Alors $\zzdl$
r\'ealise une bijection de $\Pi(G'_F,l)$ sur $\Pi(G'_L,l)$, et aussi
(mettre $l-1$ \`a la place de $l$), $\bar{\zeta}_{D_FD_L}^{l-1}$ r\'ealise
une bijection de $\Pi(G'_F,l-1)$ sur $\Pi(G'_L,l-1)$. Par la remarque
\ref{lcomutativ} et ce qu'on vient de voir, la restriction de $\zzdl$
\`a $\Pi(G'_F,l-1)$ induit $\bar{\zeta}_{D_FD_L}^{l-1}$. Nous obtenons
alors que $\zzdl$
r\'ealise une bijection de $\Pi(G'_F,l)\backslash \Pi(G'_F,l-1)$ sur
$\Pi(G'_L,l)\backslash \Pi(G'_L,l-1)$, qui est la variante
``niveau exactement \'egal \`a $l$'' de l'\'enonc\'e.

b) Supposons maintenant que $\pi$ soit de carr\'e int\'egrable.
Soit $\sigma$ une repr\'esentation (irr\'eductible et de niveau
inf\'erieur ou \'egal \`a $l$) 
se trouvant dans l'image par $\zzdl$ de la classe d'\'equivalence de
$\pi$. Comme nous l'avons dit plus haut, on a un isomorphisme d'espaces
vectoriels 
 
$$f:V_{\pi}^{K_F^l}\simeq V_{\sigma}^{K_L^l}$$ 
induit par $\zzdl$. L'isomorphisme $f$
induit un isomorphisme, dans les espaces duaux munis des actions
contragr\'edientes~: 

$$f':V_{\pi}^{'K_F^l}\simeq V_{\sigma}^{'K_L^l}.$$

Soient  $v\in V_{\pi}^{K_F^l}$ et $v'\in
V_{\pi}^{'K_F^l}$ tels que $v'(v)\neq 0$. Consid\'erons le coefficient
non nul de $\pi$

$$h_{\pi}:G'_F\to \mathbb C$$
d\'efini par 
$$g\mapsto v'\big( \pi(g)(v)\big).$$
Pour tout $g\in G'_F$, $h_{\pi}$ est constant sur $K_F^lgK_F^l$ \'egal
\`a $h_{\pi}(g)$. On a aussi
\begin{equation}\label{prima}
h_{\pi}(g)=vol(K_F^lgK_F^l)^{-1}v'\big(\pi({\bf
 1}_{K_F^lgK_F^l})(v)\big).
\end{equation} 
La
repr\'esentation $\pi$ \'etant de carr\'e int\'egrable, $|h_{\pi}|^2$ est
trivial sur $Z$ et int\'egrable sur $G'_F/Z$. Exprimons cette
propri\'et\'e \`a partir de la d\'ecomposition
 
$$G'_F=\coprod_{A\in {\mathcal{A}}_F}\coprod_{(\tilde{B};\tilde{C})\in
\TT}K_F^l\tilde{B}A\tilde{C}^{-1}K_F^l$$
(lemme \ref{3.2.6}.b)) en \'etudiant l'action par multiplication de $Z$
l\`a-dessus. Notons 
${\mathcal{A}}_F^0$ le sous-ensemble de ${\mathcal{A}}_F$ form\'e des
matrices $A=diag(\pi_{D_F}^{a_1};\pi_{D_F}^{a_2};...\pi_{D_F}^{a_r})$
telles que $a_1\in \{0;1\ ...d-1\}$. Pour toute matrice $A\in
{\mathcal{A}}_F^0$ notons ${\mathcal A}_F(A)$ l'ensemble des matrices
obtenues \`a partir de $A$ par multiplication avec une puissance de
$\pi_F=\pi_{D_F}^d$. C'est un sous-ensemble de ${\mathcal{A}}_F$. On a 

$${\mathcal{A}}_F=\coprod_{A\in {\mathcal{A}}_F^0}{\mathcal A}_F(A).$$
Montrons que, pour tout $A\in {\mathcal{A}}_F^0$, l'ensemble
$\coprod_{A'\in {\mathcal A}_F(A)}\coprod_{(\tilde{B};\tilde{C})\in
\TT}K_F^l\tilde{B}A'\tilde{C}^{-1}K_F^l$ est stable sous l'action de
$Z$ par multiplication et \'etudions cette action. Identifions $Z$
avec $F^*$. Alors, si $z$ est un \'el\'ement de $Z$, $z$ s'\'ecrit de
fa\c{c}on unique $z=\pi_F^{\alpha}x$ o\`u $x$ est un \'el\'ement de
$O_F^*$. D'autre part, on a un isomorphisme 

$$O_F^*/(1_F+P_F^l)\simeq (O_F/P_F^l)^*.$$
On peut ainsi d\'ecomposer l'action de $Z$ sur $\coprod_{A'\in
{\mathcal A}_F(A)}\coprod_{(\tilde{B};\tilde{C})\in
\TT}K_F^l\tilde{B}A'\tilde{C}^{-1}K_F^l$ puisque~:

- si $z=\pi_F^{\alpha}$, alors
  $zK_F^l\tilde{B}A'\tilde{C}^{-1}K_F^l$ s'\'ecrit
  $K_F^l\tilde{B}A''\tilde{C}^{-1}K_F^l$  
  o\`u $A''=zA'\in {\mathcal A}_F(A)$, et une simple v\'erification
  montre qu'on a aussi $\tilde{\bf T}_{F,l,A''}=\tilde{\bf T}_{F,l,A'}$, 

- si $z\in O_F^*$, alors~:

\ \ \ \ \ - si $z\in 1_F+P_F^l$, on a
$zK_F^l\tilde{B}A'\tilde{C}^{-1}K_F^l=K_F^l\tilde{B}A'\tilde{C}^{-1}K_F^l$,
tandis que  

\ \ \ \ - si $z\notin 1_F+P_F^l$,
$zK_F^l\tilde{B}A'\tilde{C}^{-1}K_F^l=K_F^l\tilde{B}'A'\tilde{C}^{'-1}K_F^l$,
o\`u $(\tilde{B}';\tilde{C}')$ est un \'el\'ement de $\TT$ diff\'erent
de $(\tilde{B};\tilde{C})$ et qui ne d\'epend que de la classe de $z$
modulo $1_F+P_F^l$. 

Finalement, dire que $|h_{\pi}|^2$ est int\'egrable sur $G'_F/Z$
revient \`a dire que la somme 
 
\begin{equation}\label{suma}
\sum_{A\in {\mathcal{A}}^0_F}\Big(\sum_{(\tilde{B};\tilde{C})\in
\TT}\big(card((O_F/P_F^l)^*)\big)^{-1} 
\end{equation}
$$
vol(K_F^l\tilde{B}A\tilde{C}^{-1}K_F^l;dg)\big(vol(1_F+P_F^l;dz)\big)^{-1}|h_{\pi}(\tilde{B}A\tilde{C}^{-1})|^2\Big)    
$$
est convergente.

Maintenant, $f(v)$ est un \'el\'ement de $V_{\sigma}^{K_L^l}$,
$f'(v')$ est un \'el\'ement de 
$V_{\sigma}^{'K_L^l}$ et l'application  
$$h_{\sigma}:G'_L\to \mathbb C$$
d\'efinie par 
$$g\mapsto f'(v')\big( \sigma(g)(f(v))\big)$$
 est un coefficient non nul de $\sigma$.
Pour tout $g\in G'_L$, $h_{\sigma}$ est constant sur $K_L^lgK_L^l$ \'egal
\`a $h_{\sigma}(g)$. On a aussi
\begin{equation}\label{adoua}
h_{\sigma}(g)=vol(K_L^lgK_L^l)^{-1}f'(v')\big(\sigma({\bf
 1}_{K_L^lgK_L^l})(f(v))\big).
\end{equation}

La fonction
$h_{\sigma}$ est de carr\'e int\'egrable modulo le centre sur $G'_L$
si et seulement si la somme  
\begin{equation}
 \sum_{A\in {\mathcal{A}}^0_F}\Big(\sum_{(\tilde{B};\tilde{C})\in
\TTl}\big(card((O_L/P_L^l)^*)\big)^{-1}
\end{equation} 
$$vol(K_L^l\tilde{B}\zzzdl(A)\tilde{C}^{-1}K_L^l;dg)\big(vol(1_L+P_L^l;dz)\big)
^{-1}|h_{\sigma}(\tilde{B}\zzzdl(A)\tilde{C}^{-1})|^2\Big)$$    
est convergente, o\`u on a tenu compte du fait que $\zzzdl$ r\'ealise
une bijection de ${\mathcal A}_F^0$ sur ${\mathcal A}_L^0$. Mais cette
somme correspond terme pour terme \`a la somme \ref{suma} puisque 

- les volumes  des sous-ensembles de $G'_F$ et $G'_L$ qui se
  correspondent sont \'egaux pour les mesures fix\'ees sur $G'_F$ et
  $G'_L$,  

- les volumes  des sous-ensembles des centres de $G'_F$ et $G'_L$ qui
  se correspondent sont \'egaux pour les mesures fix\'ees sur les
  centres $F^*$ de $G'_F$ et $L^*$ de $G'_L$, 

- $card((O_F/P_F^l)^*)=card((O_L/P_L^l)^*)$ parce que les anneaux
  $O_F/P_F^l$ et $O_L/P_L^l$ sont isomorphes, 

- par construction de l'isomorphisme $f'$ on a 
$$v'\big(\pi({\bf
  1}_{K_F^lgK_F^l})(v)\big)=f'(v')\big(\sigma({\bf 
  1}_{K_L^lgK_L^l})(f(v))\big),$$
 et alors \ref{prima} et \ref{adoua} impliquent que
$$h_{\sigma}(\tilde{B} A \tilde{C}^{-1})=
h_{\pi}(\zzdl(\tilde{B}) \zzzdl(A)   \zzdl(\tilde{C}^{-1})).$$   

On vient de montrer qu'un coefficient non nul de $\sigma$ est de
carr\'e int\'egrable sur $G'_L/Z$. Donc $\sigma$ est de carr\'e
int\'egrable.
La r\'eciproque se montre exactement de la m\^eme fa\c{c}on.

c) Une repr\'esentation est cuspidale si et seulement si elle admet un
coefficient non nul  \`a support compact modulo le centre. La d\'emonstration
est la m\^eme que celle du point b), en plus facile~: il faut remplacer
``somme convergente'' avec ``somme \`a nombre fini de termes non nuls''.
\qed
\ \\
\ \\

Soit $\psi$ un caract\`ere additif non trivial de $F$. On dit que $k$ est le
conducteur de $\psi$ si $k$ est le plus petit entier tel que $P_F^k$
soit inclus dans $ker\, \psi$. On fixe une fois pour toutes un
caract\`ere $\psi$ non trivial de $F$ {\it 
de conducteur nul}. Si $\pi$ est une repr\'esentation lisse
irr\'eductible de $G'_F$, on note $L(s;\pi)$, $\epsilon(s;\pi;\psi)$
et $\epsilon'(s;\pi;\psi)$ les fonctions de Godement-Jacquet ([GJ]). 
On sait que la fonction $\epsilon(s;\pi;\psi)$ est \'egale \`a
$q^{-ms}$, \`a multiplication par un scalaire non nul pr\`es,
o\`u $q$ est, comme avant, le cardinal du corps r\'esiduel de $F$, et
$m$ est un entier 
([GJ], th.3.3, (4)). D'apr\`es [GJ], \'equation (3.3.5), page 33,
l'entier $m$ ne d\'epend pas du choix de $\psi$ plus haut (une fois
que son conducteur est nul)~; il ne
d\'epend donc que de $\pi$ : on le note par la suite $m(\pi)$ et on
l'appelle le {\it conducteur de $\pi$}. En suivant [GJ], nous
allons utiliser une formule particuli\`ere pour la fonction
$\epsilon'(s;\pi;\psi)$, tr\`es commode pour les calculs. Supposons
que le niveau de $\pi$ soit $l$. Soit $f$ un 
coefficient de $\pi$ tel que $f$ est constant non nul sur $K_F^l$
(choisi comme dans la d\'emonstration du th. \ref{careseridica}
b)). Prenons le cas 
particulier 
$\Phi={\bf 1}_{K_F^l}$ o\`u dans le th.
3.3 de [GJ]. Nous posons donc
$$Z(s;f)=\int_{K_F^l}f(g)|N(g)|^sdg=f(1)vol(K_F^l)$$
o\`u $N$ est la norme r\'eduite. On  sait par le dit th\'eor\`eme
qu'il existe $s_0$ dans $\Bbb R$ tel que, si $s$ est 
un nombre complexe de partie r\'eelle sup\'erieure ou \'egale \`a
$s_0$, alors 
$$Z(s;\check{f})=q^{-nl}\int_{\pi_F^{-l}M_r(O_{D_F})\cap
GL_r(D_F)}\psi(tr_{M_r(D_F)/F}(g))f(g^{-1})|N(g)|^sdg$$  
converge et $Z(s;\check{f})$
est une fraction rationnelles en la quantit\'e $q^{-s}$. Par les points
(2) et (4) du m\^eme th\'eor\`eme, $\epsilon'(s;\pi;\psi)$ est une
fraction rationnelle en la quantit\'e $q^{-s}$ qui v\'erifie~:   

\begin{equation}\label{defepsilonprim}
\epsilon'(s;\pi;\psi)=(-1)^{r(d-1)}Z(1-s+(n-1)/2;\check{f})Z(s+(n-1)/2;f)^{-1}.
\end{equation}
\\

Nous adoptons les m\^emes conventions pour le corps $L$. Si $L$ est
$m$-proche de $F$, le triplet de $m$-proximit\'e
$(\pi_F;\pi_L;\lll)$ associ\'e induit
naturellement un
isomorphisme de groupes additifs $\lambda_{-m} : \pi_F^{-m}O_F/O_F\simeq
\pi_L^{-m}O_L/O_L$. Alors, si $\psi_L$ est un caract\`ere additif de
$L$, on dit que $\psi_L$ est $m${\it -proche} de $\psi$ si $\psi_L$ est de
conducteur nul et le caract\`ere induit par $\psi$ sur
$\pi_F^{-m}O_F/O_F$ et le caract\`ere induit par $\psi_L$ sur
$\pi_L^{-m}O_L/O_L$ se correspondent via $\lambda_{-m}$
(i.e. $\psi=\psi_L\circ \lambda_{-m}$).\\

Le th\'eor\`eme suivant est d'une importance capitale pour la
d\'emonstration de la correspondance~:

\begin{theo}\label{epsilonseridica}
Si $m$ est comme dans le th\'eor\`eme \ref{isomalg}, si $L$ est
un corps $m$-proche de $F$, si $\psi_L$ est un caract\`ere additif de $L$
$m$-proche de $\psi$, alors on a~:

a) si $\pi$ est une \rli\ de niveau $\leq l$ de $G'_F$, alors les
fonctions $\epsilon'(s;\pi;\psi)$  et $\epsilon'(s;\zzdl(\pi);\psi_L)$
sont \'egales, 

b) si $\pi$ est une \cusp\ de niveau $\leq l$ de $G'_F$, alors les
fonctions $\epsilon(s;\pi;\psi)$ et $\epsilon(s;\zzdl(\pi);\psi_L)$
sont \'egales ; en particulier $m(\pi)=m(\zzdl(\pi))$. 
\end{theo}
\ \\
{\bf D\'emonstration.} On peut
supposer que le niveau de $\pi$ est $l$.
 
a) Posons $\pi_L=\zzdl(\pi)$. Soient $f$ un coefficient de $\pi$ et
$f_L$ un coefficient de $\pi_L$ choisis comme dans la d\'emonstration
de la prop.\ref{careseridica}. Notons
$N_L$ la norme r\'eduite sur $G'_L$. Montrons qu'on a alors~: 

(i) $Z(s;f)=Z(s;f_L)$ et 

(ii) pour les $s$ pour lesquels $Z(s;\check{f})$ converge,
$Z(s;\check{f}_L)$ converge et on a $Z(s;\check{f})=Z(s;\check{f}_L)$.\\
\ \\
Le point (i) est \'evident. Pour montrer (ii) on prouve que 
\begin{equation}\label{Y}
\int_{\pi_L^{-l}M_r(O_{D_L})\cap
GL_r(D_L)}\psi_L(tr_{M_r(D_L)/L}(g))f_L(g^{-1})|N_L(g)|^sdg=\cr 
=\int_{\pi_F^{-l}M_r(O_{D_F})\cap
GL_r(D_F)}\psi(tr_{M_r(D_F)/F}(g))f(g^{-1})|N(g)|^sdg, 
\end{equation}
en montrant que les int\'egrales prennent la forme de deux
sommes identiques. On a une d\'ecomposition 

$$
\pi_F^{-l}M_r(O_{D_F})\cap GL_r(D_F)=\coprod_{A\in {\mathcal A}_F^+}
\coprod_{(\tilde{B};\tilde{C}) \in
\tilde{T}_{l,A,F}}K_L^l\tilde{B}A\tilde{C}^{-1}K_L^l
$$
o\`u ${\mathcal A}_F^+$ est l'ensemble des matrices dans ${\mathcal
A}_F$ dans lesquelles toutes les puissances de l'uniformisante
$\pi_{D_F}$ qui apparaissent sont sup\'erieures ou \'egales \`a $-ld$
(rappelons que $\pi_{D_F}^{d}=\pi_F$). On a une d\'ecomposition
analogue pour $\pi_L^{-l}M_r(O_{D_L})\cap GL_r(D_L)$ qu'on peut
\'ecrire 

$$
\pi_L^{-l}M_r(O_{D_L})\cap GL_r(D_L)=\coprod_{A\in {\mathcal A}_F^+}
\coprod_{(\tilde{B};\tilde{C}) \in
\tilde{T}_{l,A,F}}K_L^l\zzdl(\tilde{B})\zeta_{D_FD_L}^m(A)\zzdl
(\tilde{C}^{-1})K_L^l   
$$
tenant compte du fait que $\zeta_{D_FD_L}^m$ r\'ealise une bijection
de ${\mathcal A}_F^+$ sur ${\mathcal A}_L^+$ et que, pour tout $A\in
{\mathcal A}_F^+$, $\zzdl$ r\'ealise une bijection de
$\tilde{T}_{l,A,F}$ sur $\tilde{T}_{l,\zeta_{D_FD_L}^m(A),L}$. 

On a vu \`a la proposition \ref{careseridica} que, par la
construction de $f_L$ \`a partir de $f$, on a~: pour tout $A\in
{\mathcal A}_F$, pour tout $(\tilde{B};\tilde{C}) \in
\tilde{T}_{l,A,F}$, $f$ est  constant sur l'ensemble
$K_F^l\tilde{B}A\tilde{C}^{-1}K_F^l$ et $f_L$ est constant sur
$K_L^l\zzdl(\tilde{B})\zeta_{D_FD_L}^m(A)\zzdl(\tilde{C}^{-1})K_F^l$;
en outre les valeurs de $f$ et $f_L$ sont \'egales. C'est pareil pour les
fonctions $|N|$ et $|N_L|$ qui y sont constantes \'egales \`a
$|N(A)|$. On a aussi  

$$
vol(K_F^l\tilde{B}A\tilde{C}^{-1}K_F^l)=vol(K_L^l\zzdl(\tilde{B})\zeta_{D_FD_L}^m(A)\zzdl(\tilde{C}^{-1})K_L^l).
$$ 
Par ailleurs, si on note $U$ la matrice $(u_{ij})_{1\leq i,j\leq r}\in
GL_r(D_F)$ d\'efinie par $u_{ij}=\delta_{i,r-j}$, alors pour tout $A\in
{\mathcal A}_F$ et pour tout $(\tilde{B};\tilde{C}) \in
\tilde{T}_{l,A,F}$, la fonction $f(g^{-1})$  est  constante sur
l'ensemble  

$$
K_F^l\tilde{C}A^{-1}\tilde{B}^{-1}K_F^l
$$
\'egale \`a $f(\tilde{B}A\tilde{C}^{-1})$, et on a 

$$
K_F^l\tilde{C}A^{-1}\tilde{B}^{-1}K_F^l=
K_F^l(\tilde{C}U)(UA^{-1}U)(U\tilde{C}^{-1})K_F^l
$$ 
o\`u $UA^{-1}U\in {\mathcal A}_F$ et $(\tilde{C}U;U\tilde{B}^{-1}) \in
\tilde{T}_{l,UA^{-1}U,F}$.  Le m\^eme ph\'enom\`ene se produit aussi
sur $G'_L$ et toutes les applications des objets sur $F$ vers les
objets correspondants sur $L$ commutent \`a l'action par
multiplication de $U$, qui est une simple permutation de lignes ou
colonnes. 

Montrons encore que, si $A\in {\mathcal A}_F^+$, si $(\tilde{B};\tilde{C}) \in \tilde{T}_{l,A,F}$, alors  

-- $\psi\circ tr_{M_r(D_F)/F}$ est constante sur $K_F^l\tilde{B}A\tilde{C}^{-1}K_F^l$, 

-- $\psi_L\circ tr_{M_r(D_L)/L}$ est constante sur 
$K_L^l\zzdl(\tilde{B})\zeta_{D_FD_L}^m(A)\zzdl(\tilde{C}^{-1})K_L^l$, et 

-- on a 
$$
\psi\circ tr_{M_r(D_F)/F}(\tilde{B}A\tilde{C}^{-1})=\psi_L\circ
tr_{M_r(D_L)/L}(\zzdl(\tilde{B})\zeta_{D_FD_L}^m(A)\zzdl(\tilde{C}^{-1})).
$$

Posons $A=diag(\pi_{D_F}^{a_1},\pi_{D_F}^{a_2}...\pi_{D_F}^{a_r})$
o\`u $a_1\leq a_2...\leq a_r$. On rappelle qu'on a $A\in {\mathcal
A}_F^+$ et donc $a_1\geq -ld$. Par cons\'equent, la diff\'erence de
deux \'el\'ements de l'ensemble $K_F^l\tilde{B}A\tilde{C}^{-1}K_F^l$
est un \'el\'ement de $M_r(O_{D_F})$, donc la diff\'erence de leurs
traces r\'eduites est un \'el\'ement de $O_F$ et $\psi$ est trivial
sur $O_F$. Donc $\psi\circ tr_{M_r(D_F)/F}$ est constante sur
$K_F^l\tilde{B}A\tilde{C}^{-1}K_F^l$. De la m\^eme fa\c{c}on on montre
que $\psi_L\circ tr_{M_r(D_L)/L}$ est constante sur  
$K_L^l\zzdl(\tilde{B})\zeta_{D_FD_L}^m(A)\zzdl(\tilde{C}^{-1})K_L^l$.

Maintenant, montrer qu'on a 

$$
\psi\circ tr_{M_r(D_F)/F}(\tilde{B}A\tilde{C}^{-1})=\psi_L\circ
tr_{M_r(D_F)/F}(\zzdl(\tilde{B})\zeta_{D_FD_L}^m(A)\zzdl(\tilde{C}^{-1})),
$$
revient \`a  montrer que 

$$\psi(\pi_F^{-l} tr_{M_r(D_F)/F}(\tilde{B}(\pi_F^lA)\tilde{C}^{-1})=$$
$$=\psi_L(\pi_L^{-l}
tr_{M_r(D_L)/L}(\zzdl(\tilde{B})(\pi_L^{l}\zeta_{D_FD_L}^m(A))\zzdl(\tilde{C}^{-1}))$$ 
et, vu que $\psi_L$ et  $\psi$ sont $m$-proches, et $m\geq l$, il suffit de
montrer que l'image (bien d\'efinie) de\\
$$tr_{M_r(D_F)/F}(\tilde{B}(\pi_F^lA)\tilde{C}^{-1})$$
dans $O_F/P_F^l$ et l'image (bien d\'efinie) de 
$$tr_{M_r(D_L)/L}(\zzdl(\tilde{B})(\pi_L^{l}\zeta_{D_FD_L}^m(A))\zzdl(\tilde{C}^{-1}))$$  
dans $O_L/P_L^l$ se correspondent par l'application
$\bar{\lambda}_{FL}^l$ (induite par $\bar{\lambda}_{FL}^m$).  

Maintenant, si on choisit des repr\'esentatnts $B$ et $C^{-1}$ de
$\tilde{B}$ et de $\tilde{C}^{-1}$, puisque $\bar{\lambda}_{D_FD_L}^l$
(induit par $\bar{\lambda}_{FL}^l$) est un isomorphisme d'anneaux de
$O_{D_F}/P_{D_F}^{dl}$ sur $O_{D_L}/P_{D_L}^{dl}$, on
a 
$\zeta_{D_FD_L}^m(B\pi_F^lAC^{-1})-\zeta_{D_FD_L}^m(B)\zeta_{D_FD_L}^m(\pi_F^lA)\zeta_{D_FD_L}^m(C^{-1})\in
M_r(P_{D_L}^{ld})$ et par cons\'equent l'image dans $O_L/P_L^l$ de  
$tr_{M_r(D_L)/L}(\zeta_{D_FD_L}^m(B\pi_F^lAC^{-1}))$  est \'egale \`a
l'image dans $O_L/P_L^l$ de
$tr_{M_r(D_L)/L}(\zeta_{D_FD_L}^m(B)\zeta_{D_FD_L}^m(\pi_F^lA)\zeta_{D_FD_L}^m(C^{-1}))$.  
Il nous suffit donc de montrer que l'image dans $O_F/P_F^l$  de
$tr_{M_r(D_F)/F}(B(\pi_F^lA)C^{-1})$  et l'image dans $O_L/P_L^l$ de
$tr_{M_r(D_L)/L}(\zeta_{D_FD_L}^m(B\pi_F^lAC^{-1}))$  se correspondent
par l'application $\bar{\lambda}_{FL}^l$. Mais la trace r\'eduite d'un
\'el\'ement de $M_r(D_F)$ est la somme des traces r\'eduites des
\'el\'ements diagonaux de cette matrice et de m\^eme pour $D_L$, donc,
pour avoir (enfin) le r\'esultat voulu il nous suffit du lemme
suivant~:

\begin{lemme} Soit $x\in O_{D_F}$. Alors l'image de $tr_{D_F/F}(x)$
dans $O_F/P_F^l$ et l'image de $tr_{D_L/L}(\lambda_{D_FD_L}^m(x))$
dans $O_L/P_L^l$ se correspondent par l'isomorphisme
$\bar{\lambda}_{FL}^l$. 
\end{lemme}
\ \\
{\bf D\'emonstration.} Reprenons les notations de la section
\ref{bijform}. \'Ecrivons 
$$x=\sum_{i=0}^{d-1}\pi_{D_F}^ie_i$$
o\`u tous les $e_i$ sont dans $E$, et m\^eme dans $O_E$ puisque $x\in O_{D_F}$.
Alors, par d\'efinition,
  
$$
\lambda_{D_FD_L}^m(x)=\sum_{i=0}^{d-1}\pi_{D_F}^i\lambda_{EK}^m(e_i).
$$
L'alg\`ebre $D_F$ agit sur le $E$-espace vectoriel $D_F$ de dimension $d$
par multiplication \`a gauche et la trace r\'eduite de $x$ sur $F$ est
la trace de l'endomorphisme qui correspond \`a $x$ ; on peut la
calculer facilement en choisissant dans le $E$-espace $D_F$ la base
$\pi_{D_F}^0,\pi_{D_F}^1,... \pi_{D_F}^{d-1}$. On
trouve

$$
tr_{D_F/F}(x)=\sum_{i=0}^{d-1}\sigma_F^i(e_0).
$$ 
Pareillement, on a 

$$
tr_{D_L/L}(\lambda_{D_FD_L}^m(x))=\sum_{i=0}^{d-1}\sigma_L^i(\lambda_{EK}^m(e_0))=\sum_{i=0}^{d-1}\lambda_{EK}^m(\sigma_F^i(e_0)).  
$$
Pour tout $1\leq i \leq d-1$, l'image de $\sigma_F^i(e_0)$ dans
$O_E/P_E^l$ et l'image de $\lambda_{EK}^m(\sigma_F^i(e_0))$ dans
$O_K/P_K^l$ se correspondent par l'isomorphisme $
\bar{\lambda}_{EK}^m$, et donc l'image de $tr_{D_F/F}(x)$ 
dans $O_E/P_E^l$ et l'image de $tr_{D_L/L}(\lambda_{D_FD_L}^m(x))$ 
dans $O_K/P_K^l$ se correspondent aussi par l'isomorphisme $
\bar{\lambda}_{EK}^m$. Comme $tr_{D_F/F}(x)\in O_F$ et
$tr_{D_L/L}(\lambda_{D_FD_L}^m(x))\in O_L$ et que $
\bar{\lambda}_{EK}^m:O_E/P_E^l\simeq O_K/P_K^l$ induit $
\bar{\lambda}_{FL}^m:O_F/P_F^l\simeq O_L/P_L^l$, le lemme est
d\'emontr\'e.
\qed
\ \\

On peut donc d\'ecomposer les int\'egrales dans l'\'egalit\'e \ref{Y}
en des sommes qui se correspondent terme \`a terme et en d\'eduire que   
$$\int_{\pi_L^{-l}M_r(O_{D_L})\cap
GL_r(D_L)}\psi(tr_{M_r(D_L)/L}(g))f_L(g^{-1})|N_L(g)|^sdg=$$ 
$$=\int_{\pi_F^{-l}M_r(O_{D_F})\cap
GL_r(D_F)}\psi_L(tr_{M_r(D_F)/F}(g))f(g^{-1})|N(g)|^sdg.$$ 
\ \\

On conclut maintenant par la relation \ref{defepsilonprim} et par
l'\'egalit\'e des fonctions $Z$ (pour un nombre infini de valeurs de
$q^{-s}$) montr\'ee.\\

b) Les facteurs $\epsilon$ et $\epsilon'$ sont reli\'es par la
relation-d\'efinition~:

\begin{equation}\label{defepsilon}
\epsilon'(s;\pi;\psi)=
\epsilon(s;\pi;\psi)L(1-s;\check{\pi})L(s;\pi)^{-1}.
\end{equation}
o\`u $\check{\pi}$ est la repr\'esentation contragr\'ediente de $\pi$.\\

Si $G'_F$ n'est pas le groupe des \'el\'ements inversibles d'une
alg\`ebre \`a division, alors la fonction $L$ associ\'ee \`a une
\cusp\ est triviale. Comme la contragr\'ediente d'une \cusp\ est une
\cusp\, on conclut par le point a) ci-dessus et le th\'eor\`eme
\ref{careseridica}.c).\\

Si $G'_F$ est le groupe des \'el\'ements inversibles d'une
alg\`ebre \`a division, alors on sait que si $\pi$ est une \rep\ de
$G'_F$ ou $G'_L$ on a $m(\pi)=niv(\pi)+n-1$. Comme
$\zeta_{D_FD_L}^m$ conserve le niveau, nous savons donc d\'ej\`a qu'elle
conserve le conducteur. Autrement dit, le rapport   
\begin{equation}\label{const}
\epsilon(s;\pi_F;\psi)/\epsilon(s;\pi_L;\psi_L)=c
\end{equation}
o\`u $c$ est une constante. Reste \`a montrer que $c=1$.  

Si $\pi$ est une \rep\ de $G'_F$ ou $G'_L$ on a
$$L(1-s;\check{\pi})L(s;\pi)^{-1}=U_{\pi}(q^{-s})$$
o\`u 
\begin{equation}\label{XXX}
U_{\pi}(X)=1\ {\text {ou}}\ U_{\pi}(X)=\frac{X(1-\alpha X)}{X-\beta}
\end{equation} 
avec $\alpha$ et $\beta$ deux nombres complexes (facile a obtenir \`a
partir de [GJ], th.5.11). Les \'egalit\'es \ref{defepsilon}, \ref{const}
et le point a) impliquent alors~:
$$U_{\pi_L}/U_{\pi_F}=c.$$
Des calculs simples montrent que, quelle que soit la forme de
$U_{\pi_F}$ et $U_{\pi_L}$ 
(voir \ref{XXX}), on doit avoir $c=1$.
\qed

\newpage

\section{Cons\'equences imm\'ediates}

Dans cette section nous donnons une preuve du fait que les r\'esultats
de [Ba3], a priori valables en caract\'eristique nulle, sont vrais
ind\'ependamment de la caract\'eristique. Ce sont les th\'eor\`emes
\ref{niveau} et \ref{finitude} plus bas. Ce n'est pas la preuve la
plus simple, ni la plus naturelle qui soit, \'etant bas\'ee sur
la construction faite \`a la section pr\'ec\'edente. D'autre part, ces
r\'esultats sont un corollaire imm\'ediat de la dite construction et
sont, aussi, indispensables pour la preuve de
la correspondance de Jacquet-Langlands en caract\'eristique non
nulle. Ils \'etaient donc incontournables. 
 
Soit $F$ un corps local non archim\'edien, soit $D$ une alg\`ebre \`a
division centrale de dimension finie $d^2$ sur $F$, soit $r$ un entier
strictement positif et posons $A=M_r(D)$. On pose $n=rd$. 

\begin{theo}\label{niveau} 
a) Si $\pi$ est une repr\'esentation cuspidale de $A^*$, alors on a
$$niv(\pi)\leq \frac{m(\pi)}{r}-d+2.$$

b) Si $\pi$ est une repr\'esentation lisse
irr\'eductible quelconque de $A^*$, alors on a 
$$niv(\pi)\leq m(\pi)-n+2r.$$
\end{theo}

Soit maintenant $F$ un corps global et soit $A$ une alg\`ebre centrale
simple de 
dimension finie sur $F$. Soient $F^*(\a)$ et $A^*(\a)$ les groupes des
ad\`eles de
$F^*$ et $A^*$ respectivement. On identifie $F^*$ au centre de $A^*$ et
$F^*(\a)$ au centre de  $A^*(\a)$. Pour toute place $v$ de $F$ on note
$F_v$ et $A_v$  
les localis\'es de $F$ et de $A$ en $v$.

\begin{theo}\label{finitude} 
Soit $V$ un ensemble fini de places finies de
$F$. Si pour toute place $v\notin V$ on fixe une repr\'esentation
lisse irr\'eductible $\pi_v$ de $A_v^*$,
alors il existe au plus un nombre fini de classes d'\'equivalence de
repr\'esentations automorphes cuspidales $\pi$ de  $A^*(\a)$,  telles que, pour
tout $v\notin V$, la composante locale de $\pi$ \`a la place $v$ soit
\'equivalente \`a $\pi_v$.
\end{theo}
\ \\

La preuve des deux th\'eor\`emes est tr\`es concise~: Dans [Ba3], nous
avons montr\'e que le th.\ref{niveau}.a) implique le
th.\ref{niveau}.b) qui \`a son tour implique le
th.\ref{finitude}, et ce ind\'ependamment de la caract\'eristique. D\`es
lors, le seul probl\`eme reste le point \ref{niveau}.a). Il a \'et\'e
prouv\'e dans [Ba3] en caract\'eristique nulle. Pour boucler la
d\'emonstration nous indiquons ici la preuve dans le cas o\`u $F$ est de
caract\'eristique non nulle. Si, dans ce cas, $\pi$ est une
\cusp\ de $A^*=GL_r(D)$ de niveau $l$, alors en se pla\c{c}ant dans la
situation des 
th\'eor\`emes \ref{careseridica} et \ref{epsilonseridica} et avec les
m\^emes notations, nous pouvons appliquer le th.\ref{niveau}.a) \`a la
repr\'esentation {\it cuspidale} (par \ref{careseridica}.c))
$\zzdl(\pi)$, puisque le corps L est de caract\'eristique 
nulle. On obtient ensuite la m\^eme relation pour $\pi$, car l'application
$\zzdl$ conserve toutes les quantit\'es qui apparaissent dans
l'in\'egalit\'e (par construction pour $r$ et $n$, par
th.\ref{careseridica}.a) pour le niveau et \ref{epsilonseridica}.b)
pour le conducteur). 

Le th\'eor\`eme \ref{finitude} appliqu\'e \`a une situation
particuli\`ere (la prop.\ref{finiG'}) jouera un r\^ole important dans la
d\'emonstration de la correspondance.

\newpage

\section{Quelques r\'esultats d'analyse harmonique sur des structures proches}

Soit $F$ un corps local non archim\'edien {\it de caract\'eristique
non nulle} et $D_F$ une alg\`ebre \`a division centrale sur $F$ de
dimension $d^2$. Soit $r\in \mathbb{N}^*$. On pose $n=dr$ et
$G_F=GL_{n}(F)$ et $G'_F=GL_r(D_F)$. Comme
en caract\'eristique nulle, on veut trouver une
correspondance entre les \ecis\ de $G_F$ et les \ecis\ de $G'_F$. Le
fait essentiel pour lequel la 
d\'emonstration ne marche pas comme en caract\'eristique nulle est
qu'on n'a pas l'orthogonalit\'e des caract\`eres sur $G'_F$. Je ne
pense pas qu'on puisse l'obtenir directement comme sur $G_F$, m\^eme
si on a construit dans la section pr\'ec\'edente une situation proche
en caract\'eristique nulle, parce qu'on ne  sait pas ``relever'' les
int\'egrales orbitales, et surtout parce qu'on n'a pas
l'int\'egrabilit\'e locale des caract\`eres pour les repr\'esentations
de $G'_F$. C'est pourquoi on va \'etudier plut\^ot $G_F$ et $G'_F$ en
parall\`ele. Plus pr\'ecisement on aura \`a tout instant en m\'emoire le
carr\'e~: 

$$\diagram
{G_L&\hfl{1}{}&G'_L\cr
 & & \cr
\vfl{2}{}&&\vfl{2'}{}\cr
 & & \cr
G_F&........&G'_F\cr}
$$ 
o\`u $L$ est un corps proche de $F$ qui est de caract\'eristique
nulle. La fl\`eche 1 est la correspondance (qu'on va noter  $\cf$)
d\'ej\`a \'etablie en caract\'eristique nulle, les fl\`eches 2 et 2'
sont les applications du type $\zzz$ et $\zzd$, et on voudrait
d\'efinir une correspondance \`a la place o\`u on a mis sur le dessin
des pointill\'es. Les ennuis viennent du fait que les
correspondances horizontales et verticales sont de natures tr\`es
diff\'erentes~: pour pouvoir user de 2 et 2' on doit \^etre s\^ur que
les objets qu'on veut transf\'erer sont constants sur des ouverts
assez ``gros'', alors que pour la fl\`eche 1 on ne sait pas en quelle
mesure elle conserve la propri\'et\'e ``\^etre constant sur un ouvert
assez gros'' (que ce soit pour des fonctions, pour leurs int\'egrales
orbitales ou pour les caract\`eres de repr\'esentations). C'est \`a
l'\'etude de ce probl\`eme que la pr\'esente section est d\'edi\'ee. 
Un autre probl\`eme  d\'ecoule du fait que les    
correspondances verticales 2 et 2' sont partielles et envoyent
{\it certaines} repr\'esentations des groupes d'en bas sur {\it certaines}
repr\'esentations des groupes d'en haut, et pour obtenir des
renseignements sur une autre repr\'esentation on 
est oblig\'e de changer de corps $L$. Pour r\'egler ce probl\`eme il
faut se placer dans une situation o\`u le niveau de toutes les
repr\'esentations qui apparaissent est born\'e uniform\'ement, et
c'est ce qu'on va faire dans la section 5 (suivante) au cours de la
d\'emonstration proprement dite. 

Dans les sous-sections 1, 2 et 3 nous nous int\'eressons
seulement \`a $GL_n(F)$. C'est dans la section 4 qu'on
d\'emontre le seul r\'esultat sur les formes int\'erieures de
$GL_n(F)$ dont nous avons besoin par la suite.\\ 

\subsection{\'El\'ements proches} Soient $F$ et $L$ deux corps locaux non archim\'ediens
$m$-proches. On reprend toutes les notations de la section
2, notamment $\llam$. On note $v_F$ et $v_L$ les
valuations sur $F$ et $L$ respectivement. 

Si $a\in F^*$, $b\in L$ et
$0<l\leq m$ on dit que $a$ et $b$ sont $l$-$proches$ si
$(b-\lam(a))\in P_L^{(l+v_F(a))}$. On note alors $a\sim_l b$. On
remarquera que $b$ est alors non nul, car il a la m\^eme valuation que
$a$. On consid\`ere que les \'el\'ements nuls de $F$ et $L$ sont
$l$-proches pour tout $l$. 

Soient $a\in F^*$ et $b\in L^*$. Si $a=\pi_F^{v_F(a)}a'$ avec 
$a'\in O_F^*$ et $b=\pi_L^{v_L(b)}b'$ avec $b'\in O_L^*$,  
alors $a\sim_l b$ si   
et seulement si on a l'\'egalit\'e $v_F(a)=v_L(b)$ et que l'image par
$\lll$ de la classe de 
$a'$ dans $O_F/P_F^m$ et la classe de $b'$ dans $O_L/P_L^m$ sont 
\'egales modulo $P_L^l$.\\

PROPRI\'ET\'ES

1) Pour tout $ x\in F$, $x$ et $\lam(x)$ sont $m$-proches,

2) Pour tout $a\in F$ et $b\in L$ qui sont $l$-proches, pour tout $i\in
   \mathbb{Z}$, $\pi_F^ia$ et $\pi_L^ib$ sont $l$-proches,  

3) Si $a_1$ et $a_2$ sont dans $F$ et $b_1$ et $b_2$ sont dans $L$, si
   $a_1\sim_l b_1$ et $a_2\sim_l b_2$, alors $a_1a_2\sim_l b_1b_2$, 

4)\label{prop4} Soit $A$ un ensemble fini et pour tout $i\in A$, $a_i$
un \'el\'ement de $F$ et $b_i$ un \'el\'ement de $L$ tel que
$\sum\limits_Aa_i\neq 0$. On pose
$l'=l+v_F(\sum\limits_Aa_i)-\min\limits_A(v_F(a_i))$. Si $m\geq l'$ et
pour tout $i\in A$, $a_i\sim_{l'}b_i$, alors on a~: 
$\sum\limits_Aa_i \sim_l \sum\limits_Ab_i$.\\   
 
Les premi\`eres trois  propri\'et\'es sont triviales. Pour d\'emontrer
la quatri\`eme on \'ecrit
$$
\lam(\sum\limits_Aa_i)-\sum\limits_Ab_i=
\lam(\sum\limits_Aa_i)-\sum\limits_A\lam(a_i)+\sum\limits_A\lam(a_i)-\sum\limits_Ab_i=
$$
$$
\lam(\sum\limits_Aa_i)-\sum\limits_A\lam(a_i)+\sum\limits_A\big(\lam(a_i)-b_i\big).
$$ 
Or, pour tout $i$ on a $a_i\sim_{l'} b_i$ et donc 
$$
\lam(a_i)-b_i \in
P_F^{l+v_F(\sum\limits_Aa_i)}
$$ 
d'o\`u
$$\sum\limits_A\big(\lam(a_i)-b_i\big) \in
P_F^{l+v_F(\sum\limits_Aa_i)}.$$ 
Reste \`a montrer que
$$
\lam(\sum\limits_Aa_i)-\sum\limits_A\lam(a_i)\in
 P_F^{l+v_F(\sum\limits_Aa_i)}.
$$
En \'ecrivant, pour tout $i$,
$a_i=\pi_F^{\min_Av_F(a_i)}a'_i, a'_i\in O_F$, on a
$$
\lam(\sum\limits_Aa_i)=\pi_L^{\min_Av_F(a_i)}\lam(\sum\limits_Aa'_i)
$$
et, pour tout $i$, 
$$
\lam(a_i)=\pi_L^{\min_Av_F(a_i)}\lam(a'_i)
.$$ 
Il
suffit donc de v\'erifier que
$\sum\limits_Aa'_i-\sum\limits_A\lam(a'_i)\in P_L^{l'}$. Mais comme
les $a'_i$ sont dans $O_F$ et les corps sont $m$-proches, alors par le
fait que $\lll$ est induit par $\lam$ (sous-section \ref{bijformcorps}),
$\lam(\sum\limits_Aa'_i)-\sum\limits_A\lam(a'_i)\in P_L^{m}$ ;  comme
$m\geq l'$ le r\'esultat est prouv\'e.\\

Sur le $F$-espace vectoriel $F^n$ on consid\`ere la valuation
$v_{F^n}((a_1;a_2;...a_n))=\min\limits_{1\leq i\leq n}v_F(a_i)$. On
rappelle que, si $F$ et $L$ sont $m$-proches on \'etend l'isomorphisme
$\lam$ de fa\c{c}on naturelle, composante par composante, en un
isomorphisme de $F^n$ sur $L^n$. Si $a=(a_1;a_2;...a_n)\in F^n$ et
$b=(b_1;b_2;...b_n)\in L^n$ et si $0< l\leq m$ on dit que $a$ et $b$
sont $l$-proches si $a$ et $b$ sont nuls tous les deux ou si $a$ est non
nul et que pour tout $i$ on a $b_i-\lam(a_i)\in P_L^{l+v_{F^n}(a)}$. Si
$U$ est un $F$-espace vectoriel de dimension finie $n$ muni d'une base,
et $V$ un $L$-espace vectoriel de m\^eme dimension $n$ muni d'une
base, alors on peut identifier $U$ \`a $F^n$ et $V$ \`a $L^n$ et
parler d'\'el\'ements $l$-proches de $U$ et $V$. Cette extension de la
d\'efinition des \'el\'ements proches s'applique en particulier aux
espaces de matrices $M_n(F)$ et $M_n(L)$.\\

\subsection{\'El\'ements proches et polyn\^omes} Si $P$ est un polyn\^ome en $n^2$ variables commutatives
$X_{11}, X_{12},...X_{nn}$ \`a coefficients dans $\mathbb{Z}$, si
$M=(m_{ij})\in M_n(F)$ (ou $M_n(L)$), on pose $P(M)=P(m_{ij})\in F$
(ou $L$). \'Enon\c{c}ons quatre propositions qui nous seront
utiles par la suite:\\ 

\begin{prop}\label{21} Si $M, M'\in M_n(F)$ alors $v_{M_n(F)}(MM')\geq
v_{M_n(F)}(M) + v_{M_n(F)}(M')$.\\ 
\end{prop}

\begin{prop}\label{22} Soit $M\in GL_n(F)$. Pour tout $k>0$ on a~:
$$M+M_n(P_F^{k-\vvf (M^{-1})})\subset K_F^kMK_F^k\subset
M+M_n(P_F^{k+\vvf (M)}).$$ 
\end{prop}

\begin{prop}\label{23} Si $k>0$ est fix\'e, si $M\in GL_n(F)$, en
posant 
$$m=k-\vvf (M) -\vvf (M^{-1})$$
on a~: si $F$ et $L$ sont
$m$-proches, alors $\zzzf(M)\in \zzz(K_F^k M K_F^k)$.\\ 
\end{prop}

\begin{prop}\label{24} Supposons que $F$ est de caract\'eristique non
nulle $p$. Soit $P\in {\mathbb{Z}}[X_{11},X_{12},...X_{nn}]$. Soient
$k>0$ et $M\in M_n(F)$ fix\'es. 

a) On suppose que tous les coefficients de $P$ se trouvent dans
l'ensemble $\{1,2...p-1\}$ et que $P(M)\neq 0$. On pose $S=\{s\
\text{tel que}\ s\ \text{est un mon\^ome de}\ P\}$ et  
$$m=k+v_F(P(M))-\min\limits_{s\in S}\ v_F(s(M))-\vvf(M)-\vvf(M^{-1}).$$
 Alors, si $F$ et $L$ sont $m$-proches, pour tout $N\in
\zzz(K_F^mMK_F^m)$ on a $P(M)\sim_kP(N)$. 

b) On ne fait aucune supposition sur les coefficients de $P$, mais on
suppose toujours que $P(M)\neq 0$. \'Ecrivons $P=Q+pR$ o\`u $Q,R\in
{\mathbb{Z}}[X_{11},X_{12},...X_{nn}]$ et tous les coefficients de $Q$
se trouvent dans l'ensemble $\{1,2...p-1\}$. Soit $\mu$ le degr\'e total
de $R$, $S=\{s\ \text{tel que}\ s\ \text{est un mon\^ome de}\ Q\}$ et  
$$m=\max\{k+v_F(P(M))-\min\limits_{s\in S}\ v_F(s(M))-\vvf(M)-\vvf(M^{-1});$$
$$k+v_F(P(M))+\max\{0;-\mu v_{M_n(F)}(M)\}\}.$$ 
 Alors, si $F$ et $L$ sont $m$-proches, pour tout $N\in
\zzz(K_F^mMK_F^m)$ on a $P(M)\sim_kP(N)$. 

c) Supposons que $P(M)=0$. Alors il existe $m$ tel que, si $F$ et $L$
sont $m$-proches, pour tout $N\in \zzz(K_F^mMK_F^m)$ on ait
$v_L(P(N))\geq k$.\\ 
\end{prop}
\ \\
{\bf D\'emonstrations.}\\ 
\ \\
PROPOSITION \ref{21} :  Pour tout $1\leq i\leq n$ et tout $1\leq j \leq n$ on a $$v_F(\sum_{k=1}^{n} m_{ik}m'_{kj})\geq \min\limits_{1\leq k\leq n} \ v_F(m_{ik}) + \min\limits_{1\leq k\leq n}\ v_F(m'_{kj})$$
 d'o\`u le r\'esultat.
\qed
\ \\
\ \\
PROPOSITION \ref{22} : Si $A=M+B$ o\`u $B\in M_n(P_F^{k-\vvf
(M^{-1})})$ alors  
$$A=M(Id+M^{-1}B)\in MK_F^k$$ 
par la prop.\ref{21}, d'o\`u la premi\`ere inclusion.  

Si $A=(Id+B)M(Id+C)$ avec $B,C\in M_n(P_F^k)$, alors 
$$A=M+BMC+BM+MC$$
 o\`u $BMC,BM$ et $MC$ se trouvent dans $M_n(P_F^{k+\vvf(M)})$ par la
prop.\ref{21}, d'o\`u la deuxi\`eme inclusion.
\qed
\ \\
\ \\
PROPOSITION \ref{23} :  Soit $M=BAC$, o\`u $B\in GL_n(O_F)$, $C\in
GL_n(O_F)$ et $A\in {\mathcal{A}}_F$. On a vu qu'alors  
$$\zzz(K_F^k M K_F^k)=K_L^k\zzzf(B)\zzzf(A)\zzzf(C)K_L^k.$$
 Pour montrer que $\zzzf(M)\in \zzz(K_F^k M K_F^k)$ il suffit, par la
prop.\ref{22}, de montrer que  
$$\zzzf(BAC)-\zzzf(B)\zzzf(A)\zzzf(C)\in M_n(P_L^u)$$
 o\`u 
$$u=k-v_{M_n(L)}((\zzzf(B)\zzzf(A)\zzzf(C))^{-1})= k-v_{M_n(F)}(M^{-1})=m+v_{M_n(F)}(M).$$
On a utilis\'e pour la deuxi\`eme \'egalit\'e le fait que 

- $\zzzf(B)\in GL_n(O_L)$, 

- $\zzzf(C)\in GL_n(O_L)$ et 

- $v_L(\zzzf(A^{-1}))=v_F(A^{-1})=v(C^{-1}M^{-1}B^{-1})=v_F(M^{-1})$.
Pour les m\^emes raisons, $v_{M_n(F)}(M)=v_{M_n(F)}(A)$ et si on \'ecrit 
$$A=diag(\pi_F^{a_1};\pi_F^{a_2}...\pi_F^{a_n}),$$
 alors $v_{M_n(F)}(A)=a_1$. Donc la relation \`a montrer est 
$$\zzzf(BAC)-\zzzf(B)\zzzf(A)\zzzf(C)\in M_n(P_L^{m+a_1})\ ;$$
 ou encore~: 
$$\zzzf(B\pi_F^{-a_1}AC)-\zzzf(B)\zzzf(\pi_F^{-a_1}A)\zzzf(C)\in
M_n(P_L^{m}),$$  
qui est \'evidente, car $B,C,\pi_F^{-a_1}A\in M_n(O_F)$ et on peut
appliquer le fait que l'application $\zzz$ est induite par la
restriction de $\zzzf$.
\qed 
\ \\
\ \\
PROPOSITION \ref{24} : Le partage du premier r\'esultat de ce
th\'eor\`eme (hypoth\`ese $P(M)\neq 0$) en point a) et point b) vient
du  fait que, si $F$ et $L$ sont $m$-proches la somme \`a $l$ termes
$1+1+...+1$ dans les corps $F$ et $L$ respectivement donne des
\'el\'ements $m$-proches pour $l<p$ (voir la propri\'et\'e 4 page
\pageref{prop4}),  mais pas pour 
$l=p$. Le point c) traite du cas $P(M)=0$ o\`u le r\'esultat est de
nature diff\'erente.\\ 

a) Par la proposition \ref{23}, $\zzzf(M)\in \zzz(K_F^mMK_F^m)$. Par
 la proposition \ref{22}, $N-\zzzf(M)\in
 M_n(P_F^{m-v_{M_n(F)}(M^{-1})})$ et donc $N$ et $\zzzf(M)$ sont
 $[m-v_{M_n(F)}(M^{-1})-v_{M_n(F)}(M)]$-proches. Mais  
$$m-v_{M_n(F)}(M^{-1})-v_{M_n(F)}(M)=k+v_F(P(M))-\min\limits_{v\in
S}v_F(s(M)).$$ 
Donc les coefficients de m\^emes indices des deux matrices
respectivement sont $[k+v_F(P(M))-\min\limits_{v\in
S}v_F(s(M))]$-proches. Par la propri\'et\'e 3, pour chaque $s\in S$,
$s(M)$ et $s(N)$ sont $[k+v_F(P(M))-\min\limits_{v\in
S}v_F(s(M))]$-proches. Les coefficients devant ces mon\^omes se
trouvent dans l'ensemble $\{1,2...p-1\}$. Les \'el\'ements $1_F$ et
$1_L$ sont $m$-proches par la condition impos\'ee dans la d\'efinition
de l'application $\lam$. Alors, par la propri\'et\'e 4, pour tout
$l\in \{1,2...p-1\}$, les sommes \`a $l$ termes $1+1+...+1$ dans $F$
et $L$ respectivement sont $m$-proches. Ainsi $P(M)$ et $P(N)$ sont
des sommes non nulles d'\'el\'ements $[k+v_F(P(M))-\min\limits_{v\in
S}v_F(s(M))]$-proches deux \`a deux. Donc, par la propri\'et\'e 4,
$P(M)$ et $P(N)$ sont $k$-proches. 

b) Par le point a) et le choix de $m$, $Q(M)$ et $Q(N)$ sont
$k$-proches. Remarquons que la caract\'eristique de $F$ \'etant $p$,
$Q(M)=P(M)$. Il suffit de montrer donc que $Q(N)+pR(N)$ et $Q(N)$ sont
$k$-proches. Or, si $L$ est un corps de caract\'eristique nulle
$m$-proche de $F$, alors l'image de la somme \`a $p$ termes
$1+1+...+1$  dans $O_L/P_L^m$ est la classe de $0$, donc la valuation
de l'\'el\'ement  $1+1+...+1$ ($p$ fois 1) est sup\'erieure \`a
$m$. Par le choix de $m$ dans l'hypoth\`ese, $pR(N)\in
P_L^{k+v_L(Q(N))}$ et donc $Q(N)+pR(N)$ et $Q(N)$ sont $k$-proches.  

c) Consid\'erons le polyn\^ome $Q(X)=P(X)+1$. On a $Q(M)\neq 0$ et on
peut appliquer le point b) \`a $Q$. Or, si $Q(M)$ et $Q(N)$ sont
$k$-proches, alors $Q(N)-\zzzf(Q(M))\in P_L^k$. Mais  
$$Q(N)-\zzzf(Q(M))=(P(N)+1)-1=P(N)$$ 
d'o\`u le r\'esultat.
\qed
\ \\

\subsection{\'El\'ements proches et polyn\^omes caract\'eristiques,
cas de $GL_n$} Soit $F$ un corps local de caract\'eristique non nulle
$p$. Soit $P_M$ un polyn\^ome unitaire de degr\'e $n$ \`a coefficients
dans $F$ et s\'eparable (sans racine multiple),  et $M$ 
la matrice compagnon de $P_M$ . Soit $l>0$.

\begin{prop}\label{polcarseridica} Il existe
deux entiers, $m\geq l$ et $s$, qui ne d\'ependent que de $P_M$ et de $l$, tels que, si $L$
est un corps local de  
caract\'eristique nulle $m$-proche de $F$, on ait~: si $g$ est un 
\'el\'ement  de $G_L$ (resp. de $G_F$) dont le polyn\^ome caract\'eristique est
$s$-proche de $P_M$, alors $g$ est conjugu\'e \`a un \'el\'ement de
$\zzz(K_F^lMK_F^l)$ (resp. de $K_F^lMK_F^l$).
\end{prop}
\ \\
{\bf D\'emonstration.} On pose 
$$s=l-v_{M_n(F)}(M)-v_{M_n(F)}(M^{-1})$$
 et 
$$m=s.$$
 Par la proposition \ref{23} et le choix de $m$, $\zzzf(M)\in
\zzz(K_F^lMK_F^l)$ et donc  
$$\zzz(K_F^lMK_F^l)=K_L^l\zzzf(M)K^l_L.$$
 D'autre part, si $K_FMK_F=K_FAK_F$, $A\in {\mathcal A}_F$, alors  
$$v_{M_n(F)}(M)=v_{M_n(F)}(A)$$ 
 et
$$v_{M_n(F)}(M^{-1})=v_{M_n(F)}(A^{-1}).$$
 Comme on a $\zzz(K_F^lMK_F^l)=K_L^l\zzzf(M)K^l_L$, on en d\'eduit que 
$$\zzzf(M)\in K_L^l\zzzf(A)K^l_L$$ 
 et que 
$$v_{M_n(L)}(\zzzf(M))=v_{M_n(L)}(\zzzf(A))$$ 
 et 
$$v_{M_n(L)}(\zzzf(M)^{-1})=v_{M_n(L)}(\zzzf(A)^{-1}).$$
 Bref, 
$$v_{M_n(F)}(M)=v_{M_n(L)}(\zzzf(M))$$ 
et 
$$v_{M_n(F)}(M^{-1})=v_{M_n(L)}(\zzzf(M)^{-1}).$$
On peut donc \'ecrire~:
$$s=l-v_{M_n(L)}(\zzzf(M))-v_{M_n(L)}(\zzzf(M)^{-1}).$$
On remarque par ailleurs que $\zzzf(M)$ est la matrice compagnon de
$\zzzf(P_M)$. Donc, si $P$ est un polyn\^ome $s$-proche de $P_M$,
alors $P$ est un polyn\^ome $s$-proche du polyn\^ome caract\'eristique
$P_{\zzzf(M)}$ de $\zzzf(M)$ aussi. Donc la matrice compagnon
$Comp(P)$ de $P$ sera $s$-proche de $\zzzf(M)$. Par cons\'equent  
 
$$
Comp(P)-\zzzf(M)\in
M_n(P_L^{s+v_{M_n(L)}(\zzzf(M))})=M_n(P_L^{l-v_{M_n(L)}(\zzzf(M)^{-1})}).
$$ 
Par la proposition \ref{22} on a alors 

$$
Comp(P)\in K_L^l\zzzf(M)K^l_L=\zzz(K_F^lMK_F^l).
$$
On conclut par le fait que si $g\in GL_n(F)$ a le m\^eme polyn\^ome
caract\'eristique $P$ que $Comp(P)$, alors $g$ et $Comp(P)$ sont
conjugu\'es (car $P$ est sans racine multiple).
\qed
\ \\
\begin{prop}\label{caractseridica}  Soit $\pi$ une \care\ de
$G_F$. Soit $M$ un \'el\'ement elliptique r\'egulier de $G_F$.
Soit $P_M$ le polyn\^ome caract\'eristique de $M$. 
Il existe alors $m$ et $s$ qui ne d\'ependent
que de $\pi$ et de $P_M$ tels que, si $L$ est un corps local de
caract\'eristique  
nulle $m$-proche de $F$, on ait~: pour tout \'el\'ement $g$ de $G_L$
dont le polyn\^ome caract\'eristique est $s$-proche de $P_M$, on a  
$$\chi_{\zzz(\pi)}(g)=\chi_{\pi}(M).$$  
\end{prop}
\ \\
{\bf D\'emonstration.} On peut supposer  que $M$ est la matrice compagnon de
$P_M$, puisque les
caract\`eres sont constants sur une classe de conjugaison. Dans [Ba2]
nous montrons \`a la page 65, au cours de la d\'emonstration du
th.4.3, qu'il existe $l$ tel que  $\chi_{\pi}$ soit constant sur $
K_F^lMK_F^l$ et $m$ tel que, si $L$ est un corps $m$-proche de $F$,
alors $ \chi_{\zzz(\pi)}$ soit constant sur 
$\zzz(K_F^lMK_F^l)$, et ces deux constantes sont \'egales. Les entiers
$l$ et $m$ ne d\'ependent que de $M$ et de $\pi$. (Ce rel\`evement
local des caract\`eres ne se fait que pour des \cares\ et pour un $M$
elliptique r\'egulier, car il passe par un rel\`evement local de
l'int\'egrale orbitale d'un pseudocoefficient.)
Nous appliquons ensuite la prop.\ref{polcarseridica} pour ce $l$ et,
quitte \`a augmenter $m$, nous avons le r\'esultat. Les entiers $m$ et
$s$ obtenus ne d\'ependent que de $l$, $M$ et $\pi$, mais $l$ ne
d\'epend \`a son tour que de $M$ et $\pi$, et $M$ ne d\'epend que de
$P_M$, \'etant sa matrice compagnon. 
\qed

\subsection{\'El\'ements proches et polyn\^omes caract\'eristiques,
cas des formes int\'erieures de $GL_n$} Dans cette sous-section on
d\'emontre un r\'esultat sur 
les formes int\'erieures de $GL_n(F)$ (la proposition \ref{34}). Soit
$D_F$ une alg\`ebre \`a division centrale de dimension $d^2$ sur
$F$. Soit $E$ une extension non ramifi\'ee de dimension $d$ sur $F$
incluse dans $D_F$. On supppose qu'on a fix\'e une uniformisante
$\pi_F$ de $F$ (et de $E$ aussi), ainsi qu'une uniformisante
$\pi_{D_F}$ de $D_F$ et un g\'en\'erateur $\sigma_E$ de $Gal(E/F)$ qui
correspondent \`a $D_F$ comme dans la sous-section \ref{algadiv}. Soit $r$ un
entier strictement positif et 
$G'_F=GL_r(D_F)$. On pose $n=rd$. Chaque fois qu'on se donne $L$ un
corps local $m$-proche de $F$, on consid\`ere que le triplet
correspondant est choisi de fa\c{c}on \`a ce que l'uniformisante de
$F$ qui y appara\^\i t soit $\pi_F$. On reprend alors toutes les
notations de la section pr\'ec\'edente pour $K$, $D_L$ et tous les
objets qui leur sont associ\'es, avec une seule exception~: la
dimension de $D_F$ sur $F$ est not\'ee ici $d$ alors que dans la
section pr\'ec\'edente elle \'etait not\'ee $n$. On rappelle que la
base de voisinages $\{K_F^l\}_{l\in \mathbb{N}}$ de l'identit\'e avec
laquelle on a travaill\'e sur $G'_F$ est associ\'ee \`a la base de
voisinages $\{P_{D_F}^{dl}\}_{l\in \mathbb{N}}$ de $0$ et non pas \`a
la base de voisinages $\{P_{D_F}^{l}\}_{l\in \mathbb{N}}$. \\ 

\begin{prop}\label{31} Si $M, M'\in M_r(D_F)$ alors
$v_{M_r(D_F)}(MM')\geq v_{M_r(D_F)}(M) + v_{M_r(D_F)}(M')$.\\ 
\end{prop}

\begin{prop}\label{32} Soit $M\in GL_r(D_F)$. Pour tout $k>0$ on a~:
$$M+M_r(P_{D_F}^{d(k-v_{M_r(D_F)}(M^{-1}))})\subset K_F^kMK_F^k\subset 
M+M_r(P_{D_F}^{d(k+ v_{M_r(D_F)}(M))}).$$ 
\end{prop}

\begin{prop}\label{33} Si $k>0$ est fix\'e, si $M\in GL_r(D_F)$, en
posant $m=k-v_{M_r(D_F)}(M) -v_{M_r(D_F)} (M^{-1})$ on a~: si $F$ et
$L$ sont $m$-proches, alors $\zzzd(M)\in \zzd(K_F^k M K_F^k)$.\\ 
\end{prop}
\ \\
{\bf D\'emonstrations.} Les d\'emonstrations des propositions \ref{21}
et \ref{22} s'appliquent aux propositions \ref{31} et \ref{32} sans
changement. Pour la proposition \ref{33} il y a un petit probl\`eme
car l'uniformisante de $D_F$ ne commute pas avec tous les \'el\'ements de
$D_F$ donc il faut v\'erifier que 
$\zzzd(BAC)-\zzzd(B)\zzzd(A)\zzzd(C)\in M_n(P_L^{d(m-a_1)})$ est
toujours vrai. 
On multiplie par $\pi_{D_F}^{-a_1}$ comme dans la d\'emonstration de
la proposition \ref{23} et on \'ecrit~:  

$$
\pi_{D_L}^{-a_1}\zzzd(BAC)=\zzzd(\pi_{D_F}^{-a_1}BAC)=\zzzd\big(\sigma_E^{a_1}(B)(\pi_{D_F}^{-a_1}A)C\big)
$$  
car nous avons d\'efini $\zzzd$ de sorte qu'elle commute \`a la
multiplication par les uniformisantes (propri\'et\'e 2).  

Maintenant 
$$\pi_{D_L}^{-a_1}\zzzd(B)\zzzd(A)\zzzd(C)=\sigma_K^{a_1}\big(\zzzd(B)\big)\pi_{D_L}^{-a_1}\zzzd(A)\zzzd(C)$$\\ 
$$=\sigma_K^{a_1}\big(\zzzd(B)\big)\zzzd(\pi_{D_F}^{-a_1}A)\zzzd(C).$$
On a aussi
$$\zzzd\big(\sigma_E^{a_1}(B)\big)=\sigma_K^{a_1}\big(\zzzd(B)\big)$$
par la relation \ref{**}, page \pageref{**}. Il faut donc montrer que  
$$\zzzd\big(\sigma_E^{a_1}(B)(\pi_{D_F}^{-a_1}A)C\big)-\zzzd\big(\sigma_E^{a_1}(B)\big)\zzzd(\pi_{D_F}^{-a_1}A)\zzzd(C)\in
M_r(P_{D_F}^k).$$ 
On conclut comme dans la d\'emonstration de la proposition \ref{23},
car $\sigma_E^{a_1}(B),\ \pi_{D_F}^{-a_1}A$ et $C$ sont dans
$M_r(O_{D_F})$.\\

\begin{prop}\label{34} Soient $M'\in G'_F$ et $k\in {\mathbb{N}}$. Il
existe un entier $m$ tel que, si $F$ et $L$ sont $m$-proches, alors pour tout
$g'\in \zzd(K_{D_F}^mM'K_{D_F}^m)$ les polyn\^omes caract\'eristiques
de $M'$ et $g'$ sont $k$-proches. 
\end{prop}
\ \\
{\bf D\'emonstration.} On rappelle la proposition de la page 295,
[Pi]~:  

Soit $A$ une alg\`ebre centrale simple sur $F$ de dimension
$n^2$. Soit $E$ une extension de dimension $n$ de $F$. Alors, si on a
un morphisme d'alg\`ebres unitaires $\Psi : A\rightarrow M_n(E)$, pour
tout \'el\'ement $g$ de $A$, le polyn\^ome caract\'eristique de
$\Psi(g)$ (qui a priori a des coefficients dans $E$) a tous ses
coefficients dans $F$ et c'est le polyn\^ome caract\'eristique de
$g$. 

Dans notre cas, $A=M_r(D_F)$. Elle agit sur $D_F^r$. En \'ecrivant
$D_F=\bigoplus\limits_{0\leq i\leq d}\pi_{D_F}^iE$ on a un
isomorphisme $D_F^r\simeq E^n$ et par cons\'equent une action de
$M_r(D_F)$ sur $E^n$. On a obtenu donc un morphisme d'alg\`ebres $\Psi
: M_r(D_F)\rightarrow M_n(E)$. Le polyn\^ome caract\'eristique de $M'$
est alors \'egal au polyn\^ome caract\'eristique de $\Psi(M')$. On va
calculer ce dernier en fonction des coefficients de $M'$. Supposons
que $M'$ s'\'ecrit $M'=(m'_{ij})_{1\leq i,j\leq r}$ et que $\Psi(M')$
s'\'ecrit $(n_{st})_{1\leq s,t\leq n}$. Supposons maintenant que pour
tout $i$ et $j$, $m'_{ij}$ s'\'ecrit sur la base
$1,\pi_{D_F},\pi_{D_F}^2...\pi_{D_F}^{d-1}$ de $D_F$ sur $E$~: 
$$m'_{ij}=\sum\limits_{0\leq k\leq d-1}\pi_{D_F}^ke_{ij}^k,\
e_{ij}^k\in E.$$ 
Pour tout $0\leq l\leq d-1$, pour tout tel $e_{ij}^k$
posons $e_{ij}^{kl}=\sigma_E^l(e_{ij}^k)$. On peut  fabriquer une
matrice ${\mathcal U}(M')=(u_{vw})_{1\leq v,w\leq n}$ a $n$ lignes et
$n$ colonnes et \`a coefficients dans $E$ en posant pour tout $1\leq
v,w\leq n$~: 
$u_{vw}=e_{ij}^{kl}$ o\`u $i, j, k$ et $l$ sont d\'efinis comme suit~:
$l$ est le quotient de la division euclidienne de $v$ par $r$, $i-1$
est le reste de la division euclidienne de $v$ par $r$, $k$ est le
quotient de la division euclidienne de $w$ par $r$ et  $j-1$ est le
reste de la division euclidienne de $w$ par $r$.\\  
\ \\
{\bf Notation~:} Si $P$ est un polyn\^ome dans
${\mathbb{Z}}[X_{11},X_{12},...X_{nn}][t]$, si $A=(a_{ij})_{1\leq
i,j\leq n}$ est une matrice dans $M_n(F)$, si $x\in F$, l'\'el\'ement
$P(a_{11};a_{12};...;a_{nn};x)$ de $F$ sera not\'e abr\'eg\'e
$P(A;x)$.

\begin{lemme}\label{titi} Pour tout $1\leq s,t\leq n$ il existe un
polyn\^ome ind\'ependant  de $M'$, $P_{st}\in
{\mathbb{Z}}[X_{11},X_{12},...X_{nn}][t]$, dont le degr\'e total  en
les variables $X_{11},X_{12},...X_{nn}$ est 1, tel qu'on ait
$n_{st}=P_{st}({\mathcal U}(M');\pi_F)$. 
\end{lemme}
\ \\
{\bf D\'emonstration.} Il suffit de v\'erifier cette propri\'et\'e
pour des matrices du type $M_{i_0j_0}^{ke}=(m'_{ij})_{1\leq i,j \leq
r}$ o\`u $m'_{ij}=\delta_{i_0i}\delta_{j_0j}\pi_{D_F}^ke$ o\`u
$i_0,j_0$ sont des entiers entre $1$ et $r$, $k$ est un entier entre
$1$ et $d$, et  $e\in E$, car l'ensemble form\'e par ces matrices
engendre $M_r(D_F)$ sur $\mathbb{Z}$ et les polyn\^omes consid\'er\'es
sont de degr\'e 1 en les  $n^2$ premi\`eres variables. Soit
$d_1,d_2...d_r$ la base canonique de $D_F^r$. L'\'el\'ement
$M_{i_0j_0}^{ke}$ agit sur $D_F^r$ en envoyant $d_{i}$ sur $0$ pour
tout $i\neq i_0$ et en envoyant $d_{i_0}$ sur
$\pi_{D_F}^ked_{j_0}$. Si on se repr\'esente la matrice
$\Psi(M_{i_0j_0}^{ke})$ par blocs de taille $d\times d$, alors tous
ces blocs sont nuls \`a l'exception de celui qui se trouve dans la
position $i_0j_0$, et ce dernier est \'egal \`a $X=(x_{ij})_{1\leq
i,j\leq d}$ o\`u les $x_{ij}$ sont donn\'es par~:

\ \ \ \ \ \ \ \ - si $1\leq i\leq k$, alors
$x_{ij}=\delta_{i,j-d+k}\pi_F\sigma^{k-d+j-1}(e)$ 

\ \ \ \ \ \ \ \ - si $k+1\leq i\leq n$, alors
$x_{ij}=\delta_{i,j+k}\sigma^{j-1}(e).$ 

Le lemme est v\'erifi\'e.\qed
 
\begin{lemme} Il existe des polyn\^omes $P_0,P_1...P_{n-1}\in
{\mathbb{Z}}[X_{11},X_{12},...X_{nn}][t]$ tels que pour toute matrice
$M'$ dans $M_r(D_F)$, le coefficient de $X^i$, $1\leq i\leq n-1$, dans
le polyn\^ome caract\'eristique de $M'$ soit \'egal \`a $P_i({\mathcal
U}(M');\pi_F)$. 
\end{lemme}
\ \\
{\bf D\'emonstration.} C'est \'evident par le lemme \ref{titi} plus haut.\\
\begin{lemme}\label{toto} Les  points a), b) et c) de la proposition
\ref{24} sont v\'erifi\'es si on remplace  ${\mathbb{Z}}[X_{11},X_{12},...X_{nn}]$ par ${\mathbb{Z}}[X_{11},X_{12},...X_{nn}][t]$,
$P(M)$ par $P(M;\pi_F)$ et $P(N)$ par $P(N;\pi_L)$.  
\end{lemme}
\ \\
{\bf D\'emonstration.} La d\'emonstration marche identiquement en
tenant compte que, si $F$ et $L$ sont $m$-proches, alors $\pi_F$ et
$\pi_L$ sont $m$-proches. Une autre fa\c{c}on de d\'emontrer ce lemme
est de le voir comme un cas particulier de la proposition \ref{24}~:
on applique la proposition \ref{24} \`a $(n+1)^2$ variables.\qed
\ \\

D\'emontrons maintenant la proposition \ref{34}. On remarque que, si
$M$ et $M'$ sont dans $M_r(D_F)$, et si $M-M'\in M_r(P_{D_F}^{dh})$,
alors pour tout $i,j$, si on \'ecrit $m_{ij}=\sum\limits_{0\leq k\leq
d-1}\pi_{D_F}^ke_{ij}^k$ et $m'_{ij}=\sum\limits_{0\leq k\leq
d-1}\pi_{D_F}^k{e'}_{ij}^k$, on a pour tout $k$~:
$e_{ij}^k-{e'}_{ij}^k\in P_E^h$. Par cons\'equent, pour tout entier
$k_1$ fix\'e, il existe un entier $k_2$ tel que, si $F$ et $L$ sont
$k_2$-proches, pour tout $M\in
\bar{\zeta}^{k_2}_{D_FD_L}(K_{D_F}^{k_2}M'K_{D_F}^{k_2})$ on ait~:
${\mathcal{U}}(M)\in
\bar{\zeta}^{k_2}_{EK}(K_E^{k_1}{\mathcal{U}}(M')K_E^{k_1})$ (on a
utilis\'e les propositions \ref{22}, \ref{32}, \ref{23}, \ref{33} et
le fait que si $e,e'\in E$ sont $k$-proches alors, pour tout $\sigma
\in Gal(E/F)$, $\sigma(e)$ et $\sigma(e')$ sont $k$-proches). 

Soit maintenant $M'$ comme dans l'hypoth\`ese de la proposition
\ref{34}. On pose $N=\min\limits_{0\leq i\leq
n-1}v_{M_n(E)}(P_i({\mathcal{U}}(M')))$ qui a un sens parce qu'au
moins $P_0(({\mathcal{U}}(M')))$ est non nul (car \'egal \`a
$det(M'))$. L'entier $N$ n'est autre que la valuation du polyn\^ome
caract\'eristique de $M''$ vu comme \'el\'ement de $F^n$. En
appliquant la proposition 4 b) \`a la matrice ${\mathcal{U}}(M')\in
M_n(E)$ on trouve qu'il existe un $k_0$ tel que, si $L$ est
$k_0$-proche de $F$, pour toute matrice $M''\in
\bar{\zeta}^{k_0}_{EK}(K_{EK}^{k_0}{\mathcal{U}}(M')K_{EK}^{k_0})$,
pour tout $i$ entre $0$ et $n-1$ tel que $P_i({\mathcal{U}}(M'))\neq
0$, $P_i({\mathcal{U}}(M'))$ et $P_i(M'')$ soient $k$-proches. En
appliquant le lemme \ref{toto}c) aux polyn\^omes $P_i$ qui v\'erifient
$P_i({\mathcal{U}}(M'))=0$ on trouve qu'il existe un $k'_0$ tel que,
si $L$ est $k'_0$-proche de $F$, pour toute matrice $M''\in
\bar{\zeta}^{k'_0}_{EK}(K_{EK}^{k'_0}{\mathcal{U}}(M')K_{EK}^{k_0})$
on a $v_{M_n(K)}(P_i(M''))\geq k+N$. En posant $k_1=\max\{k_0;k'_0\}$
l'entier $m=k_2$ (voir quelques lignes plus haut pour $k_2$) v\'erifie
les propri\'et\'es requises par la proposition \ref{34}.\\ 
\ \\

\newpage

\section{Preuve de la correspondance}
\ \\

Soient $F$ un corps local non archim\'edien et $D$ une alg\`ebre \`a
division centrale sur $F$ et de dimension $d^2$. On pose $n=rd$,
$G=GL_n(F)$ et 
 $G'=GL_r(D)$.  Si $O_F$ (resp.$O_D$) est l'anneau des entiers de $F$
 (resp.$D$), on fixe des mesures de Haar sur $G$ et $G'$ telles que le
 volume de $GL_n(O_F)$ (resp.$GL_r(O_D)$) soit \'egal \`a 1.
Rappelons qu'un \'el\'ement  $g$ de
 $G$ ou $G'$  est dit {\it semisimple r\'egulier} (resp. {\it
elliptique r\'egulier}) si le polyn\^ome caract\'eristique de $g$ est
s\'eparable (resp. irr\'eductible s\'eparable) et que, si $g$ est 
 un \'el\'ement de $G$ et $g'$ un \'el\'ement de $G'$ on dit que $g$
correspond \`a $g'$ et on \'ecrit
 $g\lra g'$ si $g$ et $g'$ sont semisimples r\'eguliers et ont le
 m\^eme polyn\^ome caract\'eristique. Si $\pi$
 est une repr\'esentation lisse de longueur finie de $G$ ou $G'$,
 alors on note $\chi_{\pi}$ le caract\`ere de $\pi$ vu comme
 fonction localement constante sur l'ensemble des \ssrs\   du groupe
 en question.   
  
On note $E^2(G)$ l'ensemble des classes d'\'equivalence de \ecis\ de
$GL_n(F)$ et $E^2(G')$ l'ensemble des classes d'\'equivalence de
\ecis\ de $G'$. On fixe une fois pour toutes un caract\`ere additif non trivial $\psi_F$ de $F$,
trivial sur l'anneau des entiers de $F$. Tous les facteurs $\epsilon'$
des repr\'esentations de $G$ ou $G'$
seront calcul\'es \`a partir de $\psi_F$. Nous voulons montrer le
th\'eor\`eme suivant, dit 
correspondance de Jacquet-Langlands.\\ 

\begin{theo}\label{correspforte} Il existe une unique application~:
$${\bf{C}}:E^2(G)\to E^2(G')$$
telle que pour tout $\pi\in E^2(G)$ on ait 
\begin{equation}\label{forte}
\chi_{\pi}(g)=(-1)^{n-r} \chi_{{\bf C}(\pi)}(g')
\end{equation}
pour tous $g\lra g'$. 

L'application $\bf{C}$ est bijective. Les repr\'esentations $\pi$ et
${\bf{C}}(\pi)$ ont le m\^eme facteur $\epsilon'$.
\end{theo}

Nous voudrions montrer ce th\'eor\`eme en partant du fait qu'il a
d\'ej\`a \'et\'e montr\'e en caract\'eristique nulle (dans [DKV]) et en
utilisant 
la construction de situations proches de la section 2. L'id\'eal
serait de prouver que pour les deux groupes $G$ et $G'$ on peut
relever localement les caract\`eres des repr\'esentations au voisinage
des \'el\'ements r\'eguliers. En pratique cela est difficile \`a
envisager, et tout ce dont on dispose est un rel\`evement local des
caract\`eres seulement pour $G$ et seulement au voisinage des
\'el\'ements elliptiques, et seulement pour les \cares\ (cette
deni\`ere condition n'est pas vraiment une contrainte dans le cadre
pr\'esent). Ce 
r\'esultat, qui est une cons\'equence imm\'ediate des r\'esultats dans
[Ba2], permet n\'eanmoins de prouver le th\'eor\`eme plus faible
ci-dessous, par comparaison avec la caract\'eristique nulle et
\`a l'aide des r\'esultats prouv\'es dans la section 4~:

\begin{theo}\label{correspfaible}{\rm Correspondance faible.}\\
Il existe une bijection~:
$${\bf{C}'}:E^2(G)\to E^2(G')$$
telle que pour tout $\pi\in E^2(G)$ on ait 
$$\chi_{\pi}(g)=(-1)^{n-r} \chi_{{\bf C}(\pi)}(g')$$
pour tout $g\lra g'$ {\underline{elliptiques r\'eguliers}}.  

\end{theo}

Attirons l'attention sur le fait que cette variante faible n'est pas
satisfaisante parce qu'elle ne permet pas d'\'etendre la
correspondance de Jacquet-Langlands aux repr\'esentations {\it qui ne
sont 
pas essentiellement de carr\'e int\'egrable} (voir pour l'instant [Ba1],
ch.4). Nous ne pouvions donc pas nous arr\^eter l\`a. Pour montrer le
th.\ref{correspforte} (variante forte) en caract\'eristique non nulle,
nous remarquons  
que le th.\ref{correspfaible} permet de montrer les relations
d'orthogonalit\'e des caract\`eres des \cares\ sur $G'$ (voir
l'annexe), par 
transfert, une fois qu'on a ce r\'esultat sur $G$. Nous \'enon\c{c}ons  
donc un troisi\`eme th\'eor\`eme~:

\begin{theo}\label{ortogonalitate}
Les relations d'orthogonalit\'e des caract\`eres sont valables sur
$G$ et $G'$.
\end{theo} 
\ \\
On sait pour l'instant qu'il est vrai si la caract\'eristique du
corps de base 
est nulle ([Cl]) et qu'il est vrai pour $G$ si la caract\'eristique du corps
de base est non nulle ([Ba2]). Nous suivons ensuite le sch\'ema suivant~:

On montre que, sur un corps de base fix\'e (toute caract\'eristique),

1) le th.\ref{ortogonalitate} implique le 
th.\ref{correspforte},

2) le th.\ref{correspforte} implique le th.\ref{correspfaible}

3) le th.\ref{correspfaible} implique le th.\ref{ortogonalitate}.

La preuve de 1) est une r\'ecriture de
[DKV] telle qu'il soit clair que la caract\'eristique n'intervient
pas, en particulier en supprimant toute r\'ef\'erence aux r\'esultats
de l'appendice 
1 dans [DKV], indisponibles en caract\'eristique non nulle. 
Le 2) est trivial.  
Comme le
th.\ref{ortogonalitate} est prouv\'e pour $G$ en toute
carac\'eristique, le 3) est imm\'ediat par transfert sur $G'$. La partie
originale de la d\'emonstration est la preuve (qui utilise la
th\'eorie des corps proches 
d\'evelopp\'ee aux sections 2 et 4, ainsi que les r\'esultats de la
section 3), de l'implication  

4) si le th.\ref{correspforte}
est vrai sur tout corps de caract\'eristique nulle, alors le
th.\ref{correspfaible} est vrai sur tout corps de   
caract\'eristique non nulle. 

C'est pour pouvoir d\'emontrer cette implication qu'on a du int\'egrer
la relation sur les facteurs $\epsilon'$ \`a l'\'enonc\'e du
th.\ref{correspforte}.\\

Dans cette situation, comme le th.\ref{ortogonalitate} est connu en
caract\'eristique nulle, 1), 2), 3) et 4) impliquent que les trois
th\'eor\`emes sont vrais en toute caract\'eristique.
\ \\
Montrons que {\it le th.\ref{ortogonalitate} implique
le th.\ref{correspforte}.}
\ \\
{\bf D\'emonstration.} On fixe $D$ et on raisonne par r\'ecurrence sur
$r$. L'hypoth\`ese de 
r\'ecurrence est~:\\  

$({\bf H}_k)$ Le th\'eor\`eme \ref{correspforte} est vrai pour
$G=GL_{dk}(F)$ et $G'=GL_k(D)$.\\

Le premier pas ($k=1$), connu en caract\'eristique nulle ([Ro]) a
\'et\'e trait\'e dans [Ba2] pour le cas de la caract\'eristique non nulle. Nous
supposons que le th\'eor\`eme plus haut est v\'erifi\'e pour tout   
entier $k$ de 1 \`a $r-1$. On pose $G=GL_{dr}(F)$ et $G'=GL_r(D)$,
pour v\'erifier l'hypoth\`ese ${\bf H}_r$.\\

Soit $\pi_0$ une \care\ de $G$ de
 caract\`ere central (unitaire) $\omega$. 
Consid\'erons un corps global $\mathbb{F}$ et une alg\`ebre \`a 
division $\mathbb{D}$ sur  
$\mathbb{F}$ tels que~: 

- il existe une place $v_0$ de $\mathbb{F}$ telle que
${\mathbb{F}}_{v_0} \simeq F$ et 
${\mathbb{D}}_{v_0} \simeq M_r(D)$ ;

- aux places infinies $\mathbb{D}$ est scind\'ee ;

- \`a toute place $v$ diff\'erente de $v_{0}$ o\`u $\mathbb{D}$ est
ramifi\'ee, ${\mathbb{D}}_v$ est isomorphe \`a une alg\`ebre \`a
division sur ${\mathbb{F}}_v$.\\ 

Soient $v_0,v_1...v_m $ les places de $\mathbb{F}$ o\`u $\mathbb{D}$
est ramifi\'ee. On fixe une fois pour toutes un isomorphisme
${\mathbb{D}}_{v_0} \simeq M_r(D)$ et des isomorphismes  
${\mathbb{D}}_v \simeq M_n({\mathbb{F}}_v)$ pour toutes les places $v$
o\`u $\mathbb{D}$ est scind\'ee. Pour toute place $v$ de $\mathbb{F}$
on note $GL_n({\mathbb{F}}_v)$ par $G_v$ 
et ${\mathbb{D}}^*_v$ par $G'_v$.

        On note indistinctement $Z_v=Z_{GL_n({\mathbb{F}}_v)}\simeq
        Z_{{\mathbb{D}}^*_v}$. Les ad\`eles des groupes $GL_n$ et
        ${\mathbb{D}}^*$ sur $\mathbb{F}$ seront not\'ees  $\
        GL_n(\a(\F))$ et $\ {\mathbb{D^*}}(\a(\F))$.\\ 

       Soit $v_{m+1}$ une place finie de ${\mathbb{F}}$ o\`u
${\mathbb{D}}$ n'est pas scind\'ee.
On pose $S=\{v_0...v_{m+1}\}$. Pour tout $v\in S$ nous d\'efinissons
$\omega_v$ et $\pi_v$ de la fa\c{c}on suivante~:  

- $\omega_{v_0}=\omega$, $\pi_{v_0}=\pi_0$ ;

- $\omega_{v_i}\cong 1$ pour tout $i \in \{1,2...m+1\}$ ;

- $\pi_{v_i}\in E^2(GL_n({\mathbb{F}}_{v_i});\omega_{v_i})$ la
  repr\'esentation de Steinberg de $GL_n({\mathbb{F}}_{v_i})$ pour
  tout $i\in \{1,2...m\}$ ; 

- $\pi_{m+1}$ une repr\'esentation cuspidale de caract\`ere central
  trivial de $GL_n({\mathbb{F}}_{v_{m+1}})$. 

Notons $\tilde{\pi}$ une repr\'esentation automorphe cuspidale qui
v\'erifie $\tilde{\pi}_v\simeq \pi_v$ pour tout $v\in S$ et
$\tilde{\omega}$ son 
caract\`ere central. L'existence d'un tel $\tilde{\pi}$ n'est pas
\'evident. Dans un appendice de [He], Henniart montre ce r\'esultat
pour un groupe 
r\'eductif quelconque, mais \`a condition que les $\pi_v$ soient
cuspidales pour tout $v\in S$. Dans le cas particulier du groupe
lin\'eaire, sa d\'emonstration marche avec {\it \cares} \`a la place de
{\it \cusps}. En effet, il suffit de~:

-  rajouter \`a une autre place une composante
cuspidale, 

- remplacer tout au long de la preuve {\it coefficient} par
{\it pseudocoefficient}~;\\
pour montrer que le psudocoefficients jouent ici exactement le m\^eme
r\^ole que les coefficients dans la d\'emonstration de Henniart,
il
faut utiliser les deux arguments suivants~:

- un
pseudocoefficient d'une 
\care\ $\tau$ est une fonction sur laquelle la trace de toute
repr\'esentation g\'en\'erique non \'equivalente \`a $\tau$ s'annule
(car ou bien de carr\'e int\'egrable, ou bien non elliptique),

- toute composante locale d'une
repr\'esentation automorphe cuspidale de $GL_n$ est g\'en\'erique
([Sh],pp.190).    

Si $\omega_v$ est un caract\`ere de $Z_v$, on note $H(G_v;\omega_v)$
l'espace des fonctions $f$ sur $G_v$ localement constantes \`a support
compact modulo $Z_v$ telles que $f(zg)=\omega_v^{-1}(z)f(g)$ pour tout
$g\in G_v$ et 
tout $z\in Z_v$. Pareil pour $G'_v$. Pour toute place $v$ on pose
$K_v=GL_n(O_v)$ et $K'_v=GL_r(O_{D_v})$. Si $f\in H(G_v;\omega_v)$ et $f'
\in H(G'_v;\omega_v)$ sont \`a support dans l'ensemble des
\'el\'ements semisimples r\'eguliers, on dit que $f$ et $f'$ se
correspondent si 
leurs int\'egrales orbitales sont \'egales sur des \'el\'ements qui se
correspondent ($\forall\ g\lra g',\ \Phi(f;g)=\Phi(f';g')$) et
l'int\'egrale orbitale de $f$ est nulle sur tout \'el\'ement $g$ qui
ne correspond \`a aucun $g'\in G'$. Ici le
choix de mesures pour ces int\'egrales orbitales se fait de la fa\c{c}on
suivante~: sur $G_v=GL_n(\F_v)$ la mesure est telle que
$vol(K_v)=1$, sur $G'_v=GL_r(\D_v)$ la mesure est telle que
$vol(K'_v)=1$, sur $Z_v$ on fixe une mesure $dz$ arbitraire, sur
tout tore elliptique maximal $T$ de $G_v$ ou $G'_v$ on fixe la mesure $dt$
telle que pour la mesure quotient $dt/dz$ on ait $vol(T/Z_v)=1$, sur
tous les tores maximaux de $G_v$ on fixe des mesures arbitraires avec
la seule
condition que si deux tores sont conjugu\'es les mesures choisies se
correspondent via cet isomorphisme (c'est ind\'ependant de la conjugaison
envisag\'ee), et sur les tores maximaux de $G'_v$ on fixe des mesures
provenant des tores maximaux de $G_v$ comme dans l'Annexe. Remarquons
que les fonctions qui appara\^\i tront sont \`a support dans les
\'el\'ements semisimples r\'eguliers, et ne ``voyent'' pas les autres
points. En particulier, bien que l'on n'ait pas pos\'e de condition
sur la caract\'eristique du corps de base, pour toute fonction $f'
\in H(G'_v;\omega_v)$, il existe une fonction $f\in H(G_v;\omega_v)$
qui lui correspond et r\'eciproquement, par le th\'eor\`eme de
submersion de Harish-Chandra ([H-C]).

Soit $\tilde{\omega}$ un caract\`ere unitaire de $\F^*$. On pose
$\tilde{\omega}_v= \omega_v$ et on d\'efinit
$H(GL_n({\mathbb{F}});\tilde{\omega}))$  comme l'ensemble des fonctions
$f:GL_n(\a(\F))\to {\mathbb C}$ invariantes \`a gauche par $GL_n(\F)$, qui sont
produit sur l'ensemble des places $v$ de fonctions locales 
$f_v\in H(G_v;\omega_v)$ telles que, pour presque tout $v$, le support
de $f_v$ est inclus dans $Z_vK_v$ et  
$f_v(k)=1$ pour tout $k\in K_v$. On d\'efinit de fa\c{c}on analogue
$H(\D^*({\mathbb{F}});\tilde{\omega}))$.

Soient $\tilde{f} \in
      H(GL_n({\mathbb{F}});\tilde{\omega}))$  et 
      $\tilde{f'}\in H({\mathbb{D}}^*;\tilde{\omega})$. On
      dit que $\tilde{f}$ et $\tilde{f'}$ se correspondent et on
      \'ecrit $\tilde{f}\lra \tilde{f'}$ si pour tout $i\in
      \{0,1,2...m\}$ on a $\tilde{f}_{v_i}\lra \tilde{f'}_{v_i}$ et
      pour toute place $v$ o\`u ${\mathbb{D}}$ est scind\'ee on
      a $\tilde{f}_v=\tilde{f'}_v$ (via les isomorphismes \`a ces
      places fix\'es au d\'ebut). 

On note $\rho_0$ (resp. $\rho'_0$) les repr\'esentations de
$GL_n(\a(\F))$ (resp. 
$\D^*(\a(\F))$) dans l'espace des
formes automorphes cuspidales de caract\`ere
central $\tilde{\omega}$. \\ 
\begin{prop}\label{11}    Si $\tilde{f}\in
H(GL_n({\mathbb{F}});\tilde{\omega})$ et $\tilde{f'}\in
H({\mathbb{D}}^*;\tilde{\omega})$ sont  telles que $\tilde{f} \lra
\tilde{f'}$, et si $\tilde{f}_{v_{m+1}}$ est un coefficient d'une
repr\'esentation cuspidale de $G_{v_{m+1}}$, alors on a~:  
$$tr\rho_0(\tilde{f})=tr\rho'_0(\tilde{f'}).$$\\ 
\end{prop}

 \begin{prop}\label{12} On pose
$$V=\{v_0,v_1...v_m\}$$
 Posons $G_V=\Pi_{v\in V}G_v$, $G'_V=\Pi_{v\in V}G'_v$ et
$\omega_V=\Pi_{v\in V} \tilde{\omega}_v$. Notons $\tilde{\pi}_V$ la
repr\'esentation de $G_V$ induite par $\tilde{\pi}$ par restriction
aux places dans $V$. Si $f_V\in
\prod_{i=0}^m H(G_{v_i};\omega_{v_i})$ et $f'_V\in \prod_{i=0}^m
H(G'_{v_i};\omega_{v_i})$, on \'ecrit $f_V\lra f'_V$ si pour tout
$i\in \{1,2...m\}$ $f_{v_i}$ et $f'_{v_i}$ sont \`a support dans les
\'el\'ements elliptiques r\'eguliers et pour tout $i\in \{0,1,2...m\}$
on a $f_{v_i}\lra f'_{v_{i}}$. 
    On a alors~:
          $$tr\tilde{\pi}_{V}(\tilde{f}_V)=\sum_{\tilde{\pi}'\in U'}
m(\tilde{\pi}')tr\tilde{\pi}'_V(\tilde{f}'_V)$$ 
   o\`u $U'$ est l'ensemble des repr\'esentations automorphes
cuspidales $\tilde{\pi}'$ de ${\mathbb{D}}^*$ telles que pour tout
$v\notin V$ on ait $\tilde{\pi}'_v=\tilde{\pi}_v$, $m(\tilde{\pi}')$ est
la multiplicit\'e de $\tilde{\pi}'$ dans $\rho'_0$ et l'indice $V$
veut dire ``restriction aux places dans $V$''.\\ 
\end{prop}  

\begin{prop}\label{finiG'} 
L'ensemble $U'$ est fini.
\end{prop}

\begin{prop}\label{13} Posons $V=\{v_0,v_1...v_m\}$ comme plus
haut. Il existe alors un nombre fini $k$ d'entiers strictement
positifs $a_j$, $1\leq j\leq k$, et des repr\'esentations
irr\'eductibles $\pi'_{Vj}$ de $G'_V$ tel 
que, pour tout $f_V \in \prod_{i=0}^m H(G_{v_i};\omega_{v_i})$  et pour
tout $f'_V\in \prod_{i=0}^m H(G'_{v_i};\omega_{v_i})$  telles que
$f'_v\lra f_v$ pour
tout $v\in V$, on ait~:
 
$$tr\tilde{\pi}_V(f_V)=\sum_{j=1}^{k} a_j tr\pi'_{Vj}(f'_V)$$
\end{prop}

La prop.\ref{11} est une application classique de la formule
des traces, qui remonte \`a [JL]. Le r\'esultat tel qu'\'enonc\'e ici
est une cons\'equence de la formule des traces simple, valable
en toute caract\'eristique. Pour la preuve de \ref{12} voir par exemple 
[Fla] - c'est ind\'ependant de la caract\'eristique du corps de
base. La prop.\ref{finiG'} est prouv\'ee dans [Ba3] en 
caract\'eristique nulle et dans la section 3, th.\ref{finitude} du
pr\'esent article, en 
caract\'eristique non nulle. La prop.\ref{13}
d\'ecoule imm\'ediatement de \ref{12} et \ref{finiG'}.\\ 

\begin{prop}\label{14}  Il existe un nombre fini $k'$ d'entiers
strictement positifs $a_p$, $1\leq p\leq k'$, et de repr\'esentations
irr\'eductibles $\pi'_p$ de $G'$ tels qu'on ait~:  
\begin{equation}\label{egalitecareint}
\chi_{\pi_0}(g)=(-1)^{n-r} \sum_{p=1}^{k'} a_p \chi_{\pi'_p}(g')\ \
\ \ \forall g\lra g'.
\end{equation} 

\end{prop}
\ \\
{\bf{D\'emonstration.}}  Si $g'_V$ est un \'el\'ement de $G'_V$ et
$g_V$ est un \'el\'ement de $G_V$ qui lui correspond composante par
composante (en particulier 
les composantes de $g'_V$ sont des \ssrs\ des $G'_{v_i}$), alors on
\'ecrit $g_V\lra g'_V$. Comme
dans l'\'egalit\'e de la prop.\ref{13} il 
appara\^\i t un nombre fini de repr\'esentations, et que les
caract\`eres de ces repr\'esentations sont constantes au voisinage de
$g'_V$ et de $g_V$ respectivement, en choisissant des fonctions \`a
petit support au voisinage de ces \'el\'ements (on peut les choisir
telles qu'elles se correspondent, par le principe de submersion de
Harish-Chandra) on peut passer \`a une \'egalit\'e des
caract\`eres fonctions. En mettant alors tout du 
c\^ot\'e gauche de l'\'egalit\'e on obtient  
$$\chi_{\tilde{\pi}_V}(g_V)-\sum_{j=1}^{k} a_j \chi_{\pi'_{Vj}}(g'_V)=0\ \
\ \ \forall g_V\lra g'_V.$$  
Les composantes locales de $\tilde{\pi}_V$ aux places $v_1,v_2...v_m$
sont des repr\'esentations de Steinberg. Le caract\`ere de la
repr\'esentation de Steinberg de $G_{v_i}$ correspond par la
correspondance avec une alg\`ebre \`a division au caract\`ere de la
repr\'esentation triviale de $G'_{v_i}$. On utilise alors  
l'ind\'ependance lin\'eaire des caract\`eres sur l'ensemble des
\'el\'ements r\'eguliers de $\Pi_{i=1}^m
{\mathbb{D}}^*_{v_i}$ qui est compact modulo le centre. La nullit\'e
du coefficient du caract\`ere de la repr\'esentation triviale de ce
groupe donne la relation voulue sur les caract\`eres fonction \`a la
place $v_0$. Le
$(-1)^{n-r}$ vient apr\`es un simple calcul de la loi de
r\'eciprocit\'e et du fait que le caract\`ere de la repr\'esentation
de Steinberg de $GL_n(F_{v_i})$ (pour $i$ de $1$ \`a $m$) est \'egal
\`a $(-1)^{n-1}$ sur 
l'ensemble des \'el\'ements elliptiques r\'eguliers.\qed

\ \\
\begin{prop}\label{careint} Les repr\'esentations $\pi'_{p}$ qui
apparaissent dans 
l'\'egalit\'e de la proposition \ref{14} sont de carr\'e
int\'egrable.\\ 
\end{prop}
\ \\

Avant de prouver cette proposition, ouvrons une parenth\`ese o\`u nous
montrons 
deux propositions, A et B. La  prop.B 
s'applique d'une part pour prouver la proposition \ref{careint} plus
haut et d'autre part dans l'\'etude de la compatibilit\'e de la
correspondance de Jacquet-Langlands avec les foncteurs de
Jacquet. L'id\'ee de la 
d\'emonstration de la prop.B et les arguments sont de
[DKV], B.2.e. Nous \'evitons ici l'utilisation du transfert des
int\'egrales orbitales pos\'e par r\'ecurrence dans 
[DKV], et qui est d\'elicat en caract\'eristique non nulle. 
La prop.A n'est qu'un r\'esultat interm\'ediaire, plus
faible que la prop.B, et qui ne sert qu'\`a d\'emontrer cette derni\`ere.    

Soient $\pi_i,\  i\in\{1,2...k\}$ des \rlis\ de $G$ et $\pi'_j,\ j\in
\{1,2...k'\}$ des \rlis\ de $G'$, et supposons qu'on ait la relation~:

$${\bf (E_G)}\ \ \ \ (\sum_{i=1}^k
a_i\chi_{\pi_i})(g)=(-1)^{n-r}(\sum_{j=1}^{k'} 
a'_j\chi_{\pi'_j})(g')\ \ \ \text{pour tout}\ g\in G \lra g'\in G'$$ 
o\`u les $a_i$ et les $a'_j$ sont des nombres complexes. Sur $G$
(resp. $G'$) on fixe la paire parabolique minimale standard $(A;P)$
(resp. $(A';P')$) o\`u $A$ (resp. $A'$) est le tore diagonal et $P$
(resp. $P'$) est le groupe des matrices triangulaires sup\'erieures
inversibles. Donc, si $L$ est un \sgls\ de $G$, $L$ est le groupe des
matrices inversibles diagonales par blocs d'une taille donn\'ee. Il
correspond de fa\c{c}on biunivoque \`a une suite finie d'entiers
strictement positifs $(n_1;n_2;...n_p)$ o\`u  
$$n=n_1+n_2+...+n_p$$et les $n_i$ repr\'esentent les tailles des dits
blocs dans l'ordre, en lisant du haut \`a gauche vers le bas \`a
droite. Pareillement, \`a un \sgls\ $L'$ de $G'$ correspond de fa\c{c}on
biunivoque une suite finie d'entiers strictement positifs
$(n'_1;n'_2;...n'_p)$ o\`u  
$$r=n'_1+n'_2+...+n'_p.$$ 
On dit alors que $L$ {\it se transf\`ere} si, pour tout $i\in \{1,2,
...p\}$, $d$ divise $n_i$. Soit $L'$ le \sgls\ de $G'$ qui correspond
\`a la suite $(n'_1;n'_2;...n'_p)$ telle que, pour tout $i\in \{1,2,
...p\}$, $n'_i=n_i/d.$ On dit alors que $L'$ correspond \`a $L$ ou que
$L$ correspond \`a $L'$ ou que $L$ et $L'$ se correspondent. 
Si $P$ est un sous-groupe parabolique standard de $G$ et
$P=LU$ est une d\'ecomposition de Levi standard de $P$, on dit que $P$
{\it se transf\`ere} si $L$ se transf\`ere. Alors, si $L'$ est le
\sgls\ de $G'$ qui correspond \`a $L$, et $P'$ le sous-groupe
parabolique standard de $G'$ qui a pour sous-groupe de Levi standard
$L'$, on dit que $P$ et $P'$ se correspondent. Si $L$ et $L'$ sont deux
sous-groupes de Levi standard de $G$ et $G'$ respectivement, on
utilise la notation $L\lra L'$ pour  dire que $L$ et $L'$ se
correspondent. On adopte la m\^eme notation pour des sous-groupes
paraboliques standard qui se correspondent.   

Soit $P=LU$ un sous-groupe parabolique propre standard de $G$ qui se
transf\`ere et soit $P'=L'U'$ le sous-groupe parabolique standard  
de $G'$ qui lui correspond. On se demande si on a alors la relation
:$${\bf (E_P)}\ \ \ \ (\sum_{i=1}^k a_i\chi_{res_P^G
\pi_i})(l)=(-1)^{n-r}(\sum_{j=1}^{k'} a'_j \chi_{res_{P'}^{G'}
\pi'_j})(l')\ \ {\rm{pour} \ \rm{tout}} \ l\in L \lra l'\in  L'$$
($res_P^G$ est la restriction parabloique, ou ``foncteur de Jacquet''). 
La r\'eponse est oui, sous certaines conditions (voir prop.B), et non
en g\'en\'eral. 

Soient $Z_L$ et $Z_{L'}$ les centres de $L$ et de $L'$ respectivement. On
note $X(Z_L)$ l'ensemble des caract\`eres lisses (pas forc\'ement
unitaires) de $Z_L$ et pareil pour $Z_{L'}$. Il y a un isomorphisme
naturel entre $Z_L$ et $Z_{L'}$ et on identifiera tacitement par la
suite, via cet isomorphisme, $Z_L$ et $Z_{L'}$ ainsi que $X(Z_L)$ et
$X(Z_{L'})$. Soit 
$A_P$ l'ensemble des $\omega \in X(Z_L)$ tels que $\omega$ soit un
exposant central de l'un 
des $\pi_i$ relatif \`a $P$ (c'est-\`a-dire que $\omega $ est le
caract\`ere central d'un des sous-quotients ir\'eductibles de la
restriction parabolique de $\pi_i$ \`a $P$). Soit 
$A'_{P'}$ l'ensemble des $\omega\in X(Z_{L'})=X(Z_L)$ tels que
$\omega$ soit un exposant 
central de l'un des $\pi'_j$ relatif \`a $P'$.\\ 
\\
{\bf{PROPOSITION A.}}  {\it Si $({\bf{E_G}})$ est v\'erifi\'ee et si
les conditions (1), (2) et (3) plus bas sont satisfaites~:

(1) tous les $a_i$ sont des nombres r\'eels non nuls de m\^eme signe et
  tous les $a'_j$ sont des nombres r\'eels non nuls de m\^eme signe, 

(2) pour tout $i\in \{1,2...k\}$, ou bien $res_P^G \pi_i$ est nulle,
  ou bien tout sous-quotient irr\'eductible de
  $res_P^G \pi_i$ est une \eci,\   

(3) pour tout $j\in \{1,2...k'\}$, ou bien $res_{P'}^{G'}\pi'_j$ est
nulle, ou bien tout sous-quotient irr\'eductible de
$res_{P'}^{G'}\pi'_j$ est une \eci,\ \\ 
alors on a~:

i)$A_P$=$A'_{P'}$

ii)$(\bf{E_P})$ est v\'erifi\'ee.\\}
\ \\
{\bf{D\'emonstration.}} (i) Tout d'abord on a que pour tout $n'$ le
caract\`ere d'une \eci\ de $GL_{n'}(F)$ est constant, 
r\'eel, non nul et de signe $(-1)^{n'-1}$ sur les \'el\'ements elliptiques
r\'eguliers d'un voisinage de l'unit\'e  
dans $GL_{n'}(F)$. Cela est une cons\'equence imm\'ediate de la
correspondance entre $GL_n$ et une alg\`ebre \`a division ([Ro] et [Ba2]), 
pour ne pas citer d'autres r\'esultats plus g\'en\'eraux sur les
germes, mais qui ne sont \'ecrits qu'en caract\'eristique nulle.  
Mais alors on a, par l'hypoth\`ese de r\'ecurrence ${\bf H}_{r'}$, que,
pour tout $r'<r$, le caract\`ere d'une \eci\ de
$GL_{r'}(D)$ est constant, r\'eel, non nul et de signe $(-1)^{r'-1}$ sur les
\'el\'ements elliptiques r\'eguliers d'un voisinage de l'unit\'e dans
$GL_{r'}(D)$. Cela
s'applique en 
particulier \`a $L$ et $L'$ qui sont des produits de tels
groupes. Faisons la remarque que {\underline{\it c'est le seul 
endroit o\`u}} {\underline{\it nous avons besoin de l'hypoth\`ese de
r\'ecurrence}}. Il 
existe donc un voisinage $V$ de 
l'unit\'e dans $L$ tel que le caract\`ere de tout sous-quotient
irr\'eductible de $res_P^G \pi_i$ soit constant r\'eel de signe
$(-1)^{n-1}$ sur l'ensemble $V_e$ des \'el\'ements elliptiques
r\'eguliers de $V$ pour tout $i\in \{1,2...k\}$ et il existe un
voisinage de l'unit\'e  
$V'$ dans $L'$ tel que le caract\`ere de tout sous-quotient
irr\'eductible de $res_{P'}^{G'} \pi'_j$ soit constant r\'eel de signe
$(-1)^{r-1}$ sur l'ensemble ${V'}_e$ des \'el\'ements elliptiques
r\'eguliers de $V'$ pour tout $j\in \{1,2...k'\}$.

Soit $g'\in V'_e$ tel qu'il existe $g\in V$ avec la propri\'et\'e 
$g\lra g'$. On a \'evidemment $g\in V_e$. Si $L$ correspond \`a la
partition ordonn\'ee $(n_1;n_2;...n_p)$ de $n$, alors $g$ se
repr\'esente par un $p$-uple $(g_1;g_2;...g_p)$, o\`u $g_i\in
GL_{n_i}(F)$. Posons 
$$c(g)=\inf_{1\leq i\leq p-1}\frac{|det(g_i)|_F}{|det(g_{i+1})|_F}.$$   

Soit maintenant $z\in Z_L=Z_{L'}$. On repr\'esente $z$ par un $p$-uple
$(z_1;z_2;...z_p)$, o\`u $z_i\in F^*$. On pose
$$N(z)=\sup_{1\leq i\leq p-1}\frac{|z_i|_F}{|z_{i+1}|_F}.$$
Alors, si $N(z)<c(g)^{-1}$, on a $c(zg)<1$ et, par le th.5.2 de [Ca1],
pour toute repr\'esentation admissible $\pi$ de $G$ (ou $G'$) 
$$ \chi_{res_P^G\pi}(zg)=\chi_{\pi}(zg)\ \ \ \ \ \ 
(\text{ou}\ \chi_{res_{P'}^{G'}\pi}(zg')=\chi_{\pi}(zg')).$$ 

\'Ecrivons maintenant dans le groupe de Grothendieck des
repr\'esentations lisses de longueur finie de $G$ et respectivement $G'$~:   
$$res_P^G\pi_i = \sum_{s=1}^{k_i} \alpha _s \pi _i^s, \ \ \ \ \alpha_s
>0$$
et
$$ res_{P'}^{G'}\pi '_j = \sum_{t=1}^{k_j} \alpha'_t \pi_j^{'t},
\ \ \ \ \alpha'_t >0,$$ 
o\`u les $\pi _i^s$ et les $\pi_j^{'t}$ sont des repr\'esentations
irr\'eductibles. 

On a alors, en passant aux caract\`eres, et en utilisant la relation
$(E_G)$ pour $g$ et $z$ comme plus haut~:

$$\sum_{i=1}^k a_i \sum_{s=1}^{k_i}\alpha_s
\chi_{\pi_i^s}(zg)=(-1)^{n-r}\sum_{j=1}^{k'} a'_j \sum_{t=1}^{k'_j}
\alpha'_t \chi_{\pi_j^{'t}}(zg').$$
Si $\omega_i^s$ sont les
caract\`eres centraux des $\pi_i^s$ et $\omega_j^t$ sont les
caract\`eres centraux des $\pi_j^{'t}$, alors on a l'\'egalit\'e~: 

$$\sum_{i=1}^k a_i \sum_{s=1}^{k_i} \alpha_s \omega_i^s (z)
\chi_{\pi_i^s} (g) = (-1)^{n-r}\sum_{j=1}^{k'} a'_j \sum
_{t=1}^{k'_j}\alpha'_t \omega_j^t (z) \chi_{\pi_j^{'t}}(g').$$ 

On obtient, en regroupant, une relation~:\\

\begin{equation}\label{acum1}
\sum _{\omega \in  A_P} n_{\omega} \omega (z) = \sum_{\omega' \in
A'_{P'}} n_{\omega'} \omega' (z)
\end{equation}
o\`u, chose tr\`es importante, les $n_{\omega}$ et les $n_{\omega'}$
sont tous non nuls comme somme de nombres r\'eels non nuls et de
m\^eme signe (par exemple, $n_{\omega}=\sum a_i\alpha_s\chi_{\pi_i^s}
(g)$, o\`u la somme porte sur les couples $(i,s)$ tels que
$\omega_i^s=\omega$).  
Or, cette relation est vraie pour tout $z$ tel que $N(z)<c(g)^{-1}$. Pour
avoir $A_P = A'_{P'}$ le lemme suivant suffit~:\\ 
\ \\
\begin{lemme}\label{lemme}
Si $\omega_1, \omega_2...\omega_v$ sont des caract\`eres distincts de $Z_L$ et $a_1,a_2...a_v$ sont des nombres complexes tels qu'on ait~:
$$\forall z\in Z_L\ \text{tel que}\ N(z)<c,\ \ \ \Sigma_{i=1}^v a_i
\omega_i(z)=0,$$ 
alors on a $a_i=0$ pour tout $i$.
\end{lemme}
\ \\
{\bf D\'emonstration du lemme \ref{lemme}.}
Raisonnons par l'absurde et supposons $v$ 
minimal tel que \\ 
 
\begin{equation}\label{acum2}
\sum_{i=1}^v a_i\omega_i(z)=0
\end{equation} 
pour tout $z\in Z_L$ tel
que $N(z)<c$ et il existe au moins un $ a_i$ non nul. 
 
\'Evidemment $v\geq 2$ et tous les $a_i$ sont non nuls. Soit
$z_0\in Z_L$ tel que $N(z_0) <1$. Alors pour tout $z$ tel que $N(z)<c$
on a $N(z_0z)<c$ donc $\Sigma_{i=1}^v a_i \omega_i(z_0 z)=0$ ce qui
donne  
$\Sigma_{i=1}^v a_i \omega_i (z_0) \omega_i (z)=0$. En multipliant
\ref{acum2} par $\omega_1(z_0)$ et en faisant la diff\'erence avec
la derni\`ere relation obtenue on trouve une relation du type
\ref{acum2}  avec un $v$ strictement inf\'erieur donc une relation
dans laquelle tous les coefficients sont nuls. Ceci implique que pour
tout $i\in \{1,2..v\}$ on a $\omega_i(z_0)=\omega_1(z_0)$. Mais alors
on a en particulier~: 
$$\omega_1(z_0)=\omega_2(z_0) \ \ \ \ \ \forall z_0\ \text{tel que}\
N(z_0)<1$$ 
et aussi, les $\omega_i$ \'etant des caract\`eres, 
$$\omega_1(z_0^{-1})=\omega_2(z_0^{-1}) \ \ \ \ \ \forall z_0\
\text{tel que}\ N(z_0)<1$$ 

Comme tout \'el\'ement $h\in Z_L $ s'\'ecrit $h=xy^{-1}$ o\`u $N(x)<1$
et $N(y)<1$ (prendre $y$ tel qu'on ait simultan\'ement $N(y)<1$ et
$N(hy)<1$, et poser $x=hy$) on aboutit \`a $\omega_1=\omega_2 $ ce
qui contredit nos hypoth\`eses.
\qed \ \\ 
 
Ainsi le point (i) de la proposition est d\'emontr\'e.\\

 (ii)
 Prenons maintenant $g'\in L^{'reg}$ quelconque et $g\lra g'$. On peut
 raisonner de la m\^eme fa\c{c}on, et il est
 toujours vrai 
que pour tout $z\in Z_L$ qui v\'erifie $N(z) < c$ on a une relation du
 type  
\ref{acum1} avec la seule diff\'erence qu'on ne puisse plus garantir 
 le fait que les coefficients qui apparaissent soient tous 
 positifs. On \'ecrit cette relation~: 

$$\sum _{\omega \in A_P} n_{\omega} \omega(z) = \sum_{\omega \in
A'_{P'}} n'_{\omega} \omega(z).$$ 
\ \\
Par le point (i) on a $A_P=A'_{P'}$ et par le lemme plus haut on en d\'eduit que $n_{\omega}=n'_{\omega}$ pour tout $\omega \in A_P=A'_{P'}$. En particulier $\sum_{\omega}n_{\omega}=\sum_{\omega} n'_{\omega}$ et cette relation,
si on regarde qui \'etaient les coefficients $n_{\omega} $ et
$n'_{\omega}$, n'est autre que la relation $\bf{(E_P)}$ appliqu\'ee
\`a $g\lra g'$. Le point (ii) est d\'emontr\'e.
\qed  
\\
\ \\
{\bf{PROPOSITION B.}}\label{propB} {\it Si $\bf{(E_G)}$ est v\'erifi\'ee et si

-  tous les $a_i$ sont des nombres r\'eels non nuls de m\^eme signe et
   tous les $a'_j$ sont des nombres r\'eels non nuls de m\^eme signe
   et 

-  pour tout $i\in \{1,2...k\}$, $\pi_i$ est une \eci, \\ 
alors~:

(i) $A_P=A'_{P'}$ pour tout $P\lra P'$

(ii) $\bf{(E_P)}$ est v\'erifi\'ee  pour tout $P$ qui se transf\`ere.

(iii) Tous les $\pi'_j$ sont des \ecis.}\\
\ \\
{\bf{D\'emonstration.}} La d\'emonstration utilise plusieurs fois le crit\`ere
de Casselman ([Ca2] 4.4.6 et 6.5.1) et la prop.A plus haut. Pour
un \sgls\ $L$ de $G$ comme plus haut, le centre $Z_L$ de $L$
s'identifie \`a $(F^*)^p$ et tout caract\`ere $\omega$ de $Z_L$ 
correspond \`a une suite ordonn\'ee
$(\alpha_1(\omega);\alpha_2(\omega);...\alpha_p(\omega))$ 
de nombre complexes, par la formule~:
$$\omega ((z_1;z_2;...z_p))=\prod_{1\leq i\leq p} |z_i|_F^{\alpha_i(\omega)}\ \
\ \ \ \forall (z_1;z_2;...z_p)\in Z_L.$$
Soit $\pi$ une repr\'esentation lisse irr\'eductible de $G$. Le
crit\`ere de Casselman dit que $\pi$ est une \eci\ si
et seulement si elle a la propri\'et\'e~:

{\bf (Cas)} pour tout sous-groupe parabolique standard $P$, ou
bien $res_P^G\pi$ est nulle, ou bien, pour tout sous-quotient
irr\'eductible $\tau$ de $res_P^G\pi$, le caract\`ere central $\omega$
de $\tau$ a la propri\'et\'e~:
$$\forall i,\ 1\leq i\leq p-1, Re(\alpha_i(\omega))<
Re(\alpha_{i+1}(\omega)).$$  
Ce
crit\`ere implique d\'ej\`a~:

{\bf (***)} la restriction d'une
\eci\ \`a un sous-groupe parabolique standard est ou bien nulle, ou bien tous
ses  sous-quotients irr\'eductibles sont des repr\'esentations
essentiellement de carr\'e int\'egrable. 

Le
crit\`ere de Casselman admet aussi une variante~:  $\pi$ est une \eci\
de $G$ si 
et seulement si elle a la propri\'et\'e~:

{\bf (Cas0)} si $P$ est le sous-groupe parabolique standard de $G$ tel que
$res_P^G\pi$ soit cuspidale, alors pour tout sous-quotient
irr\'eductible $\tau$ de $res_P^G\pi$, le caract\`ere central $\omega$
de $\tau$ a la propri\'et\'e~:
$$\forall i,\ 1\leq i\leq p-1,\
Re(\alpha_i(\omega))<  Re(\alpha_{i+1}(\omega)).$$ 

Ces m\^emes r\'esultats valent pour $G'$ aussi.

Maintenant, pour chaque repr\'esentation $\pi'_j$, il existe un  sous-groupe
parabolique standard $P'_j$ de $G'$ tel que la restriction de $\pi'_j$
\`a $P'_j $ soit cuspidale. Soit $B_0$ l'ensemble des sous-groupes
paraboliques standard $P'$ de $G'$  tels que pour au moins un $j$ la
restriction de $\pi'_j$ \`a $P'$  soit cuspidale. D\'efinissons une
relation d'ordre partielle sur $B_0$ en posant $P'<Q'$ si et
seulement si $P'\subset Q'$. Pour tout $P'\in B_0$, on note $\Pi_{P'}$
l'ensemble des repr\'esentations $\pi'_j$ dont la restriction \`a $P'$
est cuspidale.\\

Soit $P'_0$ un \'el\'ement minimal de $B_0$. Soit $P_0 \lra
P'_0$. Comme $P'_0$ est minimal, pour tout $j\in
\{1,2...k'\}$, la restriction de $\pi'_j$ \`a $P'_0$ est nulle si
$\pi_j'\notin \Pi_{P'_0}$, et cuspidale si $\pi_j'\in \Pi_{P'_0}$.
Par cons\'equent, pour tout $j$, ou bien la restriction de 
$\pi'_j$ \`a $P'_0$ est nulle, ou bien tous ses sous-quotients
irr\'eductibles sont des \ecis. Il en est de m\^eme pour les
restrictions des $\pi_i$ \`a $P_0$ par {\bf {(***)}}. On peut donc  
 appliquer le point (i) de la proposition A pour en d\'eduire que
 $A_{P_0}=A'_{P'_0}$. Mais alors, {\bf{(Cas)}} appliqu\'e aux $\pi_i$
implique {\bf{(Cas0)}} pour les $\pi_j'\in \Pi_{P'_0}$. Donc tous les
\'el\'ements de $\Pi_{P'_0}$ sont des \ecis.\\  
 
On montre alors le point (iii) de notre proposition par r\'ecurrence~: posons 
$B_1=B_0\backslash \{P'_0\}$. Prenons un \'el\'ement minimal $P'_1$ de $B_1$.
Soit $P_1\lra P'_1$. Soit $j\in \{1,2...k'\}$. Alors on a trois
possibilit\'es~:

- ou bien la restriction de $\pi'_j$ \`a $P'_1$
est nulle, 

- ou bien $\pi'_j\in \Pi_{P'_1}$, 

- ou bien  $\pi'_j\in \Pi_{P'_0}$.\\
Dans le deuxi\`eme cas, les sous-quotients irr\'eductibles de
$res_{P'_1}^{G'}\pi'_j$ sont des \ecis\ parce que cuspidales, et dans
le troisi\`eme cas par {\bf (***)} et par le pas de r\'ecurrence pr\'ec\'edent.
En conclusion, si la restriction de $\pi'_j$ \`a $P'_1 
$ est non nulle, alors tous les sous-quotients irr\'eductibles
de $res_{P'_1}^G\pi'_j$ \`a $P'_1 $ sont des
\ecis. Comme du c\^ot\'e de $G$ c'est pareil (par {\bf(***)}) on
applique encore une fois le point (i) de la 
proposition pr\'ec\'edente pour avoir que $A_{P_1}=A'_{P'_1}$, et en
d\'eduire que tous les \'el\'ements de $\Pi_{P'_1}$ sont des \ecis.
Ainsi de suite ; \`a chaque 
pas, les restrictions non nulles du c\^ot\'e de $G$ sont constitu\'ees
uniquement de \ecis\ par {\bf(***)}, et du
c\^ot\'e de $G'$ les restrictions non nulles sont constitu\'es
uniquement de repr\'esentations cuspidales et de \ecis\ (par {\bf(***)}
et les pas pr\'ec\'edents). 
Apr\`es 
$card(B_0)$ pas, on est s\^ur de trouver que toutes les repr\'esentations
$\pi'_j$ 
sont essentiellement de carr\'e int\'egrable et le point (iii) est 
d\'emontr\'e. \\

Le point (iii) implique directement les points (i) et (ii) par
application de la prop.A et par {\bf(***)}. La proposition B est d\'emontr\'ee.
\qed \ \\ 
\ \\
{\bf D\'emonstration de la proposition \ref{careint}.} C'est imm\'ediat
par la prop.B (iii).
\qed 
\ \\

Soit $\omega$ le caract\`ere central de $\pi_0$. Dans l'annexe nous
d\'efinissons l'espace $L_0(G_e;\omega)$ (resp. $L_0(G'_e;\omega)$)
comme espace des fonctions 
localement constantes sur l'ensemble des \'el\'ements elliptiques
r\'eguliers de $G$ (resp. $G'$), invariantes par conjungaison par des
\'el\'ements 
de $G$ (resp $G'$) et de caract\`ere central $\omega$. \`A partir de
la bijection entre l'ensemble des classes de conjugaison des
\'el\'ements elliptiques r\'eguliers sur les deux groupes nous y
d\'efinissons une bijection $i$ de $L_0(G_e;\omega)$ sur
$L_0(G'_e;\omega)$. Il est bien connu qu'on peut d\'efinir des
sous-espaces pr\'ehilbertiens $L^2(G_e;\omega)$ et
$L^2(G'_e;\omega)$ de $L_0(G_e;\omega)$ et respectivement
$L_0(G'_e;\omega)$ et que la restriction de $i$ induit une isom\'etrie
entre ces deux espaces. Tous ces faits sont rappel\'es dans l'annexe.
Les repr\'esentations
du c\^ot\'e droit de l'\'egalit\'e \ref{egalitecareint} sont
forc\'ement de caract\`ere central $\omega$ et on peut donc \'ecrire une 
\'egalit\'e dans $L_0(G'_e;\omega)$~:
\begin{equation}\label{egalitecareint1}
i(\chi_{\pi_0})=(-1)^{n-r} \sum_{p=1}^{k'} a_p \chi_{\pi'_p}.
\end{equation}
\ \\

{\underline{Supposons maintenant que le th.\ref{ortogonalitate} est
vrai.}}  Alors les objets intervenant dans l'\'egalit\'e 
\ref{egalitecareint1} sont des \'el\'ements de 
$L^2(G'_e;\omega)$ ; comme $i(\chi_{\pi_0})$ est de norme $1$ et
les $\chi_{\pi'_p}$ sont orthogonaux deux \`a deux et de norme $1$, on
obtient 
$$1=\sum_{p=1}^{k'} a_p^2.$$
Comme les $a_p$ sont des entiers positifs, on en d\'eduit que $k'=1$
et $a_1=1$ dans la prop.\ref{14}. On a donc
$$\forall g\lra g'\ \ \ \ \ \ \ \
\chi_{\pi_0}(g)=(-1)^{n-r}\chi_{\pi'_1}(g').$$
De plus, $\pi'_1$ est la seule \eci\ de $G'$ avec cette
propri\'et\'e. En effet, si $\pi'_2$ en \'etait une autre, il y aurait
une contradiction \'evidente avec l'orthogonalit\'e des caract\`eres
sur $G'$ (on a suppos\'e que le th.\ref{ortogonalitate} \'etait
 vrai). 
On pose 
$${\bf C}(\pi_0)=\pi'_1.$$
Ainsi, pour toute $\pi\in E^2(G)$ {\it unitaire} (pour l'instant) il
existe une unique $\pi'\in E^2(G')$ telle qu'on ait      
$$\forall g\lra g'\ \ \ \ \ \ \ \
\chi_{\pi}(g)=(-1)^{n-r}\chi_{\pi'}(g').$$
On pose alors  
$${\bf C}(\pi)=\pi'.$$
Si $\pi\in E^2(G)$ n'est pas unitaire, alors on \'ecrit~:
$$\pi=|det(.)|_F^x\otimes\pi^u$$
o\`u $\pi^u$ est unitaire et $x$ est un nombre r\'eel, et on
pose~:
$${\bf C}(\pi)=|nrd_{G'/F^*}(.)|_F^x{\bf C}(\pi^u)$$
o\`u $nrd$ est la norme r\'eduite. Ainsi nous avons d\'efini une
application $\bf C$ qui v\'erifie l'\'egalit\'e \ref{forte}.\\ 
\ \\
{\it Unicit\'e}\\
Si on
suppose qu'il n'y a pas unicit\'e d'une telle application ${\bf C}$ on
trouve deux \ecis\ de $G'$ qui ont des caract\`eres \'egaux,
bien qu'elles soient non \'equivalentes. En particulier, elles ont des
caract\`eres centraux \'egaux et donc, en tordant par un caract\`ere de $G'$
on peut les rendre unitaires toutes les deux et obtenir deux
repr\'esentations de carr\'e int\'egrable non \'equivalentes, mais dont les
caract\`eres sont \'egaux. Cela contredit
l'orthogonalit\'e des 
caract\`eres sur $G'$.\\
\ \\
{\it Injectivit\'e}\\
Si ${\bf C}(\pi)={\bf C}(\tau)$ alors les caract\`eres de
$\pi$ et $\tau$ sont \'egaux sur l'ensemble des \'el\'ements semisimples
r\'eguliers de $G$ {\it qui ont un
correspondant sur} $G'$. En particulier il y a \'egalit\'e sur
l'ensemble des 
\'el\'ements elliptiques  
r\'eguliers de $G$, car ils ont tous un correspondant sur $G'$. Par
orthogonalit\'e des caract\`eres 
sur $G$, on doit avoir $\pi=\tau$.\\   
\ \\ 
{\it Les facteurs $\epsilon'$}\\
Si on suit la d\'emonstration depuis le d\'ebut, on voit qu'on a
obtenu au passage le r\'esultat suivant~: en se pla\c{c}ant dans la
situation globale d\'ej\`a d\'efinie, si $\pi$ est une \care\ de $G$,
il existe une repr\'esentation automorphe cuspidale $\tilde{\pi}$ de
$GL_n({\mathbb A}_{\mathbb F})$ et une repr\'esentation automorphe
cuspidale $\tilde{\pi}'$ de ${\mathbb D}({\mathbb A}_{\mathbb F})^*$
telles que

- aux places o\`u ${\mathbb D}$ est scind\'ee les composantes
locales de $\tilde{\pi}$ et $\tilde{\pi}'$ sont \'egales, 

- aux places
ramifi\'ees diff\'erentes de $v_0$ les composantes locales de
$\tilde{\pi}$ sont des repr\'esentations de Steinberg et les
composantes locales de $\tilde{\pi}'$ sont des repr\'esentations
triviales, 

- la composante locale de $\tilde{\pi}$ \`a la place $v_0$
est $\pi$ et la composante locale de $\tilde{\pi}'$ \`a la place $v_0$
est ${\bf C}(\pi)$. \\
Donc, pour toute place de $\mathbb F$
diff\'erente de $v_0$, les composantes locales de $\tilde{\pi}$ et
$\tilde{\pi}'$ ont des fonctions $L$, $\epsilon$ et $\epsilon'$
\'egales, aux places o\`u $\mathbb D$ n'est pas scind\'ee parce que ces
composantes locales sont \'egales, et aux autres par les calculs de
[GJ] sur les repr\'esentations de Steinberg. Donc, si on applique
l'\'equation fonctionnelle de [GJ], th.13.8, on trouve que les facteurs
$\epsilon'$ de $\pi$ et ${\bf C}(\pi)$ sont aussi \'egaux. Pour les
autres repr\'esentations (\ecis\ mais non unitaires) le r\'esultat
s'obtient en tordant par un caract\`ere non ramifi\'e, de fa\c{c}on
\`a les rendre unitaires, et en tenant compte du fait que le facteur
$\epsilon'$ subit le m\^eme d\'ecalage par ce changement.\\ 
\ \\
{\it Surjectivit\'e}\\
Supposons qu'il existe une \care\ $\pi'$ de $G'$ de caract\`ere central
$\omega$ qui ne se trouve pas dans l'image de {\bf{C}}. Par le
th.\ref{ortogonalitate}, avec les notations de l'annexe, la
restriction du caract\`ere de 
$\pi'$ \`a $G_e'$ est dans l'espace $L^2(G_e';\omega)$ et de plus il
est orthogonal \`a toute restriction \`a $G'_e$ d'un  caract\`ere
d'une autre \care\ de caract\`ere central $\omega$ de $G'$. En
particulier, il est orthogonal aux caract\`eres des repr\'esentations
qui se trouvent dans l'image de {\bf{C}} donc en identifiant
$L^2(G_e;\omega)$ et $L^2(G_e';\omega)$ par l'isomorphisme $i$
(voir annexe) et en utilisant la correspondance, 
la restriction du caract\`ere de $\pi'$ est orthogonale \`a la
restriction de tout caract\`ere d'une \care\ de $G$ de m\^eme
caract\`ere central. D'apr\`es le
corollaire 5.2 dans [Ba2], les restrictions des caract\`eres des
\cares\ de $G$ de caract\`ere central $\omega$ forment un syst\`eme
orthonormal {\it complet} de l'espace pr\'ehilbertien $L^2(G_e;\omega)$ ; le
caract\`ere de $\pi'$ 
serait alors nul sur $G'_e$. C'est impossible car $\chi_{\pi'}$ est de
norme 1 dans $L^2(G_e';\omega)$.
\qed
\ \\ 

Nous avons prouv\'e qu'ind\'ependamment de la caract\'eristique du corps
de base, le th.\ref{ortogonalitate} implique le th.\ref{correspforte}.
\ \\

Montrons maintenant que {\it la validit\'e du r\'esultat
\ref{correspforte} pour tout corps de caract\'eristique nulle implique
le th.\ref{correspfaible} pour tout corps de caract\'eristique non 
nulle}.
\ \\
\ \\
{\bf D\'emonstration.} Supposons donc que la caract\'eristique du
corps de base $F$ est non nulle. On proc\`ede par contradiction, en
transf\'erant \`a la caract\'eristique nulle. Du coup, nous
travaillerons avec plusieurs corps de base, et les objets porteront
parfois des indices indiquant sur quel corps nous les
consid\'erons. Nous prouverons une
correspondance faible (th.\ref{correspfaible}) entre les groupes $G_F$
et $G'_F$ en appliquant le th.\ref{correspforte} \`a $G_L$ et $G'_L$
o\`u $L$ est un corps de caract\'eristique nulle bien choisi (voir la 
section 2,
notamment sous-section 2.6). 

Soit $\pi$ une \care\ de 
$G_F$ et soit $\omega$ le caract\`ere central (unitaire) de $\pi$. \\ 
 
\begin{prop}\label{noncoresp} Soit $\pi'$ une \care\ de $G'_F$ de
caract\`ere central $\omega$. Supposons qu'il existe $M\in G_F$ et
$M'\in G'_F$ elliptiques r\'eguliers tel qu'on ait $M\lra M'$ et~: 
$$\chi_{\pi}(M)\neq (-1)^{n-r}\chi_{\pi'}(M').$$
Alors il existe un entier $m\geq 1$ tel que, si $L$ est un corps local
de caract\'eristique nulle $m$-proche de $F$ et $\cf$ est la
correspondance sur $L$ entre $GL_n(L)$ et $GL_r(D_L)$, alors $\pi$ et
$\pi'$ se rel\`event et on a 
$$\cf(\zzz(\pi))\neq \zzd(\pi').$$
\end{prop}
\ \\
{\bf D\'emonstration.} \\
{\it Principe}~: On veut trouver un voisinage ouvert compact $V$ de
$M$, un voisinage 
$V'$ ouvert compact de $M'$ et un entier $m$ tel qu'on ait~: 

1) $\chi_{\pi}$ est constant (\'egal \`a $\chi_{\pi}(M)$) sur $V$ ;

2) $\chi_{\pi'}$ est constant (\'egal \`a $\chi_{\pi'}(M')$) sur $V'$ ;

3) Pour tout corps $L$ qui est $m$-proche de $F$, $\zzz(V)$ et
   $\zzd(V')$ sont bien d\'efinis, et on a  

\ \ \ \ \ a) $\chi_{\zzz(\pi)}$ est constant \'egal \`a
$\chi_{\pi}(M)$ sur $\zzz(V)$ et 

\ \ \ \ \ b) Pour tout \'el\'ement $g'\in \zzd(V')$ il existe un
\'el\'ement $g\in \zzz(V)$ tel qu'on ait $g\lra g'$. 

Si on suppose que $\zzz(\pi)$ et $\zzd(\pi')$ se correspondent par la
correspondance en caract\'eristique nulle, on obtient alors par 3)a)
et b) que le caract\`ere $\chi_{\zzd(\pi')}$ est constant \'egal \`a
        $(-1)^{n-r}\chi_{\pi}(M)$ 
sur $\zzd(V')$. Mais, par 2), $\chi_{\pi'}$ est
constant sur $V'$ et en utilisant la fonction caract\'eristique ${\bf
1}_{V'}$ de $V'$ qui se rel\`eve, on aurait  
$$tr\pi'({\bf 1}_{V'})=tr\zzd(\pi')(\zzd({\bf 1}_{V'}))$$  
 ou encore 
$$\chi_{\pi'}(M')vol(V')=(-1)^{n-r}\chi_{\pi}(M)vol(\zzd(V')).$$
 C'est une contradiction, car $\zzd$ pr\'eserve le volume (remarque
\ref{remarquevolume}) et
$\chi_{\pi'}(M')\neq (-1)^{n-r}\chi_{\pi}(M)$.\\ 

On va donc constuire $V$, $V'$ et $m$ qui v\'erifient ces
conditions. Quitte \`a conjuguer $M$ on peut supposer qu'il est la
matrice compagnon de son polyn\^ome caract\'eristique. Soient $m$ et
$s$ comme dans la 
proposition \ref{caractseridica}. Par cette proposition, le
caract\`ere $\chi_{\zzz(\pi)}$ de $\zzz(\pi)$ est constant \'egal \`a
$\chi_{\pi}(M)$ sur l'ensemble des \'el\'ements de $G_L$ qui ont un
polyn\^ome caract\'eristique $s$-proche de celui de $M$. Donc le
caract\`ere de $\cf(\zzz(\pi))$ est constant \'egal \`a
$(-1)^{n-r}\chi_{\pi}(M)$ sur l'ensemble des \'el\'ements de $G'_L$
qui ont un polyn\^ome caract\'eristique $s$-proche de celui de
$M$. Mais le polyn\^ome caract\'eristique de $M$ est aussi celui de
$M'$. Posons $s=k$ dans la proposition \ref{34} et notons $m'$
l'entier (qui dans la proposition \ref{34} est not\'e $m$) associ\'e
\`a ce $k$. En posant $m''=\max\{m;m'\}$, si $L$ et $F$ sont
$m''$-proches, alors on a
$\bar{\zeta}^{m''}_{D_FD_L}(K_{D_F}^{m''}M'K_{D_F}^{m''})\subset
X_{M,s}$ o\`u $X_{M,s}$ est l'ensemble des \'el\'ements de $G'_L$ qui
ont un polyn\^ome caract\'eristique $s$-proche de celui de $M$. Quitte
\`a changer de voisinage (et \`a augmenter $m''$) on peut supposer
que $\chi_{\pi'}$ est constant sur $K_{D_F}^{m''}M'K_{D_F}^{m''}$
(car $M'$ est semisimple r\'egulier). On pose $V=K_F^mMK_F^m$,
$V'=K_{D_F}^{m''}M'K_{D_F}^{m''}$ et alors le triplet $(V;V';m'')$
v\'erifie les conditions 1), 2) et 3). On peut donc conclure comme
nous l'avons expliqu\'e plus haut.
\qed
\ \\

D\'emontrons maintenant le th\'eor\`eme. D\'efinissons d'abord ${\bf
C'}(\pi)$. On note $X_{G_F}^{\pi}$
(resp. $X_{G'_F}^{\pi}$) l'ensemble des classes d'\'equivalence de \cares\
de $G_F$ (resp. $G'_F$) de caract\`ere central $\omega$ et qui ont un
facteur $\epsilon'$ \'egal \`a celui de $\pi$ (calcul\'e pour le
caract\`ere additif fix\'e en d\'ebut de section). 
\begin{lemme} Les ensembles $X_{G_F}^{\pi}$
et $X_{G'_F}^{\pi}$ sont finis.
\end{lemme} 
\ \\
{\bf D\'emonstration.} Nous savons que le facteur $\epsilon'$ de $\pi$
s'\'ecrit 
$$\epsilon'(s;\pi;\psi_F)=\frac{Q(q^{-s})}{R(q^{-s})},$$
o\`u $Q$ et $R$ sont des polyn\^omes et $q$ est le cardinal du corps
r\'esiduel de $F$. On fixe une telle \'ecriture. 
Soit $\tau$ un \'el\'ement de $X_{G_F}^{\pi}$. 
On rappelle la formule \ref{defepsilon}~:
$$\epsilon'(s;\tau;\psi_F)=
\epsilon(s;\tau;\psi_F)L(1-s;\check{\tau})L(s;\tau)^{-1}.$$ 
On sait que 
$$L(s;\tau)=\frac{1}{S(q^{-s})}$$
et
$$L(1-s;\check{\tau})=\frac{1}{T(q^{s-1})}=q^{-s\, deg\, T}\frac{1}{T'(q^{-s})}$$
o\`u $S$, $T$ et $T'$ sont des polyn\^omes, $deg\, T=deg\, T'$. Alors,
si on pose  
$$\epsilon(s;\tau;\psi_F)=Kq^{-m(\tau)s},$$
o\`u $K$ est une constante complexe (non nulle), on a que
$$Kq^{-m(\tau)s}q^{-s\, deg\, T}S(q^{-s})R(q^{-s})=Q(q^{-s})T'(q^{-s}).$$
La comparaison des degr\'es des polyn\^omes qui interviennent montre
que le conducteur $m(\tau)$ de toute classe $\tau\in 
X_{G_F}^{\pi}$ est born\'e par $deg\, Q$. Le m\^eme raisonnement 
marche aussi pour $X_{G'_F}^{\pi}$, bien entendu. Alors le 
th.\ref{niveau}.b) implique que le niveau des \'el\'ements de
$X_{G_F}^{\pi}$ et $X_{G'_F}^{\pi}$ est uniform\'ement born\'e, et
comme ce sont des classes de \cares\ de caract\`ere central fix\'e, il
n'y en a qu'un nombre fini. \qed
\ \\

Supposons par l'absurde qu'il n'existe pas de
\care\ $\pi'$ de $G'_F$ de
caract\`ere central $\omega$ telle que
$\chi_{\pi}(g)=(-1)^{n-r}\chi_{\pi'}(g')$ pour tous $g\lra g'$
elliptiques r\'eguliers. Alors, en particulier, il n'en existe pas
dans $X_{G'_F}^{\pi}$. Comme cet ensemble est fini, on peut
trouver un $m$ suffisamment grand pour que  prop.\ref{noncoresp}
s'applique \`a tous $\pi'\in X_{G'_F}^{\pi}$. Soit alors $L$ un corps
$m$-proche de 
$F$. Soit $\psi_L$ un caract\`ere additif de $L$ $m$-proche de
$\psi$. Notons $X_{G_L}^{\pi}$
(resp. $X_{G'_L}^{\pi}$) l'ensemble des classes d'\'equivalence de \cares\
de $G_L$ (resp. $G'_L$) de caract\`ere central $\zzz(\omega)$ et qui ont un
facteur $\epsilon'$ \'egal \`a celui de $\pi$. Par le
th.\ref{careseridica}.b) et th.\ref{epsilonseridica}.a), $\zzd$
r\'ealise une bijection entre  $X_{G'_F}^{\pi}$ et
$X_{G'_L}^{\pi}$. Donc, en appliquant la prop.\ref{noncoresp}, on
obtient $\cf(\zzz(\pi))  \notin X_{G'_L}^{\pi}$. Ceci contredit le
fait d\'ej\`a prouv\'e que la correspondance en caract\'eristique
nulle conserve le 
facteur $\epsilon'$. Ainsi, il existe une \care\ $\pi'$ de $G_F$ telle
que  $\chi_{\pi}(g)=(-1)^{n-r}\chi_{\pi'}(g')$ pour tous $g\lra g'$
elliptiques r\'eguliers. \`A ce stade faisons la remarque que 

- $\zzd$
r\'ealise une bijection de $X_{G'_F}^{\pi}$ sur $X_{G'_L}^{\pi}$
(par th.\ref{careseridica}.b) et th.\ref{epsilonseridica}.a), on l'a dit
plus haut),

- $\zzz$
r\'ealise une bijection de $X_{G_F}^{\pi}$ sur $X_{G_L}^{\pi}$ (par le
m\^eme argument),

- $\cf$
r\'ealise une bijection de $X_{G_L}^{\pi}$ sur $X_{G'_L}^{\pi}$ par le
th.\ref{correspforte} en caract\'eristique nulle.\\
Donc, le cardinal de $X_{G_F}^{\pi}$ est \'egal \`a celui de
$X_{G'_F}^{\pi}$. D'autre 
part, en appliquant le m\^eme raisonnement que plus haut aux autres
\'el\'ements de $X_{G_F}^{\pi}$ on peut trouver pour chacun (au moins)
un
correspondant dans $X_{G'_F}^{\pi}$. Mais deux \'el\'ements de
$X_{G_F}^{\pi}$ qui ont des caract\`eres \'egaux sur les \'el\'ements
elliptiques r\'eguliers sont \'egaux, par orthogonalit\'e des
caract\`eres sur $G_F$ (vrai en caract\'eristique non nulle,
[Ba2]). Ceci, rajout\'e \`a l'\'egalit\'e des cardinaux, 
prouve que pour tout \'el\'ement de $X_{G_F}^{\pi}$ on peut
trouver exactement un \'el\'ement de $X_{G'_F}^{\pi}$ qui lui
correspond et r\'eciproquement. Ainsi, on peut d\'efinir une
correspondance pour les \cares. On \'etend
correspondance aux \ecis\ non unitaires en tordant par des
caract\`eres, comme en caract\'eristique nulle. Le fait que pour tout
$\pi$ comme plus haut la 
correspondance r\'ealise une bijection de $X_{G_F}^{\pi}$ sur
$X_{G'_F}^{\pi}$ implique qu'elle est injective et qu'elle
conserve le facteur $\epsilon'$. Montrons qu'elle est surjective. Cela
revient \`a montrer que les
ensembles du type $X_{G'_F}^{\pi}$ r\'ealisent une partition de
$E^2(G'_F)$, ou encore que, si $\pi'\in
E^2(G'_F)$, on peut toujours trouver  $\pi\in E^2(G_F)$ avec le m\^eme
facteur $\epsilon'$ que $\pi'$. Or, \`a facteur $\epsilon'$ fix\'e, 
on a vu que le niveau est born\'e par une constante $K$ (que ce soit
sur $G_F$ ou $G'_F$). 
Il suffit alors de partir avec $\pi'$ et de faire
le tour \`a l'envers, 
$G'_F-G'_L-G_L-G_F$, pour un corps 
$L$ suffisamment proche pour que toutes les repr\'esentations de niveau
inf\'erieur ou \'egal \`a $K$ se rel\`event ; nous savons que le
facteur $\epsilon'$ est conserv\'e sur le parcours.
Prouvons maintenant l'unicit\'e d'une telle correspondance. Comme
 la correspondance qu'on a d\'efini est {\it surjective}, s'il existait
une autre correspondance, alors deux classes
distinctes de \ecis\ de $G_F$ auraient des caract\`eres \'egaux
sur l'ensemble des \'el\'ements elliptiques r\'eguliers. Or, ceci
contredit l'orthogonalit\'e des caract\`eres sur $G_F$.       
\qed
\ \\

Maintenant on veut prouver que {\it le th.\ref{correspfaible} en
caract\'eristique 
non nulle implique le th.\ref{ortogonalitate} en caract\'eristique
non nulle}. Ceci est \'evident en transf\'erant les relations
d'orthogonalit\'e de $G$ \`a $G'$ par l'isomorphisme $i$ (voir
l'annexe) et la correspondance ${\bf C}'$.

R\'ecapitulons~: nous avons montr\'e que le th.\ref{ortogonalitate} en
caract\'eristique nulle (resp. non nulle)
implique le th.\ref{correspforte} en
caract\'eristique nulle (resp. non nulle). Nous avons montr\'e que le
th.\ref{correspforte} en 
caract\'eristique nulle implique le th.\ref{correspfaible} en
caract\'eristique  non
nulle et aussi que le th.\ref{correspfaible} en
caract\'eristique  non
nulle implique le th.\ref{ortogonalitate} en
caract\'eristique non nulle. Comme le th.\ref{ortogonalitate} est vrai
en
caract\'eristique nulle ([Cl]), tous ces th\'eor\`emes sont prouv\'es
en toute caract\'eristique.

\begin{cor}\label{orthogcarG'p} La restriction des caract\`eres des
classes de \cares\ de $G'$ de caract\`ere central $\omega$ fix\'e \`a
l'ensemble des \'el\'ements elliptiques r\'eguliers forme un syst\`eme
orthonormal {\underline{complet}} pour l'espace de Hilbert
$L^2(G'_e;\omega)$.  
\end{cor}
\ \\
{\bf D\'emonstration.} On avait obtenu ce r\'esultat pour $GL_n(F)$
(corollaire 5.2, [Ba2]) par transfert \`a partir d'un groupe compact
modulo le centre. Gr\^ace \`a l'isomorphisme $i$ entre
$L^2(G_e;\omega)$ et $L^2(G'_e;\omega)$ d\'efini dans l'annexe,
on le transf\`ere sans probl\`eme sur $G'$ par la correspondance
(th.\ref{correspforte}).
\newpage
\begin{section}{Annexe}      
\ \\

Soit $G$ un groupe r\'eductif connexe sur un corps local non
archim\'edien $F$ de caract\'eristique quelconque. Soient $Z$ le
centre de $G$ et $dz$ une mesure de Haar sur 
$Z$. On note $G_e$ l'ensemble des \'el\'ements elliptiques r\'eguliers
de $G$ (i.e. dont le polyn\^ome caract\'eristique est
irr\'eductible s\'eparable)~; sur tout tore elliptique maximal $T$ de $G$ on
consid\`ere  
une mesure de Haar $dt$ telle que le volume de $T/Z$ pour la mesure
quotient $d\overline{t}=dt/dz$ soit \'egal \`a 1 ; pour tout tore
maximal $T$ de $G$ on note $T^{reg}$ l'ensemble des \'el\'ements
r\'eguliers de $T$, $W_T$ le groupe de Weyl de $T$ et $|W_T|$ le
cardinal de ce groupe. Pour tout \'el\'ement semisimple r\'egulier $g$
de $G$ de centralisateur $T$, on note $D(g)$ la valeur absolue
normalis\'ee du d\'eterminant de l'op\'erateur  $Ad(g^{-1})-Id$
agissant sur $Lie(G)/Lie(T)$. Pour 
tout caract\`ere unitaire $\omega$ de $Z$ on note $L_0(G_e;\omega)$
l'espace des fonctions $f$ d\'efinies sur $G_e$ \`a 
valeurs dans $\Bbb{C}$ qui sont localement constantes, invariantes par
conjugaison par des \'el\'ements de $G$, et v\'erifient
$f(zg)=\omega(z)f(g)$ pour tout $g\in G_e$ et tout $z\in Z$. Soit ${\mathcal
T}_e$ un ensemble de repr\'esentants des classes de conjugaison de
tores elliptiques maximaux. On note  $L^2(G_e;\omega)$ le sous-espace
de $L_0(G_e;\omega)$ form\'e des fonctions $f$ pour lesquelles 
$$\sum_{T\in {\mathcal T}_e} |W_T|^{-1} \int_{T^{reg}/Z}
D(\bar{t})|f(\bar{t})|^2 d\overline{t}$$ 
 converge. On d\'efinit un produit scalaire dans $L^2(G_e;\omega)$ en
posant~: 
             $$<f_1;f_2>_e=\sum_T |W_T|^{-1}\int_{T^{reg}/Z}
D(\bar{t})f_1(\bar{t})\overline{f_2(\bar{t})} d\overline{t},$$ 
qui munit $L^2(G_e;\omega)$ une structure d'espace
pr\'ehilbertien. 

Clozel a montr\'e dans [Cl] que, si la caract\'eristique de $F$ est
nulle, alors, pour toute repr\'esentation $\pi$ de $G$ de carr\'e
int\'egrable et de caract\`ere central $\omega$, la restriction du
caract\`ere de $\pi$ \`a $G_e$ se trouve dans $L^2(G_e;\omega)$ et
les \'el\'ements de $L^2(G_e;\omega)$ ainsi obtenus forment une
famille orthonormale pour $<\ ; >_e$. Cette propri\'et\'e qu'ont tous
les groupes r\'eductifs en caract\'eristique nulle est appel\'ee
usuellement l'{\it orthogonalit\'e des caract\`eres}. Dans [Ba2] nous
avons prouv\'e cette propri\'et\'e dans le cas de $GL_n$ en
caract\'eristique non nulle, et une cons\'equence du pr\'esent article
est la validit\'e de l'orthogonalit\'e des caract\`eres pour toutes les
formes int\'erieures de $GL_n$  en 
caract\'eristique non nulle.

Dans le cas particulier o\`u $G/Z$ est compact, alors ind\'ependamment
de la caract\'eristique de $F$, les caract\`eres de toutes les
repr\'esentations lisses irr\'eductibles de $G$ de caract\`ere central
$\omega$ se trouvent dans $L^2(G_e;\omega)$ et forment une base
Hilbertienne de cet espace (cela se montre comme pour les groupes
compacts). C'est le cas du groupe des \'el\'ements inversibles d'une
alg\`ebre \`a division sur  $F$.\\

Si $G$ et $G'$ sont les groupes des \'el\'ements inversibles de deux
alg\`ebres centrales simples de m\^eme dimension $n^2$ sur $F$, on
peut identifier leurs centres qu'on note avec la m\^eme lettre $Z$. On
fixe une fois pour toutes une mesure sur $Z$. L'application qui \`a un
\'el\'ement de $G$ associe son polyn\^ome caract\'eristique sur $F$
r\'ealise une bijection entre les classes de conjugaison
d'\'el\'ements elliptiques r\'eguliers de $G$ et l'ensemble des
polyn\^omes unitaires irr\'eductibles s\'eparables de degr\'e $n$ \`a
coefficients 
dans $F$. De m\^eme pour $G'$. On obtient ainsi une bijection de
l'ensemble des classes de conjugaison d'\'el\'ements elliptiques
r\'eguliers de $G$ sur l'ensemble des classes de conjugaison
d'\'el\'ements elliptiques r\'eguliers de $G'$. On \'ecrit $g\lra g'$
si $g\in G$ est elliptique r\'egulier, $g'\in G'$ est elliptique
r\'egulier et $g$ et $g'$ ont le m\^eme polyn\^ome
caract\'eristique. Il y a un isomorphisme d'espaces vectoriels~: 
$$i:L_0(G_e;\omega)\simeq L_0(G'_e;\omega)$$
d\'efini par 
$$i(f)(g')=f(g)\ \ \ \ \ \ \ \ \text{si}\ g\lra   g'.$$
Maintenant, si $g\lra g'$, alors il y a un unique isomorphisme de $F$-alg\`ebres (qui dans ce
cas pr\'ecis sont des corps) $F[g]\simeq F[g']$, qui envoie $g$ sur
$g'$. On a donc un isomorphisme de groupes multiplicatifs \
$F[g]^*\simeq F[g']^*$. Mais $F[g]^*$ n'est autre que le tore maximal
elliptique $Z(g)$ qui contient $g$ et $ F[g']^*$ n'est autre que le
tore maximal elliptique $Z(g')$ qui contient $g'$. Toutes ces
consid\'erations impliquent que, si on choisit un ensemble de
repr\'esentants ${\mathcal T}_e$ des classes de conjugaison de tores
elliptiques maximaux de $G$ et un ensemble de repr\'esentants ${\mathcal
T'}_e$ des classes de conjugaison de tores elliptiques maximaux de
$G'$, il y a une bijection  
$$j:{\mathcal T}_e\to {\mathcal T'}_e$$
qui est caract\'eris\'ee par  $T\simeq j(T).$

Soit $T\in {\mathcal T}_e$. Notons $j_{T}$ l'isomorphisme de $T$ sur
$j(T)$. \'Etant un morphisme de groupes, $j_T$ transforme la mesure de
T en une mesure de Haar de $j(T)$, donc en un multiple de la
mesure qu'on avait choisi sur $j(T)$. Mais la restriction de
l'application $j_T$ \`a $Z$ est l'identit\'e, donc, vu le choix des
mesures sur $T$ et $i(T)$, on en d\'eduit que $j_T$ pr\'eserve la
mesure. Alors la restriction de l'isomorphisme $i$ \`a
$L^2(G_e;\omega)$ induit un isomorphisme d'espaces pr\'ehilbertiens de
$L^2(G_e;\omega)$ sur $L^2(G'_e;\omega)$. D'une fa\c{c}on plus
explicite, on a~: 
$$<i(f);i(h)>=<f;h>\ \ \ \ \ \ \ \ \forall f,h\in L^2(G_e;\omega).$$
\ \\

Si maintenant $G'$ est le groupe des \'el\'ements inversibles d'une
alg\`ebre centrale simple de degr\'e $n^2$ sur $F$, si $g'$ est un
\'el\'ement semisimple r\'egulier {\it quelconque} de $G'$ (i.e. dont
le polyn\^ome caract\'eristique $P_{g'}$ est s\'eparable), alors 
l'ensemble des \'el\'ements de $G'$ dont le polyn\^ome
caract\'eristique est $P_{g'}$ est la classe de conjugaison de
$g'$. L'ensemble des \'el\'ements de $GL_n(F)$ dont le polyn\^ome
caract\'eristique est $P_{g'}$ est non vide, constitu\'e
d'\'el\'ements semisimples r\'eguliers,  et forme une classe de
conjugaison dans $GL_n(F)$. On obtient ainsi une injection de
l'ensemble des classes de conjugaison d'\'el\'ements semisimples
r\'eguliers de $G'$ dans  l'ensemble des classes de conjugaison
d'\'el\'ements semisimples r\'eguliers de $GL_n(F)$. On note $g\lra
g'$ si $g\in G$ est semisimple r\'egulier, $g'\in G'$ est semisimple
r\'egulier et $g$ et $g'$ ont le m\^eme polyn\^ome
caract\'eristique. Soit $T'$  un tore maximal {\it quelconque} de $G'$
et $g'\in T^{reg}$. Si $g\lra g'$, alors il y a un unique isomorphisme
de $F$-alg\`ebres $F[g]\simeq F[g']$ qui envoie $g$ sur $g'$. On a
donc un isomorphisme de groupes multiplicatifs \ $F[g]^*\simeq
F[g']^*$. Mais $F[g']^*$ n'est autre que le tore $T'$ et $F[g]^*$ est
un tore maximal de $GL_n(F)$, le tore maximal qui contient $g$. On
obtient finalement une injection de l'ensemble des classes de
conjugaison de tores maximaux de $G'$ dans l'ensemble des classes de
conjugaison de tores maximaux de $GL_n(F)$ qui prolonge la bijection
$j$ d\'efinie plus haut entre l'ensemble des classes de conjugaison de
tores maximaux elliptiques de $G'$ et  l'ensemble des classes de
conjugaison de tores maximaux elliptiques de $GL_n(F)$. On note cette
application toujours par la lettre $j$. Le raisonnement qu'on vient de
faire implique aussi que pour tout tore $T'$ de $G'$ il existe un
isomorphisme de $T'$ sur tout \'el\'ement $T$ de la classe de
conjugaison correspondant \`a la classe de conjugaison de $T'$. Si on
fixe maintenant des mesures de Haar sur tous les tores maximaux 
de $GL_n(F)$, avec la propri\'et\'e que sur deux tels tores
conjugu\'es les mesures se correspondent via la conjugaison, on peut
obtenir par les isomorphismes plus haut des mesures de Haar sur tous
les tores maximaux non elliptiques de toutes les formes int\'erieures
de $GL_n(F)$.  
\end{section}
\newpage

\section{Bibliographie}

{\begin {itemize}

\item[{[Ba1]}] A.I.Badulescu, {\it Correspondance entre $GL_n$ et ses
formes int\'erieures en caract\'eristique positive}, th\`ese,
Univ.Paris XI Orsay, 1999.\\

\item[{[Ba2]}] A.I.Badulescu, Orthogonalit\'e des caract\`eres pour
$GL_n$ sur un corps local de caract\'eristique non nulle, {\it
Manuscripa Math.} 101 (2000), 49-70.\\

\item[{[Ba3]}] A.I.Badulescu, Un th\'eor\`eme de finitude, {\it
pr\'epublication}.\\ 

\item[{[BC]}] Z.I.Borevitch, I.R.Chafarevitch, {\it
Th\'eorie des nombres}, Monographies internationales de
math\'ematiques modernes, Gauthier-Villars Paris, 1967.\\ 

\item[{[Be]}] J.Bernstein, Le ``centre'' de Bernstein, r\'edig\'e par
Deligne, {\it Repr\'esentations des groupes r\'eductifs sur un corps
local}, Hermann, Paris 1984.\\ 

\item[{[Ca1]}] W.Casselman, Characters and Jacquet modules, {\it
Math. Ann.} 230, (1977), 101-105.\\ 

\item[{[Ca2]}] W.Casselman, Introduction to the theory of admissible
representations of reductive $p$-adic groups, pr\'epublication.\\ 

\item[{[Cl]}] L.Clozel, Invariant harmonic analysis on the Schwarz
space of a reductive $p$-adic group, {\it Proc.Bowdoin Conf.1989,
Progress in Math.Vol.101}, Birkh\"auser, Boston, 1991, 101-102.\\ 

\item[{[De]}] P.Deligne, Les corps locaux de caract\'eristique $p$,
limites de corps locaux de caract\'eristique 0, {\it Repr\'esentations
des groupes r\'eductifs sur un corps local}, Hermann, Paris 1984.\\  

\item[{[DKV]}] P.Deligne, D.Kazhdan, M.-F.Vign\'eras,
Repr\'esentations des alg\`ebres centrales simples $p$-adiques, {\it
Repr\'esentations des groupes r\'eductifs sur un corps local},
Hermann, Paris 1984.\\ 

\item[{[Fla]}] D.Flath, A comparison of the automorphic
representations of $GL(3)$ and its twisted forms, {\it Pacific J. of
Math.} 97 (1981), 373-402.\\ 

\item[{[GJ]}] R.Godement, H.Jacquet, {\it Zeta functions of simple
algebras}, S.L.N. 260 (1972).\\ 

\item[{[H-C]}]  Harish-Chandra, A submersion principle and its applications, {\it Proc.Indian Acad.Sc.} 90 (1981), 95-102. \\

\item[{[He]}] G.Henniart, La conjecture de Langlands locale pour
$GL(3)$, {\it M\'em. de la S.M.F. (nouvelle s\'erie)} 11/12, 1984.\\ 

\item[{[Ho]}] R.Howe, Harish-Chandra homomorphism for $p$-adic groups,
{\it Regional Conferences Series in Math.}, 59(1985), Amer.Math.Soc,
Providence, R.I.\\ 
 
\item[{[JL]}] H.Jacquet, R.P.Langlands, {\it Automorphic forms on
GL(2)}, L.N.M. 114, Springer-Verlag 1970.\\ 

\item[{[Ka]}] D.Kazhdan, Representations of groups over close local
fields, {\it J. Analyse Math.} 47 (1986), 175-179.\\ 

\item[{[Le]}] B.Lemaire, th\`ese, Univ. Paris Sud, 1994.\\

\item[{[Pi]}] R.S.Pierce, {\it Associative algebras}, Grad. Texts in
Math. 88, Springer-Verlag.\\  
 
\item[{[Ro]}] J.Rogawski, Representations of $GL(n)$ and division
algebras over a $p$-adic field, {\it Duke Math. J.} 50 (1983),
161-201.\\ 

\item[{[Sa]}] I.Satake, Theory of Spherical Functions on Reductive
Algebraic Groups over $p$-adic Fields, {\it Publ.Math.I.H.E.S.}18,
(1963), 1-69.\\ 

\item[{[Sh]}] J.A.Shalika, The multiplicity one theorem for $GL_n$,
{\it Ann.of Math.} 100 (1974), 171-193.\\

\item[{[We]}] A.Weyl, {\it Basic Number Theory}, Classics in Math.,
Springer-Verlag 1973.\\ 

\item[{[Ze]}]  A.Zelevinski, Induced representations of reductive
$p$-adic groups II, {\it Ann. Sci. ENS} 13 (1980), 165-210.

\end{itemize}}

\end{document}